\documentclass{amsart}
\usepackage{amssymb}   % For Latex2e
\usepackage{amsmath}
\usepackage{amsthm}
\def\CC{{\mathbb C}}

\def\GG{{\mathbb G}}
\def\HH{{\mathbb H}}

\def\PP{{\mathbb P}}

\def\RR{{\mathbb R}}

\def\ZZ{{\mathbb Z}}
\def\UU{{\mathbb U}}

\def\Acal{{\mathcal A}}
\def\Bcal{{\mathcal B}}
\def\Ccal{{\mathcal C}}

\def\Fcal{{\mathcal F}}

\def\Hcal{{\mathcal H}}
\def\Jcal{{\mathcal J}}

\def\Lcal{{\mathcal L}}
\def\Mcal{{\mathcal M}}
\def\Ncal{{\mathcal N}}
\def\Ocal{{\mathcal O}}
\def\Pcal{{\mathcal P}}

\def\Scal{{\mathcal S}}
\def\Tcal{{\mathcal T}}
\def\Qcal{{\mathcal Q}}

\def\Xcal{{\mathcal X}}

\def\ST{{\rm Sing}(\Theta_X)}
\def\Xg{{\mathcal X}_g}
\def\Ag{{\mathcal A}_g}
\def\t{\tau}
\def\th{\vartheta}

\def\dd{\partial}

\newcommand{\map}{\dasharrow}

\newcommand\proofsquare{\nobreak\hfill \hbox{%
\vrule height 5pt
\kern-.4pt
 \vbox{%
\hrule width 5pt depth0pt height.4pt
 \kern4.6pt \hrule  }
\kern-3.75pt
\vrule height 5pt}\kern1pt
\par}

\newtheorem{theorem}{Theorem}[section]
\newtheorem{lemma}[theorem]{Lemma}
\newtheorem{proposition}[theorem]{Proposition}
\newtheorem{corollary}[theorem]{Corollary}
\newtheorem{conjecture}[theorem]{Conjecture}
\newtheorem{definition-lemma}[theorem]{Definition-Lemma}
\newtheorem{claim}[theorem]{Claim}

\theoremstyle{definition}
\newtheorem{definition}[theorem]{\bf Definition}
\newtheorem{example}[theorem]{\bf Example}

\theoremstyle{remark}
\newtheorem{remark}[theorem]{\bf Remark}

\begin{document}

\title[Andreotti--Mayer Loci and the Schottky Problem]
{Andreotti--Mayer Loci and the Schottky Problem 
%\footnote{\tt amlschottky10122006, \today}
}
\author{Ciro Ciliberto}
\address{Dipartimento di Matematica, II Universit\`a di Roma,
Italy}
\email{cilibert@axp.mat.uniroma2.it}
\author{Gerard van der Geer}
\address{Korteweg-de Vries Instituut, Universiteit van
Amsterdam, Plantage Muidergracht~24, 1018 TV Amsterdam, The Netherlands}
\email{geer@science.uva.nl}

\begin{abstract}
We prove a lower bound for the codimension of the Andreotti-Mayer locus
$N_{g,1}$ and show that the lower bound is reached only for
the hyperelliptic locus in genus $4$ and the Jacobian locus in genus $5$.
In relation with the intersection of the Andreotti-Mayer loci
with the boundary of the moduli space ${\Acal}_g$
we study subvarieties of principally
polarized abelian varieties $(B,\Xi)$ parametrizing points $b$
such that $\Xi$ and the translate $\Xi_b$ are tangentially degenerate
along a variety of a given dimension.
\end{abstract}

\maketitle

\begin{section}{Introduction}\label{sec:intro}
The Schottky problem asks for a characterization of Jacobian varieties
among all principally polarized abelian varieties. In other words, it
asks for a description of the Jacobian locus ${\mathcal J}_g$ in the
moduli space $\Ag$ of all principally polarized abelian varieties
of given dimension $g$. In the 1960's Andreotti and Mayer (see \cite{AM})
pioneered an approach based on the fact that the Jacobian variety of a 
non-hyperelliptic (resp.\ hyperelliptic) curve of genus $g\geq 3$ has a 
singular locus of dimension $g-4$ (resp.\ $g-3$). 
They introduced the loci $N_{g,k}$ of principally polarized
abelian varieties $(X,\Theta_X)$ of dimension $g$ with a singular 
locus of $\Theta_X$
of dimension $\geq k$ and showed that ${\mathcal J}_g$ 
(resp.\ the hyperelliptic locus ${\mathcal H}_g$) 
is an irreducible component of $N_{g,g-4}$ (resp.\ $N_{g,g-3}$).
However, in general there are more irreducible components of $N_{g,g-4}$ so that the
dimension of the singular locus of $\Theta_X$ does not suffice to characterize
Jacobians or hyperelliptic Jacobians. The locus $N_{g,0}$ of abelian varieties
with a singular theta divisor has codimension $1$ in $\Ag$ and in a beautiful paper
(see \cite{Mumford}) Mumford calculated its class. But in general not much is known about these
Andreotti-Mayer loci $N_{g,k}$. In particular, we do not even know their
codimension. In this paper we give estimates for the codimension of these
loci. These estimates are in general not sharp, but  
we think that the following conjecture gives the sharp bound.

\begin{conjecture}\label{conj:i} If $1 \leq k \leq g-3$ and if $N$ is
an irreducible component of $N_{g,k}$ whose general point corresponds
to an abelian variety with endomorphism ring $\ZZ$ then ${\rm codim}_{\Ag}(N)
\geq \binom{k+2}{2}$. Moreover, equality holds if and only if one of the
following happens:
\begin{itemize}
\item [(i)]  $g=k+3$ and $N={\Hcal}_g$;
\item [(ii)] $g=k+4$ and  $N=\Jcal_g$.
\end{itemize}\end{conjecture}

We give some evidence for this conjecture by proving the case $k=1$.
In our approach we need to study the behaviour of the Andreotti-Mayer loci
at the boundary of the compactified moduli space. A principally
polarized $(g-1)$-dimensional abelian variety $(B,\Xi)$  parametrizes 
semi-abelian varieties that are extensions of $B$ by the 
multiplicative group $\GG_m$. This means that $B$ occurs in the 
boundary of the compactified moduli space $\tilde{\mathcal A}_g$
and we can intersect $B$ with the Andreotti-Mayer loci. 
This motivates the definition of loci $N_k(B,\Xi)\subset B$
for a principally polarized $(g-1)$-dimensional abelian variety $(B,\Xi)$. 
They are formed by the points $b$ in $B$ such that $\Xi$ and its 
translate $\Xi_b$ are tangentially degenerate along a subvariety of 
dimension~$k$.
These intrinsically defined subvarieties of an abelian variety are interesting
in their own right and deserve further study.
The conjecture above then leads to a boundary version 
that gives a new conjectural answer to the Schottky problem
for simple abelian varieties.

\begin{conjecture}
Let $k \in \ZZ_{\geq 1}$.
Suppose that $(B,\Xi)$ is a simple principally polarized abelian
variety of dimension $g$ not contained in $N_{g,i}$ for
all $i\geq k$. Then there is an irreducible component $Z$ of
$N_k(B,\Xi)$ with  ${\rm codim}_B(Z)= k+1$ if and only if one of
the following happens:
\begin{itemize}
\item [(i)] either $g\geq 2$, $k=g-2$ and $B$ is
a hyperelliptic Jacobian,
\item  [(ii)] or $g\geq 3$, $k=g-3$ and $B$ is a Jacobian.
\end{itemize}
\end{conjecture}

In our approach  we will use a special compactification
$\tilde{\mathcal A}_g$ of $\Ag$  (see \cite{Nam, Nakam, Br}).  
The points of the boundary $\partial\tilde{\mathcal A}_g=
\tilde{\mathcal A}_g-\Ag$ correspond to suitable compactifications 
of $g$--dimensional semi-abelian varieties. 
We prove Conjecture \ref{conj:i} for $k=1$
by intersecting with the boundary. For higher values of $k$, the intersection
with the boundary looks very complicated. \end{section}

\begin{section}{The universal theta divisor}\label{sec: univtheta}

Let $\pi : \Xg \to \Ag$ be the universal principally polarized abelian
variety of relative dimension $g$ over the moduli space $\Ag$ of principally 
polarized abelian varieties of dimension $g$ over $\CC$. In this paper 
we will work with orbifolds and we shall identify $\Xg$ (resp.\
$\Ag$) with the orbifold
${\rm Sp}(2g,\ZZ) \ltimes \ZZ^{2g} \backslash \HH_g \times \CC^g$ 
(resp.\ with ${\rm Sp}(2g,\ZZ)\backslash \HH_g$), where
$$
\HH_g=\{ (\tau_{ij})\in {\rm Mat}(g \times g, \CC)
: \tau=\tau^t, {\rm Im}(\tau)> 0\}
$$
is the usual \emph {Siegel  upper--half space} of degree $g$. 
The $\tau_{ij}$ with $1\leq i\leq j\leq g$ are coordinates on 
$\HH_g$ and we let $z_1,...,z_g$ be coordinates on $\CC^ g$. 

The {\it Riemann theta function} 
$\vartheta(\tau,z)$, given on $\HH_g \times \CC^g$ by 
$$
\vartheta(\tau,z)=\sum_{m \in \ZZ^g} e^{\pi i [m^t \tau m +2m^t z]},
$$
is a holomorphic function and its zero locus is an effective divisor 
$\tilde{\Theta}$ on  $\HH_g \times \CC^g$ which descends to a 
divisor $\Theta$ on $\Xg$. 
If the abelian variety $X$ is a fibre of $\pi$, then we let 
$\Theta_X$ be the restriction of $\Theta$ to $X$. Note that since
$\theta(\tau,z)$ satisfies $\theta(\tau,-z)=\theta(\tau,z)$, the
divisor $\Theta_X$ is symmetric, i.e.,
$\iota^*(\Theta_X)=\Theta_X$, where $\iota=-1_X :X\to X$ is
multiplication by $-1$ on $X$. The divisor $\Theta_X$ defines the line bundle
$\Ocal_X(\Theta_X)$, which  yields the principal polarization on $X$.
The isomorphism class of the pair $(X,\Theta_X)$ represents a 
point $\zeta$ of $\Ag$ and  we will write
$\zeta=(X,\Theta_X)$.  Similarly, it will be convenient to identify a
point $\xi$ of $\Theta$ with the isomorphism class of a
representative  triple $(X,\Theta_X,x)$, where $\zeta=(X,\Theta_X)$
represents $\pi(\xi)\in \Ag$  and $x\in \Theta_X$. 

The tangent space to
${\mathcal X}_g$ at a point $\xi$, with $\pi(\xi)=\zeta$,  
will be identified with the tangent space
$T_{X,x} \oplus  T_{\Ag,\zeta} \cong T_{X,0}\oplus {\rm
Sym}^2(T_{X,0})$. If $\xi=(X,\Theta_X,x)$ corresponds to the  
${\rm Sp}(2g,\ZZ) \ltimes \ZZ^{2g}$--orbit of a point 
$(\t_0,z_0) \in \HH_g \times \CC^g$, then the tangent space 
$T_{\Xg,\xi}$ to $\Xg$ at $\xi$ can be identified
with the tangent space to $\HH_g \times \CC^g$
at $(\tau_0,z_0)$, which in turn is naturally isomorphic to
$\CC^{g(g+1)/2+g}$, with coordinates $(a_{ij},b_{\ell})$ for
 $1\leq i,j\leq g$ and
$1\leq \ell \leq g$ that satisfy $a_{ij}=a_{ji}$. 
We thus view the $a_{ij}$'s as coordinates on the tangent space to
$\HH_g$ at $\tau_0$ and the $b_l$'s as coordinates on the
tangent space to $X$ or its universal cover.

An important remark is that by identifying the tangent space
to $\Ag$ at  $\zeta=(X,\Theta_X)$ with ${\rm Sym}^ 2(T_{X,0})$,
we can view  the projectivized tangent space
$\PP(T_{\Ag,\zeta})\cong \PP({\rm Sym}^ 2(T_{X,0}))$ as the
linear system of all dual quadrics in $\PP^ {g-1}=\PP(T_{X,0})$. In
particular, the matrix $(a_{ij})$  can be interpreted as the matrix
defining a dual quadric in the space $\PP^{g-1}$ with homogeneous
coordinates $(b_1:\ldots : b_g)$. 
Quite naturally, we will often use
$(z_1:\ldots:z_g)$ for the homogeneous coordinates in $\PP^{g-1}$.

Recall that the Riemann theta function $\th$
satisfies the {\it heat equations}
$$
\frac{\dd}{\dd z_i} \frac{\dd}{\dd z_j} \th 
= 2\pi \sqrt{-1} (1+\delta_{ij})\frac{\dd}{\dd {\t_{ij}}} \th
$$
for $1 \leq i,j\leq g$, where $\delta_{ij}$ is the Kronecker delta. We shall
abbreviate this equation as
$$
\dd_i \dd_j \th =  2\pi \sqrt{-1} (1+\delta_{ij}) \dd_{\tau_{ij}} \th,
$$
where $\dd_j$ means the partial
derivative $\dd/\dd z_j$ and $\dd_{\tau_{ij}}$ the partial derivative
$\dd / \dd \tau_{ij}$. One easily checks that also
all derivatives of $\theta$ verify the heat equations.
We refer to \cite{We1} for an algebraic interpretation of the heat equations
in terms of deformation theory.

If $\xi=(X,\Theta_X,x)\in \Theta$ corresponds to the 
${\rm Sp}(2g,\ZZ) \ltimes \ZZ^{2g}$--orbit of a point 
$(\t_0,z_0)$, then the Zariski tangent space
$T_{\Theta,\xi}$ to $\Theta$ at $\xi$ is the subspace of $T_{\Xg,\xi}\simeq
\CC^{g(g+1)/2+g}$ defined, with the
above conventions, by the linear equation

\begin{equation}\label{eq:ts} \sum_{1\leq i \leq j\leq g}
\frac 1 {2 {\pi} \sqrt{-1}(1+\delta_{ij})} a_{ij}\, \dd_i\dd_j\th(\t_0,z_0)  +
\sum_{1\leq \ell \leq g} b_{\ell} \, \dd_{\ell} \th (\t_0, z_0)
=0\end{equation} 
in the variables $(a_{ij},b_{\ell})$, $1\leq i,j\leq g$,
$1\leq \ell \leq g$. As an immediate consequence we get the result
(see \cite{SV1}, Lemma (1.2)): 

\begin{lemma}\label{lem:triplept} 
The point $\xi=(X,\Theta_X,x)$ is a singular point of $\Theta$ 
if and only if $x$ is a point of multiplicity at least 3 for 
$\Theta_X$.
\end{lemma}
\end{section}
\begin{section}{The locus $S_g$}\label{sec: essgi}

We begin by defining a suborbifold of $\Theta$ supported on the set
of points where $\pi_{|\Theta}$ fails to be of maximal rank.

\begin{definition} The closed suborbifold $S_g$ of $\Theta$ is defined
on the universal cover $\HH_g \times \CC^g$ by the $g+1$ equations
\begin{equation}\label{eq:singt}
\th (\t,z)  =0, \quad
\dd_j\th(\t,z)  =0, \quad j=1,\ldots,g. 
\end{equation}
\end{definition}

Lemma \ref {lem:triplept} implies that the support of $S_g$ 
is the union of ${\rm Sing}(\Theta)$ and of the set of 
smooth points of $\Theta$ where $\pi_{|\Theta}$ fails to be 
of maximal rank. Set-theoretically one has
$$
S_g=\{(X,\Theta_X,x)\in \Theta : x \in {\rm Sing}(\Theta_X)\}
$$
and ${\rm codim}_{\Xg}(S_g)\leq g+1$. It turns out that every
irreducible component of $S_g$ has codimension $g+1$ in $\Xg$
(see \cite{C-vdG} and an unpublished preprint by Debarre \cite{De1}). 
We will come back to this later in \S \ref{sec:amloci} and 
\S \ref{lowerbounds}.

With the above identification, the Zariski tangent space to $S_g$
at a given point $(X,\Theta_X,x)$ of $\Xg$, corresponding to the 
${\rm Sp}(2g,\ZZ)$-orbit of a point 
$(\t_0,z_0) \in \HH_g \times \CC^g$, is given by the
$g+1$ equations 

\begin{equation}\label{eq:tansp}
\begin{aligned}
\sum_{1\leq i\leq j\leq g} a_{ij} \dd_{\t_{ij}} \th (\t_0,z_0)  &=0,\cr
\sum_{1\leq i \leq j\leq g}
a_{ij} \dd_{\t_{ij}}\dd_k \th (\t_0,z_0)  + \sum_{1\leq \ell \leq g} b_{\ell}
\dd_{\ell} \dd_k \th (\t_0, z_0) &=0, \qquad 1\leq k \leq g \cr
\end{aligned}
\end{equation}
in the variables  $(a_{ij},b_{\ell})$ with $1\leq i,j, \ell \leq g$.
We will use the following notation:
\begin{enumerate}
\item $q$ is the row vector of length $g(g+1)/2$, given by 
$(\partial_{\t_{ij}} \theta (\t_0,z_0))$, with lexicographically
ordered entries;
\item $q_k$ is the row vector of length $g(g+1)/2$, 
given by $(\dd_{\tau_{ij}} \dd_k\theta (\t_0,z_0))$,  with lexicographically
ordered entries;
\item $M$ is the $g\times g$--matrix
$(\dd_i\dd_j \th (\t_{0},z_0))_{1\leq i,j\leq g}$.
\end{enumerate} 
Then we can rewrite
the equations \eqref {eq:tansp} as 
\begin{equation}\label {eq:tansp2}
a \cdot q^t =0, \quad a\cdot q_j^t + b \cdot M^ t_j =0, \, (j=1,\ldots,g),
\end{equation}
where $a$ is the vector $(a_{ij})$ of length $g(g+1)/2$, with
lexicographically ordered entries, $b$ is a vector in $\CC^g$ and $M_j$ the
$j$--th row of the matrix $M$.

In this setting, the equation \eqref {eq:ts} for the tangent space to
$T_{\Theta,\xi}$ can be written as:
\begin{equation}\label {eq:ts2}
a \cdot q^t + b\cdot \partial \th (\t_0,z_0)^ t =0
\end{equation}
where $\partial$ denotes, as usual, the gradient. 

Suppose now the point 
$\xi=(X,\Theta_X,x)$ in $S_g$, 
corresponding to $(\t_0,z_0) \in \HH_g \times \CC^g$
is not a point of ${\rm Sing}(\Theta)$. By Lemma \ref {lem:triplept} the
matrix $M$ is not zero and therefore we can associate to $\xi$ a quadric
$Q_\xi$ in the projective space
$\PP(T_{X,x})\simeq \PP(T_{X,0})\simeq \PP^ {g-1}$, namely the one defined by
the equation
$$ 
b \cdot M \cdot b^ t=0.
$$
Recall that $b=(b_1,\ldots,b_g)$ is a coordinate vector on 
$T_{X,0}$  and therefore $(b_1:\ldots:b_g)$ are homogeneous coordinates on
$\PP(T_{X,0})$. We will say that $Q_\xi$ is {\it indeterminate}, if $\xi\in
{\rm Sing}(\Theta)$.

The vector $q$ naturally lives in ${\rm Sym}^2(T_{X,0})^{\vee}$ and
therefore, if $q$ is not zero, the point 
$[q]\in \PP({\rm Sym}^2(T_{X,0})^{\vee})$
determines a quadric in  $\PP^ {g-1}=\PP(T_{X,0})$. 
The heat equations imply that this quadric coincides with $Q_\xi$. 

Consider the matrix defining the Zariski tangent space to 
$S_g$ at a point $\xi=(X,\Theta_X,x)$. We denote by $r:=r_\xi$ 
the corank of the quadric $Q_\xi$, with the convention that 
$r_{\xi}=g$ if $\xi\in {\rm Sing}(\Theta)$, i.e., if $Q_{\xi}$ is
indeterminate. If we choose coordinates on $\CC^g$ such that the first
$r$ basis vectors generate the kernel of $q$ then
the shape of the matrix $A$ of the system \eqref {eq:tansp} is
\begin{equation}\label {eq:a}
A=\left(\begin{matrix} q && 0_g \cr
q_1 && 0_g \cr
& \vdots & \cr
q_r && 0_g \cr
* && B \cr \end{matrix} \right),
\end{equation}
where $q$ and 
$q_k$ are as above and $B$ is a $(g-r)\times g$--matrix with the first
$r$ columns equal to zero and the remaining $(g-r)\times (g-r)$ matrix
symmetric of maximal rank.

Next, we characterize the smooth points $\xi=(X,\Theta_X,x)$
of $S_g$. 
Before stating the result, we need one more piece of notation. 
Given a non-zero vector $b=(b_1,\ldots,b_g)\in T_{X,0}$, we set 
$\dd_b = \sum_{\ell=1}^ g b_\ell \dd_\ell$. Define the matrix 
$\dd_b M$ as the $g\times g$--matrix
$(\dd_i\dd_j\dd_b\th(\t_0,z_0))_{1\leq i,j\leq g}$. 
Then define the quadric $\dd_bQ_\xi=Q_{\xi,b}$ of $\PP(T_{X,0})$ 
by the equation
$$
z\cdot \dd_bM\cdot z^ t=0.
$$
If $z=e_i$ is the $i$--th vector of the standard basis, one
writes $\dd_iQ_\xi=Q_{\xi,i}$ instead of $Q_{\xi,e_i}$ for $i=1,\ldots,g$. 
We will use similar notation for higher order derivatives or even for 
differential operators applied to a quadric. 

\begin{definition}
We let $\Qcal_{\xi}$ be the linear system of quadrics in $\PP(T_{X,0})$ 
spanned by $Q_\xi$ and by all quadrics $Q_{\xi,b}$ with
$b\in \ker(Q_\xi)$. 
\end{definition}
Since $Q_\xi$ has corank $r$, the system  
$\Qcal_\xi$ is spanned by $r+1$ elements and 
therefore $\dim(\Qcal_\xi)\leq r$. 
This system may happen to be empty, but then $Q_\xi$ is indeterminate,
i.e., $\xi$ lies in ${\rm Sing}(\Theta)$.  Sometimes we will use the
lower suffix $x$ instead of $\xi$ to denote quadrics and linear systems,
e.g.\  we will sometimes write $Q_x$ instead of $Q_\xi$, etc. 
By the heat equations,
the linear system $\Qcal_\xi$ is the image of the vector subspace of 
${\rm Sym}^2(T_{X,0})^{\vee}$ spanned by  the vectors $q, q_1,\ldots, q_r$.

\begin{proposition}\label {prop:sing} The subscheme  $S_g$ is
smooth of codimension $g+1$ in $\Xg$ at the point $\xi=(X,\Theta_X,x)$
of $S_g$ if and only if the following conditions are verified:

\begin{enumerate}
\item [(i)] $\xi\notin {\rm Sing}(\Theta)$, i.e., $Q_\xi$ is not indeterminate
and of corank $r<g$;
\item [(ii)] the linear system $\Qcal_\xi$ has maximal dimension $r$; in
particular, if $b_1,\ldots,b_r$ span the kernel of $Q_\xi$, then the $r+1$
quadrics $Q_\xi$, $Q_{\xi,b_1},\ldots,Q_{\xi,b_r}$ are linearly
independent.  \end{enumerate}
\end{proposition}

\begin{proof}  The subscheme $S_g$ is smooth of codimension $g+1$ in $\Xg$ at
$\xi$ if and only if the matrix $A$ appearing in \eqref {eq:a} has maximal
rank $g+1$. Since the submatrix $B$ of $A$ has rank $g-r$, the assertion
follows.  \end{proof}

\begin{corollary}\label {cor:smoothsg}
If $Q_\xi$ is a smooth quadric, then $S_g$ is smooth at
$\xi=(X,\Theta_X,x)$. \end{corollary}
\end{section}
\begin{section}{Quadrics and Cornormal Spaces}\label{sec: conormal}

Next we study the differential of the restriction to $S_g$ 
of the map $\pi: \Xg \to \Ag$ at a point  $\xi=(X,\Theta_X,x) \in S_g$. 
We are interested in the kernel and the image of $d\pi_{|S_g,\xi}$. 
We can view these spaces in terms of the geometry of $\PP^
{g-1}=\PP(T_{X,0})$ as follows: 
$$
\Pi_\xi=\PP( {\rm ker}
(d\pi_{|S_g,\xi})) \subseteq \PP(T_{X,0})
$$
is a linear subspace of  $\PP(T_{X,0})$ and 
$$
\Sigma_{\xi}=\PP( {\rm Im} (d\pi_{|S_g,\xi})^\perp) \subseteq 
\PP({\rm Sym}^2 (T_{X,0})^{\vee})
$$ 
is a linear system of quadrics in $\PP(T_{X,0})$.

The following proposition is the key to our approach; we use it
to view the quadrics as elements of the conormal space to our
loci in the moduli space.

\begin{proposition}\label{prop:deformations} Let $\xi=(X,\Theta_X,x)$ be a
 point of $S_g$. Then:

\begin{itemize}
\item [(i)] $\Pi_\xi$ is the vertex of the quadric $Q_\xi$. In
particular, if $\xi$ is a singular point of  $\Theta$, then  $\Pi_\xi$ is the
whole space $\PP(T_{X,0})$;
\item [(ii)] $\Sigma_{\xi}$ contains the linear system $\Qcal_\xi$. 
\end{itemize}
\end{proposition}

\begin{proof} The assertions follow from the shape of the matrix
$A$ in \eqref{eq:a}. 
\end{proof}

This proposition tells us that, given a point
$\xi=(X,\Theta_X,x) \in S_g$, the map $d\pi_{|S_g,\xi}$
is not injective if and only if the quadric $Q_\xi$ is singular. 

The orbifold $S_g$ is stratified by the corank of the matrix 
$(\partial_i\partial_j \theta)$.
\begin{definition}
For $0\leq k \leq g$  we define $S_{g,k}$ as the closed suborbifold 
of $S_g$ defined by the equations on $\HH_g \times \CC^g$
\begin{equation}\label {eq:equa}
\begin{aligned}
& \th (\t,z)  =0, \quad
\dd_j\th(\t,z)  =0, \quad (j=1,\ldots,g),\cr
&{\rm rk}\big( (\dd_i\dd_j\th(\t,z))_{1\leq i,j\leq g}\big) \leq g-k.\cr
\end{aligned}
\end{equation}
\end{definition}
Geometrically this means that $\xi\in S_{g,k}$  
if and only if $\dim(\Pi_\xi)\geq k-1$
or equivalently $Q_\xi$ has corank at least $k$. 
We have the inclusions
$$
S_g=S_{g,0}\supseteq S_{g,1}\supseteq \ldots
 \supseteq S_{g,g}=S_g\cap {\rm Sing}(\Theta)
$$ 
and $S_{g,1}$ is the locus where the map $d\pi_{|S_g,\xi}$ is not
injective. The loci $S_{g,k}$ have been considered also in \cite {GS}.

We have the following dimension estimate for the $S_{g,k}$.

\begin{proposition}\label {prop:exp} 
Let $1 \leq k\leq g-1$ and let $Z$ be an
irreducible component of $S_{g,k}$ not contained in 
$S_{g,k+1}$. Then we have
$$
{\rm codim}_{S_g}(Z)\leq \binom{k+1}{2}.
$$
\end{proposition}
\begin{proof}  Locally, in a
neighborhood $U$ in $S_g$ of a point $z$ of $Z \backslash S_{g,k+1}$ 
we have a morphism $f: U\to \Qcal$, where $\Qcal$ is the linear system 
of all quadrics in $\PP^{g-1}$. 
The map $f$ sends $\xi=(X,\Theta_X,x)\in U$ to $Q_\xi$. The scheme 
$S_{g,k}$ is the pull--back of the subscheme $\Qcal_k$ of $\Qcal$ formed by
all quadrics of corank $k$. Since ${\rm codim}_\Qcal
(\Qcal_k)=\binom{k+1}{2}$, the assertion follows.
\end{proof}

Using the equations \eqref{eq:equa} it is possible to make a local
analysis of the schemes $S_{g,k}$, e.g. it is possible to write down
equations for their Zariski tangent spaces (see \S \ref {higherorder}  for
the case $k=g$). This is however not particularly illuminating, and we
will not dwell on this here.

It is useful to give an interpretation of the points
$\xi=(X,\Theta_X,x) \in S_{g,k}$ in terms of singularities of the theta
divisor $\Theta_X$. 
Suppose that $\xi$ 
is such  that ${\rm Sing} (\Theta_X)$ contains a subscheme isomorphic to
${\rm Spec}(\CC[\epsilon]/(\epsilon^2))$ supported at $x$. This subscheme
of $X$ is given by a homomorphism 
$$
\Ocal_{X,x} \to \CC[\epsilon]/(\epsilon^2),\quad
f \mapsto f(x)+\Delta^{(1)}f(x)\cdot \epsilon, 
$$  
where $\Delta^{(1)}$ is a non--zero differential operator of order $\leq
1$, hence $\Delta^{(1)}=\dd_b$, for some non--zero vector $b\in \CC^ g$.
Then the condition ${\rm Spec}(\CC[\epsilon]/(\epsilon^2)) \subset {\rm
Sing} (\Theta_X)$ is equivalent to saying that $\th$ and $\dd_b \th$
satisfy the equations
\begin{equation}\label{eq:due}
f(\t_0,z_0) =0,\quad
\dd_j f(\t_0,z_0)=0, \quad 1\leq j \leq g,
\end{equation}
and this, in turn, is equivalent to the fact that the quadric
$Q_{\xi}$ is singular at the point~$[b]$.

More generally, we have the following proposition, which explains the nature
of the points in $S_{g,k}$ for $k<g$. 

\begin{proposition}\label{kerQ} Suppose that $x\in \ST$ does not lie on ${\rm
Sing} (\Theta)$. Then $\ST$ contains a scheme isomorphic to ${\rm
Spec}(\CC[\epsilon_1,\ldots, \epsilon_k]/(\epsilon_i\epsilon_j\colon 1 \leq
i,j \leq k<g))$ supported at
$x$ if and only if the quadric $Q_\xi$ has corank $r\geq k$. 
Moreover, the Zariski tangent space to $\ST$ at $x$ is the kernel space
of $Q_{\xi}$. \end{proposition}

\begin{proof}
With a suitable choice of coordinates in $X$, the condition
that the scheme ${\rm
Spec}(\CC[\epsilon_1,\ldots, \epsilon_k]/(\epsilon_i\epsilon_j\colon 1 \leq
i,j \leq k<g))$
is contained in $\ST$ 
is equivalent to the fact that the functions $\th$ and $\dd_i \th$
for $i=1,\ldots,k$ satisfy \eqref {eq:due}. But this the same as saying that
$\dd_i \dd_j \th(\t_0,z_0)$
is zero for $i=1,\ldots,k$, $j=1,\ldots,g$, and the vectors $e_i$,
$i=1,\ldots,k$, belong to the kernel of $Q_\xi$. This settles the first
assertion.

The scheme  $\ST$ is defined by the equations \eqref {eq:singt}, where $\t$
is now fixed and $z$ is the variable. By differentiating, and using the same
notation as above, we see that the equations for the Zariski tangent space
to $\ST$ at $x$ are
$\sum _{i=1} ^ g b_i \dd_i \dd_j \th(\t_0,z_0)$, $j=1,\ldots,g$
i.e., $b\cdot M=0$, which proves the second assertion. 
\end{proof}

\end{section}

\begin{section}{Curvi-linear subschemes in the singular 
locus of theta}\label{sec:curv}

A $0$--dimensional curvi-linear subscheme
${\rm Spec}(\CC[t]/(t^{N+1})) \subset X$
of length $N+1$ supported at $x$ is given by a homomorphism
\begin{equation}\label{eq:map}
\delta: \Ocal_{X,x} \to \CC
[t]/(t^{N+1}), \qquad
f \mapsto  \sum_{j=0}^N \Delta^{(j)}f(x)\cdot  t^j,
\end{equation}
with $\Delta^{(j)}$ a differential operator of order $\leq j$, $j=1,\ldots,N$,
with $\Delta^{(N)}$ non--zero, and $\Delta^{(0)}(f)=f(x)$.
The condition that the map $\delta$ is a homomorphism is equivalent to saying
that
\begin{equation}\label {eq:prod}
\Delta^{(k)}(fg)= \sum_{r=0}^k \Delta^{(r)}\, f \cdot \Delta^{(k-r)}\, g,
\quad k=0,\ldots,N
\end{equation}
for any pair $(f,g)$ of elements of $\Ocal_{X,x}$.
Two such homomorphisms
$\delta$ and $\delta'$ define the same subscheme if and only if
they differ by composition with a automorphism of $\CC
[t]/(t^{N+1})$.

\begin{lemma}\label {pr:mapdelta}
The map $\delta$ defined in \eqref
{eq:map} is a homomorphism if and only if there exist translation
invariant vector fields $D_1,\ldots, D_N$ on $X$ such that for every
$k=1,\ldots,N$ one has
\begin{equation}\label{eq:dkappa}
\Delta^{(k)}=\sum_{h_1+2h_2+\ldots+kh_k=k>0}
\frac{1}{h_1!\cdots h_k!} D_1^{h_1} \cdots D_k^{h_k}.
\end{equation}
Moreover, two $N$--tuples of vector fields $(D_1,\ldots, D_N)$ and
$(D'_1,\ldots, D'_N)$ determine the same $0$--dimensional
curvi-linear subscheme of $X$ of length $N+1$
supported at a given point $x\in X$  if and only if there are constants
$c_1,\ldots,c_N$, with $c_1\neq 0$, such that
$$
D'_i=\sum _{j=1}^ i c_j^ {i-j+1} D_j, \quad i=1,\ldots,N.
$$
\end{lemma}

\begin{proof}
If the differential operators $\Delta^{(k)}$, $k=1,\ldots,N$, are
as in \eqref {eq:dkappa}, one computes that \eqref {eq:prod} holds, hence
$\delta$ is a homomorphism.

As for the converse, the assertion trivially holds for $k=1$.
So we proceed by induction on $k$.
Write $\Delta^{(k)}=\sum_{i=1}^ k D^{(k)}_i$, where
$D^{(k)}_i$ is the homogeneous part of degree $i$, and write $D_k$
instead of  $D^{(k)}_1$. Using \eqref {eq:prod} one verifies that for every
$k=1,\ldots,N$ and every positive $i\leq k$ one has

$$iD^{(k)}_i=\sum _{j=1}^ {k-i+1} D_j D^{(k-j)}_{i-1}.$$
Formula \eqref {eq:dkappa} follows by induction and easy combinatorics.

To prove the final assertion, use the fact that an automorphism of $\CC
[t]/(t^{N+1})$ is determined by the image $c_1t+c_2t^ 2+\ldots+c_Nt^ N$
of $t$, where $c_1\neq 0$.
\end{proof}

In formula \eqref {eq:dkappa} one has $h_k \leq 1$. If
$\Delta^{(1)}=D_1$ then $\Delta^{(2)}=\frac{1}{2}D_1^2+D_2$,
$\Delta^{(3)}= (1/3!)D_1^3+(1/2)D_1D_2+D_3$ etc.

Each non-zero summand in \eqref{eq:dkappa} is of the form
$(1/h_{i_1}!\cdots h_{i_\ell}!) D_{i_1}^{h_{i_1}} \cdots
D_{i_\ell}^{h_{i_\ell}}$, where $1\leq i_1<\ldots<i_\ell\leq k$,
$i_1h_{i_1}+\ldots+i_\ell h_{i_\ell}=k$ and $h_{i_1},\ldots,h_{i_\ell}$ are
positive integers. Thus formula \eqref {eq:dkappa} can be written as

\begin{equation}\label{eq:dkappa'}
\Delta^{(k)}=\sum_{\{h_{i_1},\ldots, h_{i_\ell}\}}
\frac{1}{h_{i_1}!\cdots h_{i_\ell}!} D_{i_1}^{h_{i_1}} \cdots
D_{i_\ell}^{h_{i_\ell}},
\end{equation}
where the subscript $\{h_{i_1},\ldots,h_{i_\ell}\}$
means that the sum is taken over all $\ell$--tuples
of positive integers $(h_{i_1},\ldots,h_{i_\ell})$  with
$1\leq i_1<\cdots <i_\ell\leq k$ and $i_1h_{i_1}+\cdots +i_\ell
h_{i_\ell}=k$.

\begin{remark}\label{rem:sisteq}
Let $x\in X$ correspond to the pair
$(\t_0,z_0)$. The differential operators  $\Delta^{(k)}$,
$k=1,\ldots,N$, defined
as in \eqref {eq:dkappa} or \eqref {eq:dkappa'} have the following
property: if $f$ is a regular function such that $\Delta^{(i)}\, f$
satisfies \eqref{eq:due} for all
$i=0,\ldots,k-1$,  then one has
$\Delta^ {(k)}f(\t_0,z_0)=0$.
\end{remark}

We want now to express the conditions in order that a $0$--dimensional
curvi-linear subscheme of $X$ of length $N+1$
supported at a given point $x\in X$ corresponding to the pair $(\t_0,z_0)$
and determined by a given $N$--tuple of vector fields $(D_1,\ldots, D_N)$
lies in $\ST$. To do so, we keep the notation we introduced above.

Let us write $D_i=\sum_{\ell=1}^ g \eta_{i\ell}\partial_\ell$, so that $D_i$
corresponds to the vector $\eta_i=(\eta_{i1},\ldots,\eta_{ig})$.
As before we denote by $M$ the matrix $(\partial_i\partial_j \theta(\tau_0,z_0))$.

\begin{proposition}\label {pr:singsing}
The $0$--dimensional curvi-linear subscheme $R$ of $X$ of length $N+1$,
supported at the point $x\in X$ corresponding to the pair $(\t_0,z_0)$
and determined by the $N$--tuple of vector fields $(D_1,\ldots, D_N)$
lies in $\ST$ if and only if $x\in \ST$ and moreover for each $k=1,\ldots,N$
one has
\begin{equation}\label {eq:cond}
\sum_{\{h_{i_1},\ldots, h_{i_\ell}\}}
\frac{1}{h_{i_1}!\cdots h_{i_\ell}!} \eta_{i_\ell}\cdot
\partial_{\eta_{i_1}}^{h_{i_1}} \cdots
\partial_{\eta_{i_\ell}}^{h_{i_\ell-1}}\, M=0, \end{equation}
where the sum is taken
over all $\ell$--tuples of positive integers $(h_{i_1},\ldots,h_{i_\ell})$
with  $1\leq i_1<\cdots <i_\ell\leq k$ and $i_1h_{i_1}+\cdots +i_\ell
h_{i_\ell}=k$.\end{proposition}
\begin{proof} 
The scheme $R$ is contained in $\ST$ if and only if one has
$$
\Delta^ {(k)}\theta(\t_0,z_0) =0, \quad
\dd_j \Delta^ {(k)} \theta(\t_0,z_0)=0 \qquad
k=0,\ldots,N,\quad j=1,\ldots,g. 
$$
By Remark \ref{rem:sisteq}, this is equivalent to
$$
\theta(\t_0,z_0) =0,\quad
\dd_j \Delta^ {(k)} \theta(\t_0,z_0) =0 \qquad k=0,\ldots,N,\quad
j=1,\ldots,g.
$$
The assertion follows by the expression \eqref{eq:dkappa'} of the operators
$\Delta^ {(k)}$.
\end{proof}

For instance, consider the scheme $R_1$, supported at $x\in \ST$,
corresponding to the vector field $D_1$. Then $R_1$ is contained in $\ST$ if
and only if
\begin{equation}\label {eq:a1}
\eta_1 \cdot M=0.
\end{equation}
This agrees with Proposition \ref {kerQ}. If $R_2$ is the scheme supported
at $x$ and  corresponding to the pair of vector fields $(D_1,D_2)$, then $R_2$
is contained in $\ST$ if and only if, besides \eqref {eq:a1} one has also
\begin{equation}\label {eq:a2}
(1/2) \eta_1 \cdot \partial _{\eta_1} M + \eta_2 \cdot M=0.
\end{equation}
Next, consider the scheme $R_3$ supported
at $x$ and  corresponding to the triple of vector fields $(D_1,D_2,D_3)$.
Then $R_3$ is contained in $\ST$ if and only if, besides \eqref {eq:a1} and
\eqref {eq:a2} one has also
\begin{equation}\label {eq:a3}
(1/3!)\eta_1 \cdot \partial^2_{\eta_1}M+(1/2)\eta_2 \cdot \partial_{\eta_1}
M+ \eta_3 \cdot M=0\end{equation}
and so on.
Observe that \eqref {eq:cond} can be written in more than one way. For example
$\eta_2 \cdot \partial_{\eta_1} M=\eta_1 \cdot \partial_{\eta_2} M$ so that
\eqref {eq:a3} could also be written as
$$
(1/3!)\eta_1 \cdot \partial_{\eta_1}^2 M+(1/2)\eta_1 \cdot \partial_{\eta_2} M+
\eta_3 \cdot M=0.
$$

So far we have been working in a fixed abelian variety $X$. One can remove
this restriction by working on $S_g$ and by letting the vector
fields $D_1,\ldots, D_N$ vary with $X$, which means that we let the vectors
$\eta_i$ depend on the variables $\tau_{ij}$. Then
the equations \eqref  {eq:cond} define a subscheme $S_g(D)$
of ${\rm Sing}(\Theta)$ which, 
as a set, is the locus of all points $\xi=(X,\Theta_X,x) \in S_g$ such that
$\ST$ contains a curvi--linear scheme of length $N+1$ supported at $x$,
corresponding to the $N$--tuple of vector fields $D=(D_1,\ldots, D_N)$,
computed on $X$. 

One can compute the Zariski tangent space to $S_g(D)$ at a point
$\xi=(X,\Theta_X,x)$ in the same way, and with the same notation,
as in \S\ref{sec: essgi}.
This gives in general a complicated set of equations. However we indicate
one case in which one can draw substantial information from 
such a computation. 
Consider indeed the case in which $D_1=\ldots=D_N\neq 0$, and call $b$ the
corresponding tangent vector to $X$ at the origin, depending on the 
the variables $\tau_{ij}$.  In this case we use the
notation $D_{b,N}=(D_1,\ldots,D_N)$ and we denote by $R_{x,b,N}$ the
corresponding curvi--linear scheme supported at $x$.
For a given such $D=(D_1,\ldots,D_N)$, consider the
linear system of quadrics
$$
\Sigma_{\xi}(D)=\PP( {\rm Im} (d\pi_{|S_g(D),\xi})^\perp)
$$
in $\PP(T_{X,0})$. One has again an interpretation of these quadrics
in terms of the normal space:

\begin{proposition}\label{prop:deformations1} 
In the above setting, the space $\Sigma_{\xi}(D_{b,N})$ contains the 
quadrics $Q_\xi, \dd_bQ_\xi,\ldots,\dd_b^NQ_\xi$.
\end{proposition}
\begin{proof} The equations \eqref{eq:cond} take now the form
$$ 
\begin{aligned}
\theta(\tau,z)=0, \, \partial_i\theta(\tau,z)=0,\quad i=1,\dots,g\cr
b\cdot M=b\cdot \partial_bM=\dots= b\cdot \partial_b^{N-1} M=0.\cr
\end{aligned}
$$
By differentiating the assertion immediately follows.
\end{proof}
\end{section}

\begin{section}{Higher multiplicity points of the theta divisor}
\label{higherorder}
We now study the case of higher order singularities on the theta divisor.
For a multi-index $I=(i_1,\ldots,i_g)$ with $i_1,\ldots,i_g$ non-negative
integers we  set $z^ I=z_1^ {i_1} \cdots z_g^ {i_g}$ and
denote by $\dd_I$ the operator $\dd^ {i_1}_{1}\cdots \dd_g^ {i_g}$. Moreover,
we let $|I|= \sum_{\ell=1}^gi_\ell$, which is the {\it length} of $I$
and equals the {\it order} of the operator $\dd_I$.

\begin{definition}\label{def:sgr}
For a positive integer $r$ we let $S_g^{(r)}$ be the subscheme of
$\Xg$ which is defined on  $\HH_g \times \CC^g$ by the equations
\begin{equation}\label{eq:singrt}
\dd_{I}\th(\t,z)=0, \qquad |I|=0,\ldots,r-1.
\end{equation}
\end{definition}

One has the chain of subschemes
$$
\ldots \subseteq S_g^{(r)}\subseteq \ldots
\subseteq
S_g^{(3)}\subseteq S_g^{(2)}=S_g\subset S_g^{(1)}=\Theta
$$
and as a set
$ S_g^{(r)}=\{ (X,\Theta_X,x)\in \Theta :  \text{$x$ has multiplicity
$\geq r$ for $\Theta_X$} \} $. One denotes by ${\rm Sing}^
{(r)}(\Theta_X)$ the subscheme of ${\rm Sing}(\Theta_X)$ formed by all
points of multiplicity at least $r$.
One knows that $S_g^{(r)}=\emptyset$ as soon as $r>g$ (see
\cite {SV}). 
We can compute the Zariski tangent space to $S_g^{(r)}$ at a point
$\xi=(X,\Theta_X,x)$ in the same vein, and with the same notation, as in \S
\ref {sec: essgi}. Taking into account that $\theta$ and all its derivatives
verify the heat equations, we find the equations by replacing in (3) 
the term $\theta(\tau_0,z_0)$ by $\partial_I \theta(\tau_0,z_0)$.

As in \S \ref{sec: essgi}, we wish to give some geometrical
interpretation. For instance, we have the following lemma which
partially extends Lemma \ref{lem:triplept} or \ref{prop:sing}.

\begin{lemma}\label {lem:singg} For every positive integer $r$ the scheme
$S_g^{(r+2)}$ is contained in the singular locus of $S_g^{(r)}$.
\end{lemma}

Next we are interested in the differential of the restriction of the map
$\pi:\Xg\to \Ag$ to $S_g^{(r)}$ at a point $\xi=(X,\Theta_X,x)$ which does
not belong to $S_g^{(r+1)}$. This means that $\Theta_X$ has a point of
multiplicity exactly $r$ at $x$.
If we assume, as we may, that $x$
is the origin of $X$ , i.e. $z_0=0$, then the Taylor expansion of
$\theta$ has the form
$$
\th = \sum_{i=r}^{\infty}\th_i,
$$
where $\th_i$ is a homogeneous polynomial of degree $i$ in the variables
$z_1,\ldots,z_g$ and
$$
\theta_r= \sum_{I=(i_1,\ldots,i_g), |I|=r} \frac{1}{i_1!\cdots i_g!} 
\dd_I \theta(\t_0,z_0) z^I
$$
is not identically zero. The equation
$\theta_r=0$ defines a hypersurface $TC_\xi$ of degree $r$ in $\PP^
{g-1}=\PP(T_{X,0})$, which is the {\it tangent cone} to $\Theta_X$ at $x$.

We will denote by ${\rm Vert}(TC_\xi)$ the {\sl vertex} of $TC_\xi$, i.e.,
the subspace of $\PP^ {g-1}$ which is the locus of points of
multiplicity $r$ of $TC_\xi$. Note that it may be empty. In case $r=2$,
the tangent cone
$TC_\xi$ is the quadric $Q_\xi$ introduced in \S  \ref{sec: essgi} and ${\rm
Vert}(TC_\xi)$ is its vertex $\Pi_\xi$.

More generally, for every $s\geq r$, one
can define the subscheme $TC^ {(s)}_\xi=TC^ {(s)}_x$ of $\PP^
{g-1}=\PP(T_{X,0})$ defined by the equations
$$\theta_r=\ldots =\theta_{s}=0,$$
which is called the {\it asymptotic cone} of order $s$ to $\Theta_X$ at
$x$.

Fix a multi-index $J=(j_1,\ldots,j_g)$ of length $r-2$. For any pair $(h,k)$
with $1\leq h,k\leq g$, let $J_{(h,k)}$ be the multi-index of length $r$
obtained from $J$ by first increasing by $1$ the index $j_h$ and then by $1$
the index $j_k$ (that is, by $2$ if they coincide). 
Consider then the quadric $Q_\xi^ J$ in $\PP^{g-1}=\PP(T_{X,0})$ 
defined by the equation
$$
q_\xi^ J(z):=\sum _{1\leq h,k\leq g} \dd_ {J_{ (h,k) }} \theta
(\t_0,z_0)z_hz_h=0
$$
with the usual convention that the quadric is
indeterminate if the left hand side is identically zero. This is a {\it polar
quadric} of $TC_\xi$, namely it is obtained from $TC_\xi$ by iterated
operations of polarization. Moreover all polar quadrics are in the span
$\langle Q_\xi^ J,|J|=r-2\rangle$. We will denote by
$\Qcal^ {(r)}_\xi$ the span of all quadrics $Q_\xi^ J$, $|J|=r-2$
and $\dd_bQ_\xi^ J$, $|J|=r-2$, with equation
$$
\sum _{1\leq i,j\leq g} \dd_b\dd_{J_{(h,k)}}\theta (\t_0,z_0)z_iz_j=0,
$$
for every non--zero vector $b\in \CC^ g$ such that $[b]\in {\rm
Vert}(TC_\xi)$.

We are now interested in the kernel and the image of 
$d\pi_{|S^{(r)}_g,\xi}$.  Equivalently we may
consider the linear system of quadrics  $\Sigma^{(r)}_\xi=\PP( {\rm Im}
(d\pi_{|S^ {(r)}_g,\xi})^ \perp)$, and the subspace
$\Pi^{(r)}_\xi=\PP( {\rm ker} (d\pi_{|S^ {(r)}_g,\xi}))$ of
$\PP(T_{X,0})$. The following proposition partly extends Proposition \ref
{prop:deformations} and \ref{kerQ} and its proof is similar.

\begin{proposition}\label{prop:deformations2} 
Let $\xi=(X,\Theta_X,x)$ be a point of $S^ {(r)}_g$. Then:
\begin{itemize}
\item [(i)] $\Pi^ {(r)}_\xi={\rm Vert}(TC_\xi)$. In particular, if
$\xi\in S^ {(r+1)}_g$, then  $\Pi^ {(r)}_\xi$ is the whole space
$\PP(T_{X,0})$;
\item [(ii)] $\Sigma^{(r)}_\xi$ contains the linear system
$\Qcal^{(r)}_\xi$.
\end{itemize}
\end{proposition}

\begin{remark} \label {rem:tanspace}
{\rm As a consequence, just like in Proposition \ref{kerQ}, 
one sees that for $\xi=(X,\Theta_X,x)$ the Zariski tangent 
space to ${\rm Sing}^{(r)}(\Theta_X)$ at $x$ is contained 
in ${\rm Vert}(TC_\xi)$.}
\end{remark}

As an application, we have:

\begin{proposition}\label{prop:lower} Let $\xi=(X,\Theta_X,x)$ be a
point of  a component $Z$ of
$S^ {(3)}_g$ such that $TC_\xi$ is not a cone. Then
$\dim({\mathcal Q}^{(3)}_\xi)= g-1$ and therefore the codimension
of the image of $Z$ in ${\mathcal A}_g$ is al least $g$.
\end{proposition}

\begin{proof}  
Since $TC_\xi$ is a not a cone its polar quadrics are 
linearly independent. 
\end{proof} 

The following example shows that  the above bound is
sharp for $g=5$.

\begin{example}
Consider the locus $\Ccal$ of intermediate Jacobians of cubic threefolds in
${\mathcal A}_5$. These have an isolated triple point on their theta
divisor whose tangent cone gives back the cubic threefold. The locus
$\Ccal$ is dominated by an irreducible component of $S_5^ {(3)}$ for which
the estimate given in (i) of Proposition \ref {prop:lower} is sharp. Cf.\  \cite{Ca-M}
where Casalaina-Martin proves that the locus of intermediate Jacobians
of cubic threefolds is an irreducible component of the locus of 
principally polarized abelian varieties of dimension $5$ with a point 
of multiplicity~$\geq 3$. 
\end{example}
\end{section}

\begin{section}{The Andreotti--Mayer loci}\label{sec:amloci}

Andreotti and Mayer consider in $\Ag$ the algebraic sets of
principally polarized abelian varieties $X$ with a 
locus of singular points on $\Theta_X$ of dimension
at least $k$. More generally, we are
interested in the locus of principally polarized abelian varieties
possessing a $k$-dimensional locus of singular points of multiplicity $r$
on the theta divisor. To define these loci scheme-theoretically we consider the
morphism $\pi: \Xg \to \Ag$ and the quasi-coherent sheaf on $\Ag$
$$
\Fcal^ {(r)}_k=\bigoplus_{i=k}^{g-2} R^i\pi_*\Ocal_{S^ {(r)}_g}.
$$
\begin{definition}
For integers $k$ and $r$ with $0\leq k\leq g-2$
and $2\leq r\leq g$ we define $N_{g,k,r}$ as the support of
$\Fcal^ {(r)}_k$. We also set  $M_{g,k,r}=\pi^ *(N_{g,k,r})$, 
a subscheme of both $S_{g,k}$ and $S_g^ {(r)}$. 
We write $N_{g,k}$ and $M_{g,k}$ for $N_{g,k,2}$ and $M_{g,k,2}$. 
\end{definition}

The schemes $N_{g,k}$ are the so--called {\it Andreotti--Mayer loci} in $\Ag$,
which have been introduced in a somewhat different way in \cite {AM}.

Note that $N_{g,k,r}$ is locally defined by an annihilator ideal 
and so carries the structure of subscheme.  
Corollary \ref{cor:sgsmooth} below and results by Debarre \cite{D1} 
(see \S \ref{sec:AMbound})  imply that the scheme structure 
 at a general point of $N_{g,0}$ 
defined above coincides with the one considered by Mumford in
\cite{Mumford}. 

We now want to see that as a set $N_{g,k,r}$ is the locus of points
corresponding to $(X,\Theta_X)$ such that $\ST$ has an irreducible component
of dimension $k$ of points of multiplicity $r$ for $\Theta_X$.

\begin{lemma}
Let $X$ be an abelian variety of dimension $g$ and $W\subset X$
an irreducible reduced subvariety of dimension $n$ and let $\omega_W$
be its dualizing sheaf. Then $H^0(W,\omega_W)\neq (0)$.
\end{lemma}
\begin{proof}
Let $f: W^{\prime}\to W$ be the normalization of $W$.
We first claim the inequality
$h^0(W',\omega_{W^{\prime}})\leq h^0(W,\omega_W)$.
To see this, note that by \cite{Li}, p.\ 48 ff 
(see also \cite{H}, Exerc. 6.10, p.\ 239, 7.2, p.\ 249),
there exists a map
$f_*\omega_{W^{\prime}} \to {\mathcal Hom}(f_*\Ocal_{W^{\prime}},\omega_W)$,
hence a map 
$$
H^0(W^{\prime},\omega_{W^{\prime}})\to
H^0(W,{\mathcal Hom}(f_*{\Ocal}_{W^{\prime}},\omega_W)).
$$
Now $\Ocal_{W} \to f_*\Ocal_{W^{\prime}}$ is an injection and therefore
$H^0(W,{\mathcal Hom}(f_*\Ocal_{W^{\prime}},\omega_W)$ maps to
$H^0(W,{\mathcal Hom}(\Ocal_W,\omega_W))$
and we thus get a map $H^0(W^{\prime},\omega_{W^{\prime}}) \to
H^0(W,\omega_W)$ which is injective as one sees by looking at the
smooth part of $W$.

Let $\tilde{W}$ be a desingularization of $W^{\prime}$.
According to  \cite{Ku} we have $h^0(\tilde{W},\Omega_{\tilde{W}}^n) \leq
h^0(W^{\prime},\omega_{W^{\prime}})$. Since $\tilde{W}$ maps to $X$
we have $h^0(\tilde{W},\Omega_{\tilde{W}}^1)\geq n$. If $h^0(\tilde{W},
\Omega_{\tilde{W}}^n)$ were 
$0$ then $\wedge^n H^0(\tilde{W},\Omega_{\tilde{W}}^1)
\to H^0(\tilde{W},\Omega_{\tilde{W}}^n)$ would be the zero map contradicting
the fact that $W$ has dimension~$n$.
\end{proof}

\begin{corollary}
We have $(X,\Theta_X) \in N_{g,k,r}$ if and only if 
$\dim ({\rm Sing}^{(r)}(\Theta_X)) \geq k$.
\end{corollary}
\begin{proof} By the previous lemma and Serre duality
for a reduced irreducible subvariety $W$ of dimension $m$ in $X$ 
it follows that $H^{m}(W,{\Ocal}_W)\neq (0)$ and we know 
$H^k(W,{\Ocal}_W)=(0)$ for $k>m$. This implies the corollary.
\end{proof}

There are the inclusions
$$
N_{g,k,r}\subseteq N_{g,k,r-1},\quad N_{g,k,r}\subseteq N_{g,k-1,r}.
$$
If $p=(n_1,\ldots,n_r)$ with $1\leq n_1 \leq \ldots \leq n_r<g$
and $n_1+\ldots + n_r=g$ is a partition of $g$ we write $\Acal_{g,p}$
for the substack of $\Ag$ corresponding to principally polarized
abelian varieties that are a product of $r$ principally polarized abelian 
varieties of dimensions $n_1,\ldots,n_r$. We write $r(p)=r$ for the length
of the partition and write ${\Acal}_{g,[r]}$ for the substack $\cup_{r(p)=r}
{\Acal}_{g,p}$ of $\Ag$ corresponding to pairs $(X,\Theta_X)$ 
isomorphic as a polarized abelian variety to the product of $r$ 
principally polarized abelian varieties.
One has the stratification
$$
{\Acal}_{g,[g]}\subset {\Acal}_{g,[g-1]} \subset \cdots
 \subset {\Acal}_{g,[2]}.
$$
We will denote by:

i)  $\Pi_g=\cup_{p} {\Acal}_{g,p}$ the locus
of {\it decomposable} principally polarized abelian varieties;

ii) $\Ag^{({\rm NS})}$ the locus of classes of {\it non-simple
abelian varieties}, i.e., of principally polarized abelian varieties of
dimension $g$ which are isogenous to a product of abelian varieties of
dimension smaller that $g$;

iii) $\Ag^{{\rm End}\neq \ZZ}$ 
the locus of classes of {\it singular abelian varieties},
i.e., of principally polarized abelian varieties whose endomorphism ring is
larger than $\ZZ$.

\begin{remark}\label{rem:fac} 
Note the inclusions
$\Pi_g\subset \Ag^{({\rm NS})}\subset \Ag^{{\rm End} \neq \ZZ}$.
The locus $\Pi_g$ is reducible with irreducible components
${\Acal}_{g,p}$  with $p$ running through the partitions $g=(i,g-i)$
of $g$ for $1\leq i\leq g/2$ and we have 
${\rm codim}_{\Ag}({\Acal}_{g,(i,g-i)})=i(g-i)$.
In contrast to this  $\Ag^{{\rm End}\neq \ZZ}$ and 
$\Ag^{({\rm NS})}$ are the union of infinitely countably many irreducible 
closed subsets of $\Ag$ of codimension at least
$g-1$, the minimum codimension being achieved for 
families of abelian varieties that are isogenous to products of 
an elliptic curve with an abelian variety of dimension $g-1$ 
(compare with \cite{cvdg}).
\end{remark}

We recall a result from \cite{Ko} and the main result from \cite{EL}. 

\begin{theorem}\label{th:mainel} 
For every integer $r$ with  $2\leq r \leq g$ one has:
\begin{itemize}
\item[(i)] $N_{g,k,r}=\emptyset$ if $k>g-r$;
\item [(ii)] $N_{g,g-r,r}={\Acal}_{g,[r]}$, i.e., $(X,\Theta_X) \in 
N_{g,g-r,r}$ is an $r$-fold product.
\end{itemize}
\end{theorem}
Hence, for every integer $r$ such that $2\leq r\leq g$, one has the
stratification
$$
N_{g,0,r}\supset N_{g,1,r}\supset \ldots \supset N_{g,k,r}\supset \ldots 
\supset N_{g,g-r,r}={\Acal}_{g,[r]},
$$
whereas for every integer $k$ such that $0\leq k\leq g-2$, one has 
the stratification
$$
N_{g,k,2}\supset N_{g,k,3}\supset \ldots\supset N_{g,k,r}\supset \ldots
\supset N_{g,k,g-k}={\Acal}_{g,[g-k]}.
$$
\end{section}
\begin{section}{Lower bounds for the codimension of Andreotti-Mayer loci}\label{lowerbounds}
The results in the previous sections give information about the Zariski
tangent spaces to these loci and this will allow us to
prove bounds on the dimension of the Andreotti--Mayer loci, which is our
main objective in this paper. 

We start with the results on tangent spaces. We need some notation.
\begin{definition}
Let $\zeta=(X,\Theta_X)$ represent a point in $N_{g,k,r}$.
By $\Lcal_{g,k,r}(\zeta)$ we denote the
linear system of quadrics $\PP(T_{N_{g,k,r},\zeta}^\perp)$, 
where $T_{N_{g,k,r},\zeta}$ is the Zariski tangent space and where we view 
$\PP(T_{{\Acal}_{g,\zeta}})$ as a space of quadrics as in Section 
\ref{sec: univtheta}. As usual, we may drop the index~$r$ if $r=2$ and 
write ${\Lcal}_{g,k}(\zeta)$ for ${\Lcal}_{g,k,2}(\zeta)$. 
\end{definition}

Notice that
$$
{\rm dim}_{\zeta}(N_{g,k,r})\leq \binom{g+1}{2}-\dim
(\Lcal_{g,k}(\zeta))+1.
$$ 
\begin{definition}
For $\zeta=(X,\Theta_X) \in N_{g,k,r}$ we denote by 
${\rm Sing}^{(k,r)}(\Theta_X)$ the locally closed subset 
$$
{\rm Sing}^{(k,r)}(\Theta_X)=
\{ x \in \ST : \dim_x({\rm Sing}^{(r)}(\Theta_X))\geq k \}.
$$ 
Moreover, we define $\Qcal_\zeta^{(k,r)}$ to be the linear system of 
quadrics in $\PP^ {g-1}=\PP(T_{X,0})$ spanned by the union of all 
linear systems $\Qcal^{(r)}_\xi$  with $\xi=(X,\Theta_X,x)$ and 
$x\in {\rm Sing}^{(k,r)}(\Theta_X)$.
\end{definition}  

Propositions \ref{prop:deformations} and \ref{prop:deformations2} imply
the following basic tool for giving upper
bounds on the dimension of the Andreotti--Mayer loci.

\begin{proposition}\label{prop:tanam} 
Let $N$ be an irreducible component of $N_{g,k,r}$ with its reduced
structure.  If $\zeta = (X, \Theta_X)$ is a general point of $N$ 
then the projectivized conormal  space to $N$ at $\zeta$, viewed as a
subspace of  $\PP({\rm Sym}^2(T_{X,0}^{\vee}))$, contains the linear 
system $\Qcal_\zeta^{(k,r)}$. 
\end{proposition}
\begin{proof}
Let $M$ be an irreducible component of $\pi^*N$ in $S_{g,k}^{(r)}$.
If $\xi$ is  smooth point of $M$ then the image of the Zariski tangent
space to $M$ at $\xi$ under $d\pi$ is orthogonal to $\Qcal_{\xi}^{(r)}$
for all $x \in {\rm Sing}^{k,r}(\Theta_X)$. Since we work in characteristic
$0$ the map $d\pi$ is surjective on the tangent spaces for general
points $m\in M$ and $\pi(m) \in N$. Therefore the result follows from 
Propositions \ref{prop:deformations} and \ref{prop:deformations2}.
\end{proof}

We need a couple of preliminary results. First we state a well
known fact, which can be proved easily by a standard dimension count.

\begin{lemma}\label{lem:dimcount} Every hypersurface of degree $d\leq 2n-3$
in $\PP^n$ with  $n\geq 2$ contains a line. \end{lemma}

Next we prove the following:

\begin{lemma}\label {lem:polarity} Let $V$ in $\PP^ s$ be a hypersurface
of degree $d\geq 3$. If all polar quadrics of $V$ coincide, then $V$ is
a hyperplane $H$ counted with multiplicity $d$, and the polar quadrics
coincide with $2H$.
\end{lemma}

\begin{proof}  If $d=3$ the assertion follows from general properties of duality 
(see \cite {zak}, p. 215) or  from an easy calculation. 

If $d>3$, then the result, applied to the cubic polars of $V$, tells us that all these
cubic polars are equal to $3H$, where $H$ is a fixed hyperplane. This immediately
implies the assertion.\end{proof}

The next result has been announced in \cite{C-vdG}. 

\begin{theorem}\label{th:amlower} Let $g\geq 4$ and let $N$ be an
irreducible component of $N_{g,k}$ not contained in $N_{g,k+1}$. Then:
\begin{itemize}
\item [(i)] for every positive integer $k\leq g-3$, one has ${\rm
codim}_{\Ag}(N)\geq k+2$, whereas  ${\rm codim}_{\Ag}(N)= g-1$ if $k=g-2$;
\item [(ii)] if $N$ is contained in $N_{g,k,r}$ with $r\geq 3$,
 then ${\rm codim}_{\Ag}(N)\geq k+3$;
\item [(iii)] if  $g-4\geq k\geq g/3$,
then ${\rm codim}_{\Ag}(N)\geq k+3$.

\end{itemize}
\end{theorem}

\begin{proof} By Theorem \ref{th:mainel} and Remark \ref{rem:fac}, we
may assume $k<g-2$. 
By definition, there is some irreducible component $M$ of
$\pi_{S_g}^{-1}(N)$ with $\dim(M)=\dim(N)+k$ which dominates $N$ via
$\pi$. We can take a general point $(X,\Theta_X,z)\in M$ so that 
$\zeta=(X, \Theta_X)$ is a general point in $N$. By Remark \ref {rem:fac}
we may assume $X$ is simple. 

Let $R$ be the unique $k$--dimensional
component of $\pi_{|S_g}^{-1}(\zeta)$ containing $(X,\Theta_X,z)$. Its
general point is of the form  $\xi=(X,\Theta_X,x)$ with $x$ the general
point of the unique $k$--dimensional component of $\ST$ containing
$z$, and $x$ has multiplicity $r$ on $\Theta_X$.
By abusing notation, we may still denote this component by $R$.
Proposition \ref{prop:deformations2} implies that the linear
system of quadrics $\PP(T_{N,\zeta}^ \perp)$ contains all polar quadrics 
of $TC_\xi$ with $\xi=(X,\Theta_X,x)\in R$. 

Thus we have a rational map
$$
\phi: (\PP^{g-1})^ {r-2}\times R \map \Qcal^ {(r)}_\zeta
$$
which sends the general point $(b,\xi):=(b_1,\dots,b_{r-2},\xi)$ 
to the polar quadric $Q_{b,\xi}$ of $TC_\xi$ with respect to 
$b_1,\dots,b_{r-2}$. It is useful to remark that
the quadric $Q_{b,\xi}$ is a cone with vertex containing the projectivized 
Zariski tangent space to $R$
at $x$ (see Remark \ref{rem:tanspace}).

\begin{claim}\label{cl:finite} [The finiteness property.] 
For each $b\in  (\PP^{g-1})^ {r-2}$ the map $\phi$ restricted to 
$\{b\}\times R$ has finite fibres.
\end{claim} 

\begin{proof} [Proof of the claim] Suppose the assertion is not true. Then
there is an irreducible curve $Z\subseteq R$ such that, for 
$\xi=(X,\Theta_X,x)$ corresponding to the general point in $Z$, one has
$Q_{b,\xi}=Q$. Set $\Pi={\rm Vert}(Q)$, which is a proper subspace of $\PP^
{g-1}$.

Consider the {\sl Gauss map}
$$
\gamma=\gamma_Z: Z\map \PP^ {g-1}=\PP(T_{X,0})
$$
which associates to a smooth point of $Z$ its projectivized tangent
space. Then Proposition \ref {prop:deformations2} implies that $\gamma(\xi)\in R\subseteq 
 {\rm Vert}(Q_\xi)=\Pi$. Thus $\gamma(Z)$ is degenerate in $\PP^ {g-1}$ and this
yields that $X$ is non--simple, cf.\ \cite{Ran}. 
This is a contradiction which proves the claim. 
\end{proof}

Claim \ref{cl:finite} implies that the image of the map $\phi$ has dimension
at least $k$, hence ${\rm codim}_{\Ag}(N_{g,k})\geq k+1$. 
To finish the proof of part (i), we need the following information.
 
\begin{claim}\label{cl:dom} [The non--degeneracy property.] The image 
of the map $\phi$ does not contain any line.
\end{claim}
\begin{proof} [Proof of the claim]
Suppose the claim is false. Take a line in the image of the map $\phi$, 
and let $\Lcal$ be the corresponding pencil of quadrics. 
By Proposition \ref{prop:deformations2}  and Remark \ref{rem:tanspace},
the general quadric in $\Lcal$ has rank $\rho\leq
g-k$. Then part (i) of Segre's Theorem \ref {thm:segre} in 
\S \ref{sec:segre} below implies that the Gauss image $\gamma(Z)$ of any
irreducible component of the curve $Z=\phi^{-1}(\Lcal)$ is degenerate. This
again leads to a contradiction. 
\end{proof}

To prove part (ii), we want to prove that  $\dim(\Qcal^{(r)}_\zeta)>k+1$. 
Remember that the image of $\phi$ has dimension at least $k$ by 
Claim \ref{cl:finite}. If the image has dimenson at least $k+1$, then by 
Claim \ref{cl:dom} it cannot be a projective space and therefore 
$\dim(\Qcal^{(r)}_\zeta)>k+1$.  So we can assume that the image has 
dimension $k$.  Therefore each component of the fibre $F_Q$
over a general point $Q$ in the image has dimension $(g-1)(r-2)$.
 
Consider now the projection of $F_Q$ to $R$. If the image is positive 
dimensional then there is a curve $Z$ in $R$ such that the image of 
the Gauss map of $Z$ is contained in the vertex of $Q$. Then $X$ 
is non--simple, a contradiction (see the proof of Claim \ref{cl:finite}). 
 
Therefore the image of $F_Q$ on $R$ is constant, equal to a point
$\xi$, hence $F_Q= (\PP^{g-1})^{r-2}\times \{\xi\}$. 
By Lemma \ref{lem:polarity}
there is a hyperplane $H_\xi$ such
that $TC_\xi=rH_\xi$, and $Q=2H_\xi$. Therefore we have a rational map
$$
\psi: R\map{\PP^ {g-1}}^{\vee} 
$$
sending $\xi$ to $H_\xi$. We notice that the image of $\phi$ is then 
equal to the $2$--Veronese image if the image of $\psi$. 
 
By Claim \ref{cl:finite}, the map $\psi$ has finite fibres. By an argument 
as in Claim \ref{cl:dom}, we see that the image of $\psi$ does not contain
a line, hence is not a linear space. Thus it spans a space of dimension
$s\geq k+1$. Then its $2$--Veronese image, which is the image of $\phi$
spans a space of dimension at least $2s\geq 2k+2> k+2$. 

To prove part (iii), by part (ii), we can assume $r=2$. It suffices to show
that $\dim(\Qcal_\zeta)>k+1$, where $\Qcal_\zeta:=\Qcal^ {(2)}_\zeta$.
Suppose instead that  $\dim(\Qcal_\zeta)=k+1$ and set $\Sigma=\phi(R)$.  
We have to distinguish two cases:
\begin{itemize}
\item [(a)] not every quadric in $\Qcal_\zeta$ is singular;
\item [(b)] the general quadric in $\Qcal_\zeta$ 
is singular, of rank $g-\rho<g$.
\end{itemize}
In case (a), consider the discriminant $\Delta
\subset \Qcal_\zeta$, i.e. the scheme of singular quadrics in $\Qcal_\zeta$.
This is a hypersurface of degree $g$, which, by Proposition \ref{kerQ},
contains $\Sigma$ with multiplicity at least $k$.
Thus $\deg(\Sigma)\leq g/k\leq 3$ and $\Sigma$ contains some line, so that
we have the corresponding pencil $\Lcal$ of singular quadrics. 
By Claim \ref{cl:dom}, one arrives at a contradiction. 

Now we treat case (b). Let $g-h$ be the rank of the general quadric
in $\Sigma$. One has $g-h\leq g-k$, hence $k\leq h$. Moreover one has
$g-h\leq g-\rho$, i.e. $\rho\leq h$. 

Suppose first $\rho=h$, hence $\rho\geq k$. Let $s$ be the dimension of the
subspace 
$$
\Pi:={ \bigcap_{Q\in \Qcal_\zeta, {\rm rk}(Q)= g-\rho} {\rm Vert}(Q)}.
$$
By applying part (iii) of Segre's Theorem \ref{thm:segre} 
to a general pencil contained in $\Qcal_\zeta$, we deduce that
\begin{equation}\label {eq:limi}
3\rho\leq g+2s+2.
\end{equation}
%%%
\begin{claim}\label {cl:ff} One has $s<r-k$.
\end{claim}
%%%
\begin {proof} Suppose $s\geq \rho-k$. If $\Pi'$ is a general subspace
of $\Pi$ of dimension $\rho-k$, then its intersection with
$\PP(T_{R,x})$, where $x\in R$ is a general point, is not empty. Since
$\rho-k<g-k$ and  $X$ is simple, this is a contradiction (see \cite{Ran},
Lemma II.12). 
\end{proof}
By \eqref {eq:limi} and Claim \ref {cl:ff} we deduce that $\rho+2k\leq g$ and
therefore $3k\leq g$, a contradiction. 

Suppose now $\rho<h$. Then part (iv) of Segre's Theorem \ref{thm:segre} yields
$$
\deg(\Sigma)\leq  \frac{g-2-\rho}{h-\rho}.
$$
The right hand side is an increasing function of $\rho$, thus
$\deg(\Sigma)\leq g-h-1\leq 2k-1$ because $g\leq 3k\leq 2k+h$. 
By Lemma \ref{lem:dimcount} the locus $\Sigma$ contains a
line, and we can conclude as in the proof of the non--degeneracy 
property \ref{cl:dom}. 

\end{proof}

The following corollary was proved independently by Debarre \cite{De1} and
includes a basic result by Mumford \cite{Mumford}.  

\begin{corollary}\label{cor:debarrre}  
Let $g\geq 4$. Then: 
\begin{itemize}
\item [(i)] every
irreducible component $S$ of $S_g$ has codimension $g+1$ in $\Xg$, hence
$S_g$ is locally a complete intersection in $\Xg$;
\item [(ii)] if $\xi=(X,\Theta_X,x)$ is a general point of $S$, then
either $\Theta_X$ has isolated singularities or $(X,\Theta_X)$ is a product
of an elliptic curve with a principally polarized abelian variety
of dimension $g-1$;
\item  [(iii)] every irreducible component $N$ of $N_{g,0}$ has
codimension $1$ in $\Ag$
\item [(iv)] if $(X,\Theta_X)\in N$ is a general
element of an irreducible component $N$ of $N_{g,0}$, then $\Theta_X$ has
isolated singularities. Moreover for every point $x\in {\rm
Sing}(\Theta_X)$, the quadric $Q_x$ is smooth and independent of $x$. 
\end{itemize} 
\end{corollary}

\begin{proof} 
Take an irreducible component $S$ of $S_g$ and
let $N$ be its image via the map $\pi: \Xg \to \Ag$. 
Then there is a maximal integer 
$k\leq g-2$ such that $N$ is contained in $N_{g,k}$.
Suppose first that $1\leq k \leq g-3$. Then by Theorem \ref{th:amlower} the 
codimension of $N$ in $\Ag$ is at least $k+2$. This implies that the
codimension of $S$ in $\Xg$ is at least $g+2$, which is impossible because
$S_g$ is locally defined in $\Xg$ by $g+1$ equations.

Then either $k=0$ or $k=g-2$. In the former case $S$ maps to an irreducible
component $N$ of $N_{g,0}$ which is a proper subvariety of 
codimension 1 in $\Ag$. 

If $k=g-2$, then by Theorem  \ref{th:mainel} the polarized abelian variety
$(X,\Theta_X)$ is a product of principally polarized
abelian varieties. The resulting abelian variety has a 
vanishing thetanull and so $N$ is contained in 
the divisor of a modular form (a product of thetanulls, cf.\ e.g.\
\cite{Mumford}, p.\ 370) which is a component of $N_{g,0}$.  
Assertions (i)--(iii) follow by a dimension count 
(see Remark \ref{rem:fac}).
Part (iv) follows by Propositions \ref {prop:deformations} and \ref
{prop:tanam}. 
\end{proof}

Finally, we show a basic property of $S_g$. 

\begin{theorem}\label{thm:sgsmooth} 
The locus $S_g$ is reduced.
\end{theorem}

\begin{proof}  The assertion is well known for $g\leq 3$, using the theory
of curves. We may assume $g\geq 4$.

Let $S$ be an irreducible component of $S_g$, let $N=\pi(S)$
and $k\leq g-2$ the maximal integer such that $N$ is contained in 
$N_{g,k}$.
As in the proof of Corollary \ref{cor:debarrre}, one has either $k=0$
or $k=g-2$. Assume first $k=0$, and let
$\xi=(X,\Theta_X,x)\in S$ be a general point.
We are going to prove that ${\rm Sing}(\Theta_X)$ is
reduced of dimension 0.

First we prove that $x$ has multiplicity $2$ for $\Theta_X$.
Suppose this is not the case and $x$ has multiplicity $r\geq 3$. 
Since $N_{g,0}$ has codimension $1$ all polar quadrics of $TC_\xi$ 
are the same quadric, say $Q$ (see Proposition \ref{prop:deformations2}
and \ref{prop:tanam}). 
By Lemma \ref{lem:polarity} the tangent cone $TC_\xi$ is a hyperplane $H$
with multiplicity $r$. Again by Proposition \ref {prop:deformations2} 
all the derivatives of $Q$ with respect to points $b\in H$  coincide with $Q$. 
By Proposition \ref{pr:singsing} the scheme $D_{b,2}$ supported at
$x$ is contained in ${\rm Sing}(\Theta_X)$. By taking into account
 Proposition \ref  {prop:deformations1}  and repeating the same
argument we see that this subscheme can be indefinitely extended to a 
$1$--dimensional subscheme, containing $x$ and contained in 
${\rm Sing}(\Theta_X)$. This implies that the corresponding 
component of $N_{g,0}$ is contained in $N_{g,1}$, 
which is not possible since the codimension of $N_{g,1}$ is at least~$3$.

If $x$ has multiplicity $2$ the same argument shows that the quadric $Q_\xi$
is smooth. By Corollary \ref{cor:smoothsg},
$\xi$ is a smooth point of $S_g$ and this proves the assertion. 

Suppose now that $k=g-2$. Then by Theorem \ref{th:mainel}
$N$ is contained in the locus of products ${\Acal}_{g,(1+g-1)}$
and for dimension reasons equal to it and then the result
follows from a local analysis with theta functions.
\end{proof}

\begin{corollary}\label{cor:sgsmooth} 
If $(X,\Theta_X)$ is a general point in a component of $N_{g,0}$ 
then $\Theta_X$ has finitely many double points
with the same tangent cone which is a smooth quadric.
\end{corollary}

\end{section}
%%%%%%%%%%%%%%%%%%
%%%%%%%%%%%%%%%%%%
\begin{section}{A conjecture}
As shown in \cite{C-vdG}, part (i) of Theorem \ref{th:amlower}
is sharp for $k=1$ and $g=4$ and $5$. 
However, we do not expect Theorem \ref{th:amlower} to be sharp in general. 
For example, as indicated in  \cite{C-vdG}, it is never sharp for $k=1$ 
and $g\geq 6$, or for $k\geq 2$. In \cite{C-vdG} we made the following 
conjecture, which is somehow a natural completion of Andreotti--Mayer's 
viewpoint in \cite{AM} on the Schottky problem.

Recall the {\it Torelli morphism} $t_g: \Mcal_g\to \Ag$ which maps the
isomorphism class of a curve $C$ to the isomorphism class of its
principally polarized jacobian $(J(C),\Theta_C)$. As a map of orbifolds 
it is of degree $2$ for $g\geq 3$ since the general abelian variety has
an automorphism group of order $2$ and the general curve one of order $1$. 
We denote by $\Jcal_g$ the {\it jacobian locus} in $\Ag$, i.e.,
the Zariski closure of $t_g(\Mcal_g)$ in $\Ag$ and by ${\Hcal}_g$ 
the {\it hyperelliptic locus} in $\Ag$, that is, 
the Zariski closure in $\Ag$ of $t_g(H_g)$, where 
 $H_g$ is the closed subset of $\Mcal_g$ consisting of the isomorphism
classes of the hyperelliptic curves. By Torelli's theorem we have
$\dim(\Jcal_g)=3g-3$ and $\dim(\Hcal_g)=2g-1$ for $g\geq 2$.
\begin{conjecture}\label {conj:mainconj} If $1 \leq k \leq g-3$ and if $N$ is
an irreducible component of $N_{g,k}$ not contained in 
${\Acal}_g^{{\rm End}\neq \ZZ}$, then ${\rm codim}_{\Ag} (N)
\geq \binom{k+2}{2}$. Moreover, equality holds if and only if one of the
following happens:

\begin{itemize}
\item [(i)]  $g=k+3$ and $N={\Hcal}_g$;
\item [(ii)] $g=k+4$ and  $N=\Jcal_g$. 
\end{itemize}\end{conjecture}

By work of Beauville \cite {B2} and Debarre \cite {D1} the
conjecture is true for $g=4$ and $g=5$. Debarre  \cite{D0,D1} gave
 examples of components of   $N_{g,k}$ for which the bound in 
Conjecture \ref{conj:mainconj} for the codimension in $\Ag$ fails, 
but they are contained in  $\Ag^{{\rm End}\neq \ZZ}$, since the 
corresponding abelian varieties are isogenous to products.

Our main objective in this paper
will be to prove the conjecture for $k=1$.

\begin{remark} \label{rem:ob}
The question about the dimension of the Andreotti--Mayer loci is related to
the one about the loci $S_{g,k}$ introduced in \S \ref{sec: essgi}. Note
that $M_{g,k}=\pi^*(N_{g,k})$ is a subscheme of $S_{g,k}$. Let $N$ be an
irreducible component of $N_{g,k}$ not contained in $N_{g,k+1}$, let $M$
be the irreducible component of 
$M_{g,k}$ dominating $N$ and $Z$ an irreducible component of $S_{g,k}$
containing $M$. We now give an heuristic argument. Recalling Proposition
\ref{prop:exp}, we can consider $Z$ to be {\it well behaved} if 
${\rm codim}_{S_g}(Z)=\binom{k+1}{2}$. Since $M$ is contained in $Z$, 
one has also  ${\rm codim}_{S_g}(M)\geq \binom{k+1}{2}$ and therefore
${\rm codim}_{\Ag} (N) \geq \binom{k+2}{2}$, which is the first assertion
in Conjecture \ref {conj:mainconj}. In this setting, the equality holds if
and only if $Z$ is {well behaved} and $M=Z$. 

On the other hand, since $M$  is gotten from $Z$
by imposing further restrictions, one could
expect that $M$ is, in general, {\it strictly} contained in $Z$ and therefore that
${\rm codim}_{\Ag} (N) > \binom{k+2}{2}$. 

In this circle of ideas, it is natural to ask if $\Jcal_{g}$ [resp.
$\Hcal_{g}$] is dominated  by a well behaved component of 
$S_{g,g-4}$ [resp. $S_{g,g-3}$]. This is clearly the case if $g=4$. 

A related, but not equivalent, question is whether $\Jcal_{g}$ [resp.
$\Hcal_{g}$] is contained in some irreducible subvariety 
of codimension $c<\binom{g-2}{2}$ [resp. $c<\binom{g-1}{2}$] in $\Ag$, 
whose general point corresponds to a principally polarized
abelian variety $(X,\Theta_X)$ with ${\rm Sing}(\Theta_X)$ containing a
subscheme isomorphic to  ${\rm Spec}(\CC[\epsilon_1,\ldots, \epsilon_k]/
(\epsilon_i\epsilon_j\colon 1 \leq i,j \leq k))$ with $k=g-4$ [resp. $k=g-3$].  

One might be tempted to  believe that an affirmative answer 
to the first question implies a negative answer to the second. 
This is not the case. Indeed  $\Hcal_{4}$ is contained in the 
locus of jacobians of curves with an effective even theta characteristic. 
In this case $S_{4,1}$ has two well behaved irreducible components 
of dimension $8$, one dominating $\Hcal_{4}$ with fibres of dimension $1$,
the other one dominating $\Jcal_{4}\cap \theta_{4,0}$,
where $\theta_{4,0}$ is the \emph{theta--null locus}, see 
\S \ref {sec:AMbound} and \S \ref{g=4} below and \cite {GS}).
These two components intersect along a $7$--dimensional locus in 
$S_4$ which dominates $\Hcal_{4}$. 
 
Note that it is not always the case that a component of  
$S_{g,k}$ is well behaved. 
For example, there is an irreducible component of $S_{g,g-2}$, the one 
dominating ${\Acal}_{g,(1,g-1)}$, which is also an irreducible component 
of $S_g$.
\end{remark}
\end{section}                                                            
%%%%%%%%%%%%%%%%%%%%
%%%%%%%%%%%%%%%%%%%%
\begin{section}{An Example for Genus $g=4$}\label{g=4}
In this section we illustrate the fact that the quadrics associated
to the singularities of the theta divisor may provide more information
than obtained above. The locus $S_4$ in the universal family ${\Xcal}_4$
consists of three irreducible components: i) one dominating $\theta_{0,4}$,
the locus of abelian varieties with a vanishing theta-null; we call it $\Acal$; 
ii) one dominating the Jacobian locus ${\Jcal}_4$; we call it $\Bcal$;
iii) one dominating ${\Acal}_{4,3+1}$, the locus of products of an elliptic
curve with a principally polarized abelian variety of dimension $3$.

The components $\Acal$ and $\Bcal$ have codimension $5$ in ${\Xcal}_4$ and they
intersect along a locus $\Ccal$ of codimension $6$. The image of ${\Ccal}$ in 
${\Acal}_4$ is the locus of Jacobians with a half-canonical $g^1_3$.
The quadric associated to the singular point of order $2$ of $X$ in 
$\Theta_X$ has corank $1$ at a general point. 
We refer to \cite{GS} for a characterization
of the intersection $\theta_{0,4}\cap {\Jcal}_4$.

\begin{proposition} If $\xi=(X,\Theta_X,x)$ is a point of $\Ccal$ then 
${\rm Sing}(\Theta_X)$ contains a scheme of length $3$ at $x$.
\end{proposition}
\begin{proof} The scheme ${\rm Sing}(\Theta_X)$ contains a scheme
${\rm Spec}(\CC[\epsilon])$ at $x$, say corresponding to the tangent
vector $D_1\in \ker q_x$, cf.\ Section \ref{sec:curv}. 
In order that $\xi$ is a singular point of $S_4$ we must have that 
the quadrics $q$ and $q_1$ associated to $x$ and $D_1$ (cf.\ equation
(\ref{eq:a}) and Proposition \ref{prop:sing}) 
are linearly dependent. By Proposition 
\ref{pr:singsing} this implies that ${\rm Sing}(\Theta_X)$ contains 
a scheme of length $3$ at $x$.
\end{proof}

\begin{proposition}
The two components $\theta_{0,4}$ and ${\Jcal}_4$ of $N_{4,0}$ in ${\Acal}_4$
are smooth at a non-empty open subset of their intersection and are
tangent there.
\end{proposition}
\begin{proof} 
Let $\xi=(X,\Theta_X,x)$ be a general point of $\Ccal$ so that the
singular point $x$ of $\Theta_X$ defines a quadric $q_x$ of rank $3$.
Set $\eta=(X,\Theta_X)\in \Acal_4$.
From the equations \eqref{eq:tansp} and \eqref{eq:a}, one deduces
that the tangent cone to $N_{4,0}$ at $\eta$ is supported on $q_x^{\perp}$.
Suppose, as we may by applying a suitable translation, 
that $x$ is the origin $0$ in $X$. 
Then $\theta_{0,4}$ is locally defined by the equation $\theta(\tau,0)$
in $\Acal_4$ (see \cite {GS}, \S 2). 
Hence $\theta_{0,4}$ is smooth at $\eta$ with tangent space
$q_x^{\perp}$. Now  the jacobian locus $\Jcal_4$ is also smooth at $\eta$
by the injectivity of the differential of Torelli's 
morphism at non--hyperelliptic 
curves (see \cite{GH}) and the assertion follows. 
\end{proof}

\end{section}
%%%%%%%%%%%%%%%%%%%%
%%%%%%%%%%%%%%%%%%%%
\begin{section}{The Gauss map and tangencies of theta divisors}\label{sec:gauss}

Let $(B,\Xi)$ be a polarized abelian variety of dimension $g$, where $\Xi$
is an effective divisor on $B$. As usual we let  $B=\CC^g/\Lambda$, with
$\Lambda$ a $2g$--dimensional lattice, and $p: \CC^g\to B$ be the
projection. So we have coordinates $z=(z_1,\ldots,z_g)$ in $\CC^g$ and
therefore on $B$, and we can keep part of the conventions and 
the notation used so far.  
Let $\xi=\xi(z)$ be the Riemann theta function whose divisor on $\CC^g$
descends to $\Xi$ via $p$.

If $\Xi$ is reduced, then the Gauss map of $(B,\Xi)$ is the morphism
$$
\gamma=\gamma_{\Xi} \colon \Xi-{\rm Sing}(\Xi) \to
 \PP(T_{B,0}^{\vee}), \quad
x \mapsto \PP(t_{-x}(T_{\Xi,\xi})),
$$
where $t_{-x}(T_{\Xi,x})$
is the tangent space to $\Xi$ at $\xi$ translated to the origin. If
$\Xi_b=t_{b}(\Xi)$ is the translate of $\Xi$ by the point $b\in B$
defined by the equation $\xi(z-b)=0$
then for $x\in \Xi-{\rm Sing}(\Xi)$
the origin is a smooth point for $\Xi_{-x}$ and
$\gamma(x)=\PP(T_{\Xi_{-x} ,0})$. 

As usual we have natural homogeneous coordinates $(z_1:\ldots :z_g)$ in 
$\PP(T_{B,0})=\PP^{g-1}$ and therefore dual coordinates in
the dual projective space $\PP(T_{B,0}^{\vee})={\PP^ {g-1}}^{\vee}$. 
Then the expression  of $\gamma$ in coordinates equals
$\gamma(p(z))=(\partial_1\xi(z):\ldots :\dd_g\xi(z))$ 
with $\partial_i=\partial/\partial z_i$.

The following lemma is well known. 

\begin{lemma}\label{lem:fin} 
Let $(B,\Xi)$ be a simple abelian variety
of dimension $g$ and  $\Xi$ an irreducible effective divisor on $B$. 
Then the map $\gamma_{\Xi}$ has finite fibres. Moreover, for
a smooth point $x\in \Xi$ there are only finitely many $b\in B$ such that
$\Xi_b$ is smooth at $x$ and tangent to $\Xi$ there.
\end{lemma}

\begin{proof} 
Suppose $\gamma_{\Xi}$ does not have finite fibres. Then there is an 
irreducible curve $C$ of positive geometric genus contained 
in the smooth locus of $\Xi$ which is contracted by $\gamma_\Xi$ to a point
of ${\PP^ {g-1}}^{\vee}$ corresponding to a hyperplane $\Pi\subset \PP^{g-1}$. 
Then the image of the Gauss map $\gamma_C$ of $C$ lies in $\Pi$.
This implies that $C$ does not generate $B$ and contradicts the fact
that $B$ is simple. As to the second statement, if it does not hold
there exists an irreducible curve $C$ such that for all $b \in C$
the divisor $\Xi_b$ is smooth at $x$ and tangent to $\Xi$. Then
the curve $C'=\{ x-b : b \in C\}$ is contained in $\Xi$
and contracted by $\gamma_{\Xi}$. This contradicts the fact that
$\gamma_{\Xi}$ has finite fibres as we just proved.
\end{proof}

\begin{definition}\label{def:Th} 
Let $h$ be a natural number with $1\leq h \leq g-1$.
The subscheme $T_h(B,\Xi)$ of $B\times B^h$ is defined in  $\CC^g \times
(\CC^g)^h$ with coordinates $(z,u_1,\ldots,u_h)$ by the equations 
\begin{equation}\label{eq:loc}
\begin{aligned}
\xi(z)=0, \quad \xi(z-u_1)=0, \quad \ldots \quad, \, \xi(z-u_h)&=0,\cr
{\rm  rk}
\left( \begin{matrix} 
\partial_1\xi(z) & \cdots & \partial_g\xi(z)\cr
\partial_1\xi(z-u_1) & \cdots & \partial_g\xi(z-u_1)\cr
& \vdots & \cr
\partial_1\xi(z-u_h) & \cdots & \partial_g\xi(z-u_h)\cr
\end{matrix}\right) 
&\leq h. 
\end{aligned}
\end{equation}
\end{definition}
The projection to the first factor induces a morphism 
$p_1: T_h(B,\Xi)\to \Xi$. 
Note that for $i=1,\ldots,h$ the variety $E_i$ of codimension 
$h+1\leq g$ in $B\times B^h$ defined by the equations
$$
\xi(z)=0 \quad \xi(z-u_j)=0 \quad j\neq i \quad {\rm and} \quad u_i=0, 
$$
is contained in $T_h(B,\Xi)$. 
Moreover, the expected codimension of 
the irreducible components of $T_h(B,\Xi)$ is $g+1$.  
We will say that an irreducible
component $T$ of $T_h(B,\Xi)$ is {\it regular} if:
\begin{itemize}
\item[(i)] $T\neq E_i$  for $i=1,\ldots,h$;
\item[(ii)] on a non-empty open subset of $T$ all the rows of the matrix 
in \eqref{eq:loc} are non--zero.
 \end{itemize} 
In particular, if $T$ is regular then $p_1(T)\not\subseteq {\rm Sing}(\Xi)$.

\begin{proposition}\label{prop:expcod} 
If $B$ is simple each regular component of $T_h(B,\Xi)$ has dimension
$hg-1$. 
\end{proposition}

\begin{proof}  Let us first assume $h=1$ and denote a component of 
$T_h(B,\Xi)$ by $T$. By composing $p_1$ with the Gauss map $\gamma=\gamma_{\Xi}$, 
we obtain a rational map $\phi:T\map \PP(T_{B,0}^{\vee})$. 
We shall prove that $\phi$ has finite fibres.
Let $v$ be a point in the image of $\phi$ coming from a point
in the open subset as in Definition \ref{def:Th}, (ii). 
By Lemma \ref{lem:fin} 
there are only finitely many smooth points $z_1,\ldots,z_a\in \Xi$ 
such that $\gamma(z_i)=v$ for $i=1,\ldots,a$. For each $1\leq i \leq a$
we consider the theta divisor defined by $\xi(z_i-u)=0$. Again by Lemma 
\ref{lem:fin} there are only finitely many points
$u_{i1},\ldots,u_{i\ell_{i}}$ in it such that $(z_i,u_{ij})$ may 
give rise to a point on $T$.  So $\phi$ has finite fibres. Thus it is
dominant and $T$ has dimension $g-1$. 

Now we prove the assertion by induction on $h$. Consider the projection
$q: T\to B\times B^ {h-1}$ by forgetting the last factor $B$. 
If the image $T'$ of $T$ is contained in
$T_{h-1}(B,\Xi)$, then it is contained in a regular component of
$T_{h-1}(B,\Xi)$, hence by induction the codimension of $T$ is at least $g+1$,
while we know from the equations that it is at most $g+1$, and the assertion
follows. 

Suppose $T'$ is not contained in $T_{h-1}(B,\Xi)$. 
Let $\UU(h-1,g-1)\to \GG(h-1,g-1)$ be the universal bundle over the 
Grassmannian of $(h-1)$--planes in $\PP^{g-1}$. 
Then we define a rational map $\psi: T\map \UU(h-1,g-1)$ 
in the following way. If $(z,u_1,\ldots,u_{h})$ is
a general point in $T$, then by the assumption that 
$T^{\prime}\not\subseteq T_{h-1}(B,\Xi)$ 
the hyperplanes $\gamma(z),\gamma(z-u_1),\ldots,\gamma(z-u_{h-1})$ 
are linearly independent and we define 
$$
\psi(z,u_1,\ldots,u_{h})=(\langle
\gamma(z),\gamma(z-u_1),\ldots,\gamma(z-u_{h-1})\rangle, \gamma(z-u_h)).
$$
We claim that Lemma \ref{lem:fin} implies that the general fibre
of $\psi$ has dimension  $\leq h(h-1)$. Indeed, in a fibre
of $\psi$ the points $\gamma(z), \gamma(z-u_1), \ldots, \gamma(z-u_{h-1})$
vary in a $(h-1)$-dimensional space giving at most $h(h-1)$ parameters 
and if $\gamma(z),\gamma(z-u_1),\ldots,\gamma(z-u_{h})$ are fixed there are for
$(z,u_1,\ldots,u_h)$ only finitely many possibilities.
Thus $\dim(T)$ is bounded from above by
$$
\dim (\UU(h-1,g-1))+h(h-1)=h(g-h)+(h-1)+h(h-1)=hg-1,
$$
and this proves the assertion.
\end{proof}

We will also consider the closed subscheme $T^0_h(B,\Xi)$ of
$T_h(B,\Xi)$ which is defined in  $\CC^g \times (\CC^g)^h$ by the
equations 
$$
\xi(z)=0, \quad\xi(z-u_i)=0, \quad i=1,\ldots,h
$$ 
and
$$
{\rm  rk}\left( \begin{matrix} \partial_1\xi(z) & \cdots & \partial_g\xi(z)\cr
\partial_1\xi(z-u_1) & \cdots & \partial_g\xi(z-u_1)\cr
 & \vdots & \cr
\partial_1\xi(z-u_h) & \cdots & \partial_g\xi(z-u_h)\cr\end{matrix}
\right)=0. 
$$
Finally we will consider the closed subset $\Tcal_h(X,B)$ which is the union
of $T^0_h(X,B)$ and of all regular components of $T_h(B,\Xi)$. Look at
the projection 
$$
p=p_2: \Tcal(X,B)\to B^ h.
$$

\begin{definition}
We define $N_{0,h}(B,\Xi)$ to be the image of $\Tcal_h(X,B)$ under the 
map~$p$. More generally,
for each integer $k$ we define
$$
N_{k,h}(B,\Xi):=\overline{\{ u=(u_1,\ldots,u_h) \in B^h\colon u_1,\ldots,
u_h\neq0, \quad \dim (p_2^{-1}(u)) \geq k \}}. 
$$
\end{definition}
Roughly speaking,
$N_{k,h}(B,\Xi)$ is the closure of the
set of all $(u_1,\ldots,u_h)\in B^ h$ such that 
$\Xi$ contains an irreducible subvariety $V$ of dimension $n\geq k$ and for
all $z \in V$ either:
\begin{itemize}
\item [(a)] $z$, $z-u_1,\ldots,z-u_h$ are smooth points of $\Xi$ and
$\gamma(z), \gamma(z-u_1), \ldots, \gamma(z-u_h)$ are linearly dependent, or
\item  [(b)] all the points $z$, $z-u_1,\ldots,z-u_h$ are contained in ${\rm
Sing}(\Xi)$. 
\end{itemize}
In case (a) we say that the divisors $\Xi, \Xi_{u_1},\ldots,\Xi_{h_h}$,
all passing through $z$ with multiplicity one, are {\it tangentially
degenerate} at $z$. 

In case $h=1$ we have only two divisors which are just
{\it tangent} at $z$. 
We will drop the index $h$ if $h=1$. Thus, $N_0(B,\Xi)$ is the set
of all $u$ such that $\Xi$ and $\Xi_u$ are either tangent, or both
singular, somewhere. Note that, if $\Xi$ is symmetric, i.e., if $\xi(z)$ is
even, as happens for the Riemann theta function, then $N_0(B,\Xi)$
contains the divisorial component $2\Xi:=\{2\xi \colon \xi \in \Xi\}$, since
$\gamma(\xi)=\gamma(-\xi)$ for all smooth points $\xi \in \Xi$. 

One has the following result by Mumford (see \cite{Mumford},
Proposition 3.2).
                               
\begin{theorem} \label{thm:Mumford}
If $(B,\Xi)$ is a principally polarized abelian variety of dimension $g$
with $\Xi$ smooth, then $N_0(B,\Xi)$ is a divisor on $B$
algebraically equivalent to $\frac {(g+2)!} 6 \, \Xi$.
\end{theorem}

As we will see later, the following
result is related to Conjecture \ref{conj:mainconj}. 

\begin{proposition}\label{prop:partial} 
Suppose that $(B,\Xi)$ is a simple principally polarized abelian 
variety of dimension $g$. Assume that
$(B,\Xi)\notin N_{g,k}$. Then for $k\geq 0$ and $1\leq h \leq g-1$
and for every irreducible component $Z$ of $N_{k,h}(B,\Xi)$ 
one has ${\rm codim}_{B^ h} (Z)\geq k+1$. 
\end{proposition}

\begin{proof}  
By the definition of $N_{k,h}(B,\Xi)$ and by the fact that
${\rm Sing}(\Xi)$ has dimension $<k$ an irreducible component 
of $N_{k,h}(B,\Xi)$ can only be contained in the image of a 
regular component of $T_h(B,\Xi)$. The assertion follows by 
Proposition \ref{prop:expcod} and the fact that the fibres of $p_2$ 
on a regular component have dimension $\geq k$.
\end{proof}
\end{section}
%%%%
\begin{section}{Properties of the loci $N_{k,h}(B,\Xi)$}
We will prove now a more precise result in the spirit of
Proposition \ref{prop:partial}.

\begin{proposition}\label{prop:singolare} Suppose that $(B,\Xi)$ is a
simple, principally polarized abelian variety of dimension $g$. Assume that
$N_{k,h}(B,\Xi)$ is positive dimensional for some $h\geq 1$ and $k\geq 1$.
Then $(B,\Xi)\in N_{g,k-1}$.
\end{proposition}

\begin{proof} We may assume $(B,\Xi)\not\in N_{g,k}$, otherwise there is
nothing to prove. We may also assume $h\geq 1$ is minimal under the
hypothesis that $N_{k,h}(B,\Xi)$ has positive dimension.

Let $(u_1,\ldots,u_h)$ be a point in $N_{k,h}(B,\Xi)$, so that $\Xi_0:=\Xi$ and
$\Xi_i:=\Xi_{u_i}, i=1,\ldots,h$, are tangentially degenerate along a
non-empty subset $V^0$ of an irreducible subvariety $V$ of  $B$ of
dimension $k$ such that $\Xi_i$ for $i=0,\ldots,h$ are smooth
at all points $v \in V^0$.

For $j=0,\ldots,h$ the element 
$s_j^{(i)}:=\partial_i\xi(z-u_j)$ (with the convention that $u_0=0$)
is a section of $\Ocal(\Xi_j)$ when restricted to $V$  since $\Xi_j$ contains $V$.
We know that for given $j$ not all $s_j^{(i)}$ are identically zero
on $V^0$.
Our assumptions on the tangential degenerateness and minimality
tell us that there exist non--zero rational functions $a_{j}$ such that
we have a non-trivial relation
 
\begin{equation}\label{eq:relazione}
\sum_{j=0}^h a_{j} s_{j}^{(i)} =0 \quad {\rm  for} \quad  i=1,\ldots,g.
\end{equation}
Suppose that the $a_j$ are regular on all of a desingularization $f:W\to V$
of $V$.  Then they are constant and the relation  $\sum a_{j} s_{j}^{(i)} =0$
holds on the whole of $W$. By writing relation \eqref {eq:relazione} in 
different patches which trivialize the involved line bundles, 
and comparing them, we see that, if the transition functions are not all 
proportional, then we can shorten the original relation by subtracting two of
them. This would contradict the minimality assumption. 
Therefore,  we have that all of the 
divisors $f^*(\Xi -\Xi_{u_j})_{|V}$ with $j=1,\ldots h$ are 
linearly equivalent. 
Since $u=(u_1,\ldots,u_h)$ varies in a subvariety of positive
dimension this implies that the map ${\rm Pic}^0(B) \to {\rm Pic}(W)$
has a positive dimensional kernel, and this is impossible since $B$
is simple.  So there exists an index $j$ such that the function $a_j$
has poles  on  a divisor $Z_j$ of $W$.  
Note now that a point which is non--singular for all the divisors
$\Xi_j$, $j=1,\dots,h$, is certainly not a pole for the functions $a_j$. 
Therefore, for each $j=1,\dots, h$, there an $\ell$ depending of $j$
such that  $Z_j$ is contained in the divisorial part of
the scheme $f^ *({\rm Sing}(\Xi_\ell))$. Moreover, since $Z_j$ moves in a
linear system on $W$, it cannot be contracted by the birational morphism $f$.
This proves that $\Xi_\ell$ is singular along a variety of dimension $k-1$
contained in $V$, which proves the assertion.
\end{proof}

\begin{corollary}\label {cor:partial} Suppose that $(B,\Xi)$ is a
simple, principally polarized abelian variety of dimension $g$. Assume that
$(B,\Xi)\notin N_{g,0}$, that is, $\Xi$ is smooth.
Then every irreducible component $Z$ of
$N_{0,h}(B,\Xi)$ is a divisor of $B^ h$.  \end{corollary}

\begin{proof} 
Each irreducible component $Z$ of $N_{0,h}(B,\Xi)$ is
dominated by a regular component $T$ of $T_h(B,\Xi)$, which has
dimension $hg-1$ by Proposition \ref {prop:expcod}. The map $p: T\to Z$ is
finite by Corollary  \ref{lem:fin}. The assertion follows.
\end{proof}

\begin{remark}\label{rem:piu} 
If we have two divisors $\Xi_0,\Xi_1$ which are tangentially
degenerate along an irreducible $k$--dimensional variety $V$ whose general
point is smooth for both $\Xi_0,\Xi_1$, then $\Xi_0,\Xi_1$ are {\it both}
singular along some $(k-1)$--dimensional variety contained in $V$.
This can be easily proved by looking at the relation \eqref {eq:relazione}
in this case, and noting that the polar divisors $Z_j$ is contained
in $f^*({\rm Sing}(\Xi_j))$, $j=0,1$.
\end{remark}

\begin{remark}\label{rem:casi} 
Suppose $(B,\Xi)\notin N_{g,0}$. Then
$N_0(B,\Xi)$ is described by all differences of pairs of points of $\Xi$
having the same Gauss image.

Suppose $(B,\Xi)\in N_{g,0}-N_{g,1}$ and assume $\{x,-x\}={\rm
Sing}(\Xi)$ have multiplicity $2$ and the quadric $Q_x=Q_{-x}$ is smooth.
It may or may not be the case that  $b=-2x\in N_1(B,\Xi)$. In any case,
we claim that
$N_1(B,\Xi)-\{-2x\}$ is contained in the set of all differences of
points in $\gamma_\Xi^ {-1}(Q_x^ *)$ with $x$. Let us give a sketch of
this assertion, which provides, in this case, a different argument for the
proof of Proposition \ref{prop:partial}.

If $b\in N_1(B,\Xi)-\{-2x\}$, there is a curve $C\subset \Xi$ such that
for $t\in C$ general $\gamma_\Xi(t)=\gamma_\Xi(t+b)$. Along $C$ the
divisors $\Xi$ and $\Xi_b$ are tangent, hence the curve contains $x$ (see
Remark \ref{rem:piu}). Note that the curve $C$ is smooth at
$x$. Indeed, locally at $x$, the divisor $\Xi$ is a quadric cone of rank
$g-1$ in $\CC^ g$ with vertex $x$, whereas $\Xi_b$ is a hyperplane through
$x$, and they can only be tangent along a line.

Thus it makes sense to consider the image of $x$ for
$\gamma_C$, which is  a point on $Q_x$. The point $x+b\in
C+b\subset \Xi$ is smooth for $\Xi$ and $\gamma(x+b)$ is clearly tangent
to $Q_x$ at $\gamma_C(x)$. 
\end{remark}

\end{section}
\begin{section}{On $N_{g-2}$ and $N_{g-3}$ for Jacobians}
The following result shows that the bound in Proposition \ref
{prop:partial} is sharp in the case $h=1$.
Recall that a curve $C$ is called \emph{bielliptic} if it is a double cover of
an elliptic curve. 

\begin{proposition}\label{prop:jac}  
Let $C$ be a smooth, ireducible projective curve of genus $g$ and let 
$(J=J(C),\Theta_C)$ be its principally polarized jacobian. Then one has:
\begin{itemize}
\item[(i)]
$\{\Ocal_C(p-q) \in J:  p,q\in C\} \subseteq N_{g-3}(J,\Theta_C)$;
\item [(ii)] either $C$ is bielliptic or the
equality holds in (i);
\item [(iii)] if $C$ is hyperelliptic, then 
$$
\{\Ocal_C(p-q) \in J: 
p,q\in C, h^ 0(C,\Ocal_C(p+q))=2\}= N_{g-2}(J,\Theta_C).
$$
\end{itemize}
\end{proposition}

\begin{proof} We begin with (i). We assume that $C$ is not
hyperelliptic; the hyperelliptic case is similar and can be left to
the reader. We may identify $C$ with its canonical image in $\PP^{g-1}$. 
Moreover, we identify $J$ with ${\rm Pic}^{g-1}(C)$ and 
$\Theta_C\subset {\rm Pic}^{g-1}(C)$
with the set of effective divisor classes of degree $g-1$. 
Then the Gauss map $\gamma_C:=\gamma_{\Theta_C}$ can be geometrically 
described as the map which sends the class of an effective divisor $D$ 
of degree $g-1$ such that $h^ 0(C,\Ocal_C(D))=1$ to the hyperplane in 
$\PP^{g-1}$ spanned by $D$ (see \cite {GH}, p. 360). 

Take two distinct points $p_1,p_2$ on $C$. Then
$|\omega_C(-p_1-p_2)|$ is a linear series of degree $2g-4$ and dimension
$g-3$. For $i=1,2$ we let $V_i$ be the subvariety of 
$\Theta_C$ which is the Zariski closure of the set of all divisor
classes of the type $E+p_i$, where $E$ is
a divisor of degree $g-2$  contained in some divisor of
$|\omega_C(-p_1-p_2)|$. Clearly $\dim(V_i)=g-3$, hence $V_i$ is not contained
in ${\rm Sing}(\Theta_C)$ which is of dimension $g-4$.
If $u$ denotes the divisor class $p_2-p_1$ then $x\mapsto x+u$ defines an
isomorphism $V_1 \cong V_2$. Moreover, for $x$ in an non-empty open subset of 
$V_1$ we have $\gamma_C(x)=\gamma_C(x+u)$. This proves (i).

Conversely, assume there is a point $u\in {\rm Pic}^ 0(C)-\{0\}$, and a
pair of irreducible subvarieties $V,V'$ of $\Theta_C$ of dimension $g-3$
such that $x \mapsto x+u$ gives a birational map from $V_1$ to $V_2$.
For $x$ in a non-empty open subset $U \subset V_1$ 
we have $\gamma_C(x)=\gamma_C(x+u)$. If $D$ and $D^{\prime}$ are 
the effective divisors of degree $g-1$ on $C$ corresponding to 
$x$ and $x+u$  and if $E$ is the least common divisor of $D$ and $D'$
then by the geometric interpretation of the Gauss map $\gamma_C$  there 
is an effective divisor $F$ with ${\rm deg}(E)={\rm deg}(F)$ such that 
$D+D^{\prime}-E+F\equiv K_C$. Thus $(2D-E)+F\equiv K_C-u$. 

Consider the linear series $|K_C-u|$,
which is a $g_{2g-2}^ {g-2}$. If this linear series has a base point $q$,
then there is a point $p$ such that $K_C-u-q\equiv K_C-p$, i.e.
$D'-D\equiv u\equiv p-q$, proving (ii). So we may assume that
$|K_C-u|$ has no base point. If $C$ is not
bielliptic, then  $|K_C-u|$ determines a birational map of $C$ to a curve
in $\PP^{g-2}$. On the other hand it contains the $(g-3)$--dimensional
family of divisors of the form $(2D-E)+F$, which are singular along the divisor
$D-E$.  This is only possible if ${\rm deg}(D-E)=1$, i.e., $u\equiv
D'-D\equiv p-q$, with $p,q\in C$. But in this case $|K_C-u|$ has the base
point $q$, a contradiction. This proves (ii).

Assume now $C$ is hyperelliptic.  Let 
$p_1+p_2$ be an effective divisor of the $g^1_2$ on $C$
with $p_1\neq p_2$. Then  $|\omega_C(-p_1-p_2)|$ is a linear series of 
degree $2g-4$ and dimension $g-2$.  For $i=1,2$ we let $V_i$ 
be the subvariety of $\Theta_C$ which is the Zariski closure of the 
set of all divisors classes of the type $E+p_i$, where $E$ is a divisor 
of degree $g-2$ contained in a divisor of
$|\omega_C(-p_1-p_2)|$. The variety $V_i$ 
has dimension $g-2$ and is not contained in ${\rm Sing}(\Theta_C)$,
which is of dimension $g-3$. The translation over $u$ induces an
isomorphism $V_1 \cong V_2$ and for $x$ in a non-empty subset $U$ of $V_1$ we
have $\gamma_C(x)=\gamma_C(x+u)$. Hence the left hand side in
(iii) is contained in $N_{g-2}(J,\Theta_C)$.

Finally, assume there is a point $u\in {\rm Pic}^0(C)-\{0\}$, and a pair
of irreducible subvarieties $V$, $V'$ of $\Theta_C$ of dimension $g-2$ such
that translation by $u$ gives a rational map $V \map V^{\prime}$  with
$\gamma_C(x)=\gamma_C(x+u)$ on a non-empty open subset of $V$. 
Let $D$, $D'$ be the effective divisors of degree $g-1$ on $C$ 
corresponding to $x$ and $x+u$. As in the proof of
part (ii), let $E$ be the least common divisor of $D$ and $D'$. 
In the present situation the linear series $|K_C-u|$ of dimension $g-2$
contains the $(g-2)$--dimensional family of divisors of the form $2D-E$, 
which are singular along the divisor $D-E$. This means that $2(D-E)$ is in 
the base locus of $|K_C-u|$. This is only possible if $D=E+q$ for 
some point $q\in C$, and $K_C-u-2q\equiv (g-2)(p+q)$, 
where $p$ is conjugated to $q$ under the
hyperelliptic involution. In conclusion, we have $u\equiv p-q$ and the
equality in (iii) follows.
\end{proof}

\begin{remark}\label {rem:biell} The hypothesis that $C$ is not bielliptic
is essential in (ii) of Proposition \ref {prop:jac}. Let in fact $C$ be a
non-hyperelliptic bielliptic curve which is
canonically embedded in $\PP^ {g-1}$. Let $f: C\to E$ be the bielliptic
covering. 
One has $f_*\Ocal_C\simeq \Ocal_E\oplus \xi^{\vee}$, with 
$\xi^{\oplus 2}\simeq \Ocal_E(B)$, where $B$ is the branch divisor of $f$.

Let $u\in {\rm Pic}^ 0(E)-\{0\}$ be a general point, which we can 
consider
as a non--trivial element in ${\rm Pic}^ 0(C)$ via the inclusion $f^ *: {\rm
Pic}^ 0(E)\to {\rm Pic}^ 0(C)$ 
Note that $f$ is ramified, hence $f^*$ is injective. 
We want to show that $u\in N_{g-3}(J,\Theta_C)$, 
proving that equality does not hold in (i) in this case. 

The canonical image of $C$ is contained in a cone $X$
with vertex a point $v\in \PP^ {g-1}$ over the curve $E$ embedded in $\PP^
{g-2}$ as a curve of degree $g-1$ via the linear system $|\xi|$. 
Let us consider the
subvariety $W$ of  $|\xi|$ consisting of all divisors $M\in
|\xi|$ such that there is a subdivisor $p+q$ of $M$ with $p,q\in E$
and $p-q\equiv u$.  It is easily seen that $W$ is irreducible of dimension
$g-3$. 

Notice that for $M$ general in $W$ one has $M=p+q+N$, with $N$ effective
of degree $g-3$. Therefore we may write $K_C\equiv f^ *(M)\equiv f^
*(p)+f^ *(q)+F+F'$, where $F, F'$ are disjoint, effective divisors of degree
$g-3$ which are exchanged by the bielliptic involution. 

We let $V$ be the
$(g-3)$--dimensional subvariety of $\Theta_C$ described by all
classes of divisors $D$ of degree $g-1$ on $C$ of the form $D=f^ *(p)+F$,
as $M$ varies in $W$.  Then $D+u\equiv D':=f^*(q)+F$ and $D$ and $D'$ span
the same hyperplane through $v$ in $\PP^ {g-1}$. Therefore, if $x\in V$
is the point corresponding to $D$, one has $\gamma_C(x)=\gamma_C(x+u)$.
This proves that $u\in
N_{g-3}(J,\Theta_C)$.
\end{remark}
\end{section}

%%%
%%%
\begin{section}{A boundary version of the Conjecture}
\label{boundaryversion}
We will now formulate a conjecture.  As we will see later,
it can be considered as a  {\sl boundary version} of Conjecture
\ref{conj:mainconj} (see also Proposition
\ref{prop:partial}).

\begin{conjecture}\label{conj:bound} 
Suppose that $(B,\Xi)$ is a simple principally polarized abelian 
variety of dimension $g$. Assume that $(B,\Xi)\notin N_{g,i}$ for 
all $i\geq k\geq 1$. Then there is an irreducible component $Z$ of 
$N_k(B,\Xi)$ with  ${\rm codim}_B(Z)= k+1$ if and only if one of 
the following happens:
\begin{itemize}
\item [(i)] either $g\geq 2$, $k=g-2$ and $B$ is
a hyperelliptic Jacobian,
\item  [(ii)] or $g\geq 3$, $k=g-3$ and $B$ is a Jacobian.
\end{itemize}
\end{conjecture}

One implication in this conjecture holds by Proposition \ref{prop:jac}.
Note that the conjecture would give an answer for simple abelian varieties
to the Schottky problem that asks for a characterization of Jacobian
varieties among all principally polarized abelian varieties.
For related interesting questions, see \cite {PP}.
%%%
%%%
\end{section}
\begin{section}{Semi-abelian Varieties of Torus Rank One}\label{sec:rank1}
Let $(B,\Xi)$ be a principally polarized abelian variety of dimension $g-1$.
The polarization $\Xi$ gives rise to the isomorphism
$$
\phi_{\Xi}: B\to \hat{B}={\rm Pic}^0(B), \, b \to {\Ocal}_B(\Xi-\Xi_b).
$$
and we shall identify $B$ and $\hat{B}$ via this isomorphism.
Thus an element $b \in \hat{B}\cong B$ determines a line bundle
$L=L_b={\Ocal}_B(\Xi-\Xi_b)$ with trivial first Chern class. We can associate
to $L$ a \emph{semi-abelian variety} $X=X_{B,b}$, namely the $\GG_m$-bundle
over $B$ defined by $L$ which is an algebraic group since it coincides
with the theta group ${\mathcal G}_b:={\mathcal G}(L)$ of $L$ 
(cf.\ \cite{MumfordAV}, p.\ 221). 
This gives the well-known equivalence between $\hat{B}$ and 
${\rm Ext}(B,\GG_m)$, the group of
extension classes of $B$ with $\GG_m$ in the category of algebraic groups
(see \cite{serre}, p.\ 184).

Both the line bundle $L$ and the $\GG_m$-bundle
$X$ determine a $\PP^1$-bundle $\PP=\PP(L\oplus {\Ocal}_B)$ over $B$ with 
projection $\pi: \PP \to B$ and two sections
over $B$, say $s_0$ and $s_{\infty}$ given by the projections 
$L\oplus {\Ocal}_B \to L$ and $L\oplus {\Ocal}_B \to {\Ocal}_B$. 
If we set $P_0=s_0(B)$ and 
$P_{\infty}=s_{\infty}(B)$ then we can identify $\PP-P_0-P_{\infty}$
with $X$. By \cite{H}, Proposition \ 2.6 on page 371 we have ${\Ocal}_{\PP}(1)
\cong {\Ocal}_{\PP}(P_0)$.
We can complete $X$ by considering the non-normal
variety $\overline{X}=\overline{X}_{B,b}$ obtained
by glueing $P_0$ and $P_{\infty}$
by a translation over $b\in B\cong \hat{B}$. On $\PP$ we have the
linear equivalence $P_0-P_{\infty}\equiv \pi^{-1}(\Xi-\Xi_b)$. 
We set $E:=P_{\infty}+
\pi^{-1}(\Xi)$ and put $M_{\PP}={\Ocal}_{\PP}(E)$. 
This line bundle restricts
to ${\Ocal}_B(\Xi)$ on $P_0$ and to ${\Ocal}_B(\Xi_b)$ on $P_{\infty}$, 
and thus descends
to a line bundle $M=M_{\overline{X}}$ on $\overline{X}$. 
We have $\pi_*({\mathcal O}_{\PP}(E))={\mathcal O}_B(\Xi)\oplus
{\mathcal O}_B(\Xi_b)$ and $H^0(\PP,M_{\PP})$ is generated
by two sections with divisors $P_{\infty}+\pi^{-1}(\Xi)$ and
$P_{0}+\pi^{-1}(\Xi_b)$. One concludes that $H^0(\overline{X},M)$
corresponds to the sections of $M_{\PP}$ such that translation over $b$
carries its restriction to $P_0$ to the restriction to $P_{\infty}$.
It follows that $h^0(\overline{X},M)=1$ with effective divisor 
$\overline{\Xi}$, which is called the {\sl generalized theta divisor} on 
$\overline{X}$.

Analytically we can describe a section of $\Ocal_{\bar X}(\overline \Xi)$ 
on the universal cover $\CC \times \CC^{g-1}$ by a function
$$
\xi(\tau,z)+u\, \xi(\tau,z-\omega),
$$
where $\omega \in \CC^{g-1}$ represents
$b \in B=\CC^{g-1}/\ZZ^{g-1}+\tau \ZZ^{g-1}$,
$\xi(\tau,z)$ is Riemann's theta function for $B$ and $u=\exp(2\pi i \zeta)$
is the coordinate on $\CC^*$. This is called the 
{\it generalized theta function} of $\overline X$.

Let $D$ be the Weil divisor on $\overline{X}$ that is the image of $P_0$
(or, that is the same, of $P_{\infty}$).
We consider a locally free subsheaf $T_{\rm vert}$ of the tangent sheaf to $\overline{X}$,
namely the dual of the sheaf $\Omega^1(\log D)$ of rank $g$.
If $d \in D$ is a point such that on the normalization $\PP$ near
the two preimages $z_1,\ldots,z_{g-1},u$ and $z_1+b,\ldots,z_{g-1}+b,v$
are local coordinates such that $u=0$ defines $P_0$ (resp.\ $v=0$ defines
$P_{\infty}$) with $uv=1$ then $T_{\rm vert}$ is generated by
$\partial/\partial z_1,\ldots,\partial/\partial z_{g-1}, u\partial/\partial u
-v\partial/\partial v$. Here $z_1,\ldots,z_{g-1}$ are coordinates on $B$.
We interpret local sections of $T_{\rm vert}$ as derivations.
In particular, if an effective Cartier divisor $Y$ of $\overline X$ 
has local equation $f=0$, then for each local section $\partial$ of 
$T_{\rm vert}$, the restriction to $Y$ of
$\partial f$ is a local section of 
$\Ocal_{\overline X}(Y)/\Ocal_{\overline X}$.
Then the subscheme ${\rm Sing}_{\rm vert}(\Xi)$ of 
${\overline{\Xi}}$
is locally defined by the $g$ equations
\begin{equation}\label{singvertdef}
\partial_i f=0  \quad 
\hbox{\rm modulo $f$ in $O_{\overline{X}}(\overline{\Xi})/O_{\overline{X}}$}
\end{equation}
with $f=0$ a local equation of $\overline{\Xi}$ and $\partial_i$ local
generators of $T_{\rm vert}$.

The equations for ${\rm Sing}_{\rm vert}(\overline \Xi)$ on $X$ 
are thus given by
$$
\begin{aligned}
\xi(\tau,z)+u\, \xi(\tau,z-\omega)=& 0, \\
\xi(\t,z-\omega)=& 0, \\
\partial_i\xi(\tau,z)+u \partial_i \xi(\tau,z-\omega)=& 0,
\quad (1\leq i \leq g-1). \\ \notag
\end{aligned}
$$
The points in ${\rm Sing}_{\rm vert}(\overline{\Xi})$ are of two sorts
depending on whether they lie on the double locus $D$ of $\overline{X}$
or not. 
The singular points of ${\rm Sing}(\overline \Xi)$ 
on $X=\overline X-D$ are the points $(z,u)$, with $u\neq 0$,
which are zeros of $\xi(\t,z)$ and $\xi(\t,z-\omega)$ and such that
$\gamma(\tau,z)=-u\gamma(\t,z-\omega)$. That is, geometrically, these
correspond under the projection on $B$ to
the points on $B$ where $\Xi$ and $\Xi_b$ are tangent to each other.
To describe the singular points of ${\rm Sing}(\overline \Xi)$ on $D$, 
we consider the composition
$$
\phi \colon B \cong P_0 \to \PP {\buildrel \nu \over \longrightarrow}
\overline{X},
$$
where $\nu$ is the normalisation. Then we have
$ \phi^{-1}({\rm Sing}_{\rm vert}(\overline \Xi))= {\rm Sing}(\Xi)$
and the same if we identify $B$ with $P_\infty$.

Points of ${\rm Sing}_{\rm vert}(\overline{\Xi})$ determine again quadrics
in $\PP^{g-1}$ as follows. Note that the projective space $\PP(T_{X,0})$
contains a point
$P_b$ corresponding to the tangent space $T_{\GG_m}\subset T_X$
of the algebraic torus $\GG_m$ at the origin.
Recall that we write $\gamma(\t,z)$ for the row vector
$$
\gamma(\t,z)=(\dd_1\xi,\ldots,\dd_{g-1}\xi)(\t,z).
$$
Then a singular point determines a, possibly indeterminate, quadric 
defined by the matrix
\begin{equation}\label{eq:matrice}
\left(
\begin{matrix}
0 & \gamma(\t,z-\omega) \\
\gamma(\t,z-\omega)^t & M \\
\end{matrix}
\right)
\end{equation}
with $M$ the $(g-1)\times (g-1)$ matrix $(\dd/\dd \t_{ij} \xi(\t,z)+
 u\dd/\dd \t_{ij} \xi(\t,z-\omega))$. Note that we have
$\gamma(\t,z)=-u\gamma(\t,z-\omega)$.
The quadric passes through the point $P_b$.
For a point on $D$ the quadric is
a cone with vertex $P_b$ over a quadric in $\PP^{g-2}$ given by~$M$.

\begin{remark}\label{rem:triplepts} 
{\rm The above considerations show that a point in 
${\rm Sing}_{\rm vert}(\overline{\Xi})$ has to be regarded as a
point of multiplicity larger than $2$ if the matrix (\ref{eq:matrice})
vanishes identically. This can happen only if $z$ and $z+\omega$ 
are both singular for $\Xi$.}
\end{remark}

\end{section}
%%%
%%%
\begin{section}{Standard Compactifications of Semi-Abelian Varieties}\label
{sec:stand}
%%%
%%%
Let $(B,\Xi)$ be a principally polarized abelian variety.
We assume now that $\dim(B)=g-r$ with $r\geq 1$ 
and extend the considerations of the previous section.

The extensions of $B$ by $\GG_m^r$ are parametrized by ${\rm Ext}^1(B,\GG_m^r)
\cong \hat{B}^r$. To a point $(b_1,\ldots,b_r) \in \hat{B}^r$
one associates the $\GG_m^r$-extension $X=X_b$ obtained
as the fibre product of theta groups
${\mathcal G}_{b_1}\times_B \cdots \times_B {\mathcal G}_{b_r}$.

One of the type of degenerations of abelian varieties 
that we shall encounter are special compactifications
of semi-abelian varieties. We shall call them 
{\sl standard compactifications of torus rank $r$}.
Let $b=(b_1,\ldots,b_r)\in \hat{B}^r$. The algebraic group $X=X_b$
sits in a $\PP^1\times_B \cdots \times_B \PP^1$-bundle $\pi: \PP\to B$ 
that is obtained as the fibre product over $B$ of the 
$\PP^1$-bundles $P_{b_i}=\PP(L_{b_i} \oplus {\Ocal}_B)$.
The complement $\PP-X$ is a union of $2r$ divisors $\sum_{i=1}^r \Pi_0^{(i)}
+\Pi_{\infty}^{(i)}$, where $\Pi_0^{(i)}$ (resp.\ $\Pi_{\infty}^{(i)}$)
is given by taking $0$ (resp.\ $\infty$) in the $i$-th fibre coordinate, with projections
$\pi_{i,0}, \pi_{i,\infty}$ to $B$.

We now define a non-normal variety obtained from $\PP$
by glueing $\Pi_0^{(i)}$ with 
$\Pi_{\infty}^{(i)}$ for $i=1,\ldots,r$.
This identification depends on a 
$r\times r$-matrix $T=(t_{ij})$
with entries from $\GG_m$ such that $t_{ii}=1$ and $t_{ij}=t_{ji}^{-1}$.
Let $s_0^{(i)}: B \to P_{b_i}$ (resp.\ $s_{\infty}^{(i)}$)
be the zero-section (infinity section)
of $P_{b_i}$. We glue the point 
$$
(\beta,x_1,\ldots,x_{i-1},s_0^{(i)}(\beta),x_{i+1},\ldots, x_r)
$$
on $\Pi_0^{(i)}$ with the point
$$
(\beta+b_i,t_{i,1}x_1,\ldots,t_{i,i-1}x_{i-1},
s_{\infty}^{(i)}(\beta),t_{i,i+1} x_{i+1},\ldots, t_{i,r}x_r)
$$
on $\Pi_{\infty}^{(i)}$. We denote the resulting variety by
$\overline{X}$. It depends on the parameters $b\in \hat{B}^r$ and
$t \in {\rm Mat}(r\times r, \GG_m)$. 

We have the linear equivalences $\Pi_{0}^{(i)}-\Pi_{\infty}^{(i)}
\equiv \pi^{*}(\Xi-\Xi_{b_i})$.
We set $E=\Pi_{\infty}+\pi^{*}(\Xi)=\sum \Pi_{\infty}^{(i)}
+\pi^*(\Xi)$ and $E_i=\Pi_{\infty}-\Pi^{(i)}_{\infty}$
and $M:=M_{\PP}=\Ocal_{\PP}(E)$. This line bundle restricts to
$\Ocal_{\Pi_{0}^{(i)} } (E_i+\pi_{i,0}^ *(\Xi))$ on $\Pi_{0}^{(i)}$ and to
$\Ocal_{\Pi_{\infty}^{(i)}} (E_i+\pi_{i,\infty}^ *(\Xi_{b_i}))$ on
$\Pi_{\infty}^{(i)}$. Thus, by the definition of the glueing, $M$ descends to a
line bundle ${\overline M}:=M_{\overline{X}}$ on $\overline{X}$.
We have
$$
\begin{aligned}
\pi_* (M)=& \big (\otimes_{i=1}^ r
(\Ocal_B\oplus L_i^{-1})\big)\otimes \Ocal_B(\Xi) \cr
\cong & \oplus_{k=1}^ r \, \big( \oplus_{0\leq i_1<...<i_k\leq r}
\,  \Ocal_B(\Xi_{b_{i_1}+...+b_{i_k}})\big).\cr
\end{aligned}
$$
Hence we have $h^0(\PP,M)=2^r$. As in the preceding section one sees
that only a $1$-dimensional space of sections descends to sections
of $\overline{M}$ on $\overline{X}$.
In terms of coordinates $(\zeta_1,\ldots,\zeta_r,z_1,\ldots,z_{g-r})$
on the universal cover of $X$, where $(z_1,\ldots,z_{g-r})\in \CC^ {g-r}$ 
are coordinates on the universal cover of $B$, 
a non-zero section of $\overline{M}$ is given by
$$
\sum_{I\subseteq \{1,\ldots,r\}} u_I\, t_I \, \xi(\tau,z-\omega_I),
$$
where $I$ runs through the subsets of $\{1,\ldots,r\}$,
$u_I=\prod_{i\in I} u_i$ with $u_i=\exp(2\pi \zeta_i)$,
$t_I=\prod_{i,j \in I, i<j} t_{ij}$, $b_I=\sum_{i\in I} b_i$ and 
$\omega_I\in \CC^ {g-r}$ represents $b_I\in B$. 
This is the {\sl generalized theta function} of $\overline X$, 
whose zero locus is the {\sl generalized theta divisor} 
$\overline \Xi$ of $\overline X$.

Next we look at the singular points of $\overline \Xi$. 
All points in $\overline \Xi\cap D$ are singular points of
$\overline \Xi$. However, just as in the rank one case in the
preceding section we will in general disregard these singularities 
of $\overline \Xi$, and we will only look at the so called 
{\it vertical singularities}, which we  are going to define now
(cf.\ \cite{Mumford}, \S 2).

The locally free subsheaf $T_{\rm vert}$ of rank $g$ of the tangent 
sheaf $T_{{\overline X}}$ is the dual of $\Omega^1(\log D)$.
Its pull back to $\PP$ is generated, in the $(u,z)$--coordinates, 
by the differential operators 
$u_i\partial/\partial u_i - v_i \partial/\partial v_i$
with $u_iv_i=1$
for $i=1,\ldots,r$ and $\partial_j=\partial/\partial z_j$
with  $j=1,\ldots,g-r$.
We interpret local sections of $T_{\rm vert}$ as derivations as above
and define the scheme ${\rm Sing}_{\rm vert}({\overline \Xi})$ of 
{\it vertical singular points} of $\overline \Xi$ as the subscheme 
of ${\overline \Xi}$ defined by the equations (\ref{singvertdef})
with $f=0$ a local equation of $\overline \Xi$ for
all local sections $\partial \in T_{\rm vert}$. 
This is independent of the choice of a local equation.

\begin{lemma}\label{stand}
Let $(\overline{X},\overline{\Xi})$ 
be a standard compactification of a semi-abelian variety
$X$ of torus rank $r$ with abelian part $(B,\Xi)$. If $\dim({\rm Sing}_{\rm
vert} (\overline{\Xi}))\geq 1$ then $(B,\Xi) \in N_{g-r,0}$.
\end{lemma}
\begin{proof}
The compactification $\overline{X}$ is a stratified space and the 
(closed) strata are (standard) compactifications of semi-abelian extensions 
of $(B,\Xi)$ of torus rank $s$ with $0 \leq s \leq r$. In view of 
the relations $\sum_{I} u_I t_I \partial_i \xi(z-\omega_I)=0$ the vertical
singularities of $\overline{\Xi}$ correspond to points where
$2^s$ translates of $\Xi$ are tangentially degenerate.
Therefore we have $\dim(N_{1,h}(B,\Xi)) \geq 1$ 
with $h=2^s$.
By Lemma \ref{prop:singolare} it follows that $(B,\Xi) \in N_{g-r,0}$.
\end{proof}
\end{section}
\begin{section}{Semi-abelian varieties of torus rank two}\label{sec:rank2} 

In the compactification of the moduli space of principally polarized abelian
varieties of dimension $g$ we shall encounter two types of degenerations
of torus rank $2$. The first of these is a standard compactification 
introduced above and its normalization is a $\PP^1\times \PP^1$-bundle
over a principally polarized abelian variety of dimension $g-2$.
For such a standard compactification the equations for 
${\rm Sing}_{\rm vert}(\Xi)$ are given in terms of the $(u,z)$-coordinates
by the system (we write $t$ instead of $t_{1,2}$; note that $t\neq 0$)
$$
\begin{aligned}
\xi(z)-tu_1u_2\xi(z-\omega_1-\omega_2)=&0,\\
u_1\xi(z-\omega_1)+tu_1u_2\xi(z-\omega_1-\omega_2)=&0, \\ 
u_2\xi(z-\omega_2)+tu_1u_2\xi(z-\omega_1-\omega_2)=&0,\\
\partial_i\xi(z)+u_1 \partial_i\xi(z-\omega_1)+u_2\partial_i\xi(z-\omega_2)+
tu_1u_2\partial_i \xi(z-\omega_1-\omega_2)=&0,\\
i=1,\ldots, g-2. \\
\end{aligned}
$$
From this and the analogous equations in the $v$-coordinates (with $u_iv_i=1$)
we see that the vertical singular points of $\overline \Xi$ are
essentially of three types:

(i) A point $x\in D$ that is the image via
$\varphi:X\to \overline X$ of a point in $\Pi_{1,0}\cap\Pi_{2,0}\cong B$,
(i.e.\ in the $(u,z)$ coordinates one has $u_1=u_2=0$) is a vertical
singular point of $\overline \Xi$ if and only if it corresponds 
to a singular point of $\Xi$ on $\Pi_{1,0}\cap\Pi_{2,0}\cong B$ and to a
singular point of $\Xi_{b_1+b_2}$ on $\Pi_{1,\infty}\cap\Pi_{2,\infty}\cong
B$.

(ii) A point $x\in D$ which is the image via $\varphi$ of
a point in $\Pi_{j,0}$ but not of a point in $\Pi_{3-j,0}$, 
(i.e., in the $(u,z)$ coordinates one has $u_i=0,u_{3-j}\neq 0$, for
a $j=1,2$) is a vertical singular point of
$\overline \Xi$ if and only if 
$$
\begin{aligned}
\xi(\t,z)=0, \quad
\xi(\t,z-\omega_{3-j})=&0,\\
\partial_i\xi(\tau,z)+
u_{3-j}\partial_i\xi(\t,z-\omega_{3-j})=&0,\quad i=1,\ldots, g-2. \\
\end{aligned}
$$
i.e.\ if and only if $z$ and $z-b$ are in $\Xi$ and
$\gamma_{\Xi}(z)=\gamma_{\Xi}(z-b)$.

(iii) A point $x\notin D$, (i.e.\ in the $(u,z)$
coordinates one has $u_1\neq 0\neq u_2$) 
is a vertical singularity if and only if $z$ 
is a singular point of
the divisor $H\in |2\Xi_{b_1+b_2}|$ defined by the equation
$$
\xi(\tau,z-\omega_1)\xi(\tau,z-\omega_2)=t\, \xi(\tau,z)
\xi(\t,z-\omega_1-\omega_2).
$$
By the way, this occurs even in case (ii) above.  
Note also that, by the above equations, 
the existence of a vertical singularity
implies that the theta divisors
$\Xi$, $\Xi_{b_1}$, $\Xi_{b_2}$ and $\Xi_{b_1+b_2}$ are
tangentially degenerate at some point $x$ of $B$, i.e.,
$z \in N_{0,3}(B,\Xi)$.

We call {\it of type} (i), (ii) or (iii)
the singular points of $\overline \Xi$ according to whether cases (i), (ii)
or case (iii) occur. 

\begin{remark}\label {rem:coni2} A point $x$ in ${\rm Sing}_{\rm vert}({\overline
\Xi})$ again determines a quadric $Q_x$ in $\PP^{g-1}$.  It is useful to
remark that:

(a) in case (i) the quadric $Q_x$ is a cone with vertex the line
$L_{\bf b}:=\langle P_{b_1},P_{b_2}\rangle$ 
given by the tangent space to the toric
part over the quadric $Q_z$ in $\PP^ {g-3}$ which corresponds to the 
singular point $z$ of $\Xi$;

(b) in case (ii), say we are at a point with $u_1=0, u_2\neq 0$.
Then $Q_x$ is a cone with vertex $P_{b_1}$ over the quadric in the
hyperplane $u_1=0$ with matrix  
$$
\left(
\begin{matrix}
0 & -\gamma(\t,z)^ t \\
-\gamma(\t,z) & M \\
\end{matrix}
\right)
$$
with 
$
\gamma(\t,z)=(\dd_1\xi,\ldots,\dd_{g-2}\xi)(\t,z)
$
and the matrix $M$ is given by
$$M=(\dd/\dd \t_{ij} \xi(\t,z)+ u\dd/\dd \t_{ij}
\xi(\t,z-\omega))_{1\leq i,j\leq g-2}.$$
\end{remark}

In  \S \ref {sec:rank1} we saw that all rank $1$ compactifications 
of ${\GG}_m$--extensions of a principally polarized abelian variety $B$
form a compact family $\hat{B}$. This is
no longer the case in the higher rank case. This is where semi--abelian
varieties of {\sl non--standard type} come into the picture. This will
depend of choices. It is good to see this in some detail in the rank 2 case.

Given a principally polarized abelian variety $(B,\Xi)$ of dimension
$g-2$, all standard rank 2 compactifications of $(B,\Xi)$ are of the form
$({\overline X},{\overline \Xi})$ with ${\overline X}={\overline
X}_{B,b,t}$ with $b=(b_1,b_2)\in B\times B$ and $t\in \CC^*$. Thus the
parameter space may be identified (up to dividing by automorphisms)
with the total space of the  Poincar\'e bundle $\Pcal \to B\times B$ 
minus the $0$--section $P_0$. It is then natural to compactify this 
by looking at the associated $\PP^1$--bundle and by glueing on it the 
$0$--section $P_0$ with the infinity section $P_\infty$. This in fact 
works and it is explained in \cite{Mum}, \S 7, and in \cite {Nam}. 
We describe next the new objects that arise.

We denote by $L_i$ the line bundle associated to $b_i$, for $i=1,2$.
We consider again the $\PP^1\times \PP^1$--bundle $P$ on
$B$ as in \S \ref {sec:stand} and in the glueing operation 
described in \S \ref {sec:stand}, we let $t=t_{12}$ tend to $0$ 
(or equivalently to $\infty$). 
Letting $t\to \infty$, one contracts $\Pi_{1,0}$ and $\Pi_{2,0}$ 
to the section $A=\Pi_{1,\infty}\cap \Pi_{2,\infty}\cong B$, and 
$\Pi_{1,\infty}$ and $\Pi_{2,\infty}$ to the
section $\Delta=\Pi_{1,0}\cap \Pi_{2,0}\cong B$. 
In order to properly describe the glueing process, we have first to
blow--up the two sections $A$ and $\Delta$ in $P$. Let us do that. 
Let $w:\tilde P\to P$ be the blow--up, on which we have the following 
divisors:
$\alpha$ is the exceptional divisor over $A$ and
$\delta$ is the exceptional divisor over $\Delta$;
$\beta, \gamma, \epsilon, \zeta$ are the proper transforms on 
$\tilde{P}$ of $\Pi_{1,\infty}, \Pi_{2,0}, \Pi_{1,0}, \Pi_{2,\infty }$,
respectively.

We will abuse notation and denote by the same letters the restrictions of
these divisors on the general fibre $\Phi$ of $\tilde P$ over $B$, which is
a $\PP^1\times \PP^1$ blown--up at two points, hence a $\PP^2$ blown--up
at three points.
Note that $\alpha, \beta, \gamma, \delta, \epsilon, \zeta$ are 
$\PP^1$--bundles over $B$ and one has
\begin{equation}\label {eq:tab}
\begin{matrix}
\alpha\cong \PP(L_1^{\vee}\oplus L_2^{\vee}), & \gamma\cong 
\PP(L_2\oplus \Ocal_B),& \epsilon\cong \PP(L_1\oplus
\Ocal_B), \\
\delta\cong \PP(L_1\oplus L_2), & \zeta\cong \PP((L_1\otimes L_2)\oplus
L_2),& \beta\cong \PP((L_1\otimes L_2)\oplus L_1).\\ 
\end{matrix}
\end{equation}
At this point one could be tempted to suitably glue $\alpha$ with $\gamma$
and $\epsilon$ and $\delta$ with $\beta$ and $\zeta$. This however, as
\eqref{eq:tab} shows,  does not work. The right construction is instead the
following.

One considers two $\PP^2$--bundles $\phi_i:P^{\sharp}_i\to
B$, $i=1,2$, associated to the vector bundles $L_1\oplus L_{2}$
and $L^{\vee}_1\oplus L^{\vee}_{2}$ on $B$, i.e.,
$$
P^{\sharp}_1=\PP(  L^{\vee}_1\oplus L^{\vee}_{2}\oplus\Ocal_B),\quad
P^{\sharp}_2=\PP(L_1\oplus L_{2}\oplus\Ocal_B).
$$
There are three
relevant $\PP^1$--subbundles of the bundles $P_i^{\sharp}$, $i=1,2$, 
namely
\begin{equation}\label {eq:tab2}
\begin{matrix}
\bar \alpha= \PP(L_1^{\vee}\oplus L_2^{\vee}), & \bar \gamma =
\PP(L_2^{\vee}\oplus\Ocal_B), & \bar \epsilon= \PP(L_1^{\vee}\oplus\Ocal_B)
 & \text{in} \quad  P_1^{\sharp}\\
\bar \delta= \PP(L_1 \oplus L_2), & \bar \zeta =\PP(L_1
\oplus\Ocal_B), & \bar \beta= \PP(L_2\oplus\Ocal_B) & 
\text{in}\quad P_2^{\sharp}.\\
\end{matrix}
\end{equation}
As \eqref {eq:tab} and \eqref {eq:tab2} show, we can glue $P$ with $P_1^
\sharp$ and $P_2^ \sharp$ in such a way that
$\alpha$ and $\delta$ are respectively glued to $\bar \alpha$ and
$\bar \delta$; $\epsilon$ is glued to $\bar \epsilon$ and 
$\beta$ to $\bar \beta$  with a shift by $-b_1$, and
$\zeta$ to $\bar \zeta$ and 
$\gamma$  to $\bar \gamma$ with a shift by $-b_2$.
The resulting variety is denoted by $\overline X=\overline X_{B,b}$.
As usual, we will denote by $D$ its singular locus.
On $\tilde P$ we have the line bundle
$$
\tilde
M=w^*\Ocal_P(\Pi_{1,\infty}+\Pi_{2,\infty})\otimes\Ocal_{\tilde
P}(\Xi-\alpha),
$$
where we write $L$ instead of $w^*(\pi^*(L))$ for a line bundle, or
divisor, $L$ on $B$. With similar notation, one has
\begin{equation}\label {eq:relat}
\begin{aligned}
\tilde M\cong \Ocal_{\tilde P}(\alpha+\beta+\zeta+\Xi)
&\cong \Ocal_{\tilde P}(\delta+\epsilon+\zeta+\Xi_{b_1})\cr
&\cong \Ocal_{\tilde P}(\delta+\gamma+\beta+\Xi_{b_2}).\cr
\end{aligned}
\end{equation}
One has $\Ocal_{P_1^ \sharp}(1)=\Ocal_{P_1^ \sharp}(\bar\alpha)$ and the
following linear equivalences
\begin{equation}\label{eq:rel2}
\bar \alpha-\bar\gamma\equiv L_2^{\vee},\quad
\bar\alpha-\bar\epsilon\equiv L_1^{\vee}
\end{equation}
where again we write $L$ instead of $\phi_i^*(L)$, $i=1,2$, for a line
bundle, or divisor, $L$ on $B$ (see again \cite  {H}, Proposition 2.6,
p. 371). From \eqref {eq:relat} and
\eqref {eq:rel2} one deduces that $\tilde M$ glues with the line bundle
$M_1^{\sharp}=\Ocal_{P_1^{\sharp}}(\bar\alpha+\Xi)$ 
and with the line bundle
$M_2^{\sharp}=\Ocal_{P_2^{\sharp}}(\Xi)$, to give a line bundle $\overline M$
on  $\overline X$.
In the evident coordinates
$(u,z)=((u_1,u_2),(z_1,\ldots,z_{g-2}))$, which can be considered
as coordinates on $\tilde P-(\alpha\cup \cdots \cup \zeta)$,
the sections of $\tilde M$ can be expressed as
$$
\begin{aligned}
a\, \xi(\tau,z)+a_1u_1\, \xi(\tau,z-\omega_1)+a_2u_2\, \xi(\tau,z-\omega_2)
\end{aligned}
$$
with $a, a_1, a_2$ complex numbers. By taking into
account the glueing conditions, one sees that only a one--dimensional
subspace $V$ of $H^ 0(\tilde P,\tilde M)$ gives rise to a space of
sections of $H^ 0(\overline X,\overline M)$; $V$ is generated by
\begin{equation}\label{eq:tetaeq}
\xi(\tau,z)+u_1\xi(\tau,z-\omega_1)+u_2\xi(\tau,z-\omega_2).
\end{equation}
in the $(u,z)$-coordinates. Note that \eqref {eq:tetaeq} is just gotten 
from the generalized theta function in \ref{sec:stand} by 
letting $t=t_{12}$ tend to $0$.

In conclusion, one has $h^ 0(\overline X, \overline M)=1$, hence there
is a  unique effective divisor ${\overline \Xi}={\overline
\Xi}$ which is the zero locus of a non--zero section of $H^
0(\overline X,\overline M)$.  

As in the standard case, we parametrize an open subset of ${\rm
Pic}^ 0(\overline \Xi)$ with points in the union of
$P-\bigcup_{i=1,2;h=1\infty} \Pi_{i,h}$ with $P_1^
\sharp-(\bar\alpha\cup\bar\gamma\cup \bar\epsilon)$ and
$P_1^\sharp-(\bar\delta\cup\bar\zeta\cup \bar\beta)$.
We can define the {\it vertical singularities} of the divisor
${\overline \Xi}$, whose equations, in the $(u,z)$ coordinates, take the
form
\begin{equation}\label{eq:vertrank2}
\begin{aligned}
\xi(\tau,z)=0,\quad
u_1\xi(\tau,z-\omega_1)=0,\quad
u_2\xi(\tau,z-\omega_2)=&0,\\
\partial_i\xi(\tau,z)+
u_1\partial_i\xi(\t,z-\omega_1)+u_2
\partial_i\xi(\tau,z-\omega_2)=&0,\quad i=1,\ldots, g-2. \\
\end{aligned}
\end{equation}
Again the vertical singular points of $\overline \Xi$ are
essentially of three types:

(i) Consider a point $x\in D$ which is the image via
$\varphi:X\to \overline X$ of a point in $\Pi_{1,0}\cap\Pi_{2,0}\cong B$,
i.e. in the $(u,z)$ coordinates one has $u_1=u_2=0$. Then this is a vertical
singular point of $\overline \Xi$ if and only if it corresponds
to a singular point of $\Xi$ on $\Pi_{1,0}\cap\Pi_{2,0}\cong B$ and to a
singular point of $\Xi_{b_1+b_2}$ on $\Pi_{1,\infty}\cap\Pi_{2,\infty}\cong
B$.

(ii) Consider a point $x\in D$ which is the image via $\varphi$ of
a point in $\Pi_{i,0}$ but not of a point in $\Pi_{3-i,0}$,
i.e. in the $(u,z)$ coordinates one has $u_i=0,u_{3-i}\neq 0$, for
an $i=1,2$. If $i=1$, this is a vertical singular point of
$\overline \Xi$ if and only if
$$
\begin{aligned}
\xi(\t,z)=0, \quad
\xi(\t,z-\omega_2)=&0,\\
\partial_i\xi(\tau,z)+
u_2 \, \partial_i\xi(\t,z-\omega_2)=&0,\quad i=1,\ldots, g-2. \\
\end{aligned}
$$
i.e., if and only if $z$ and $z-b_2$ are in $\Xi$ and
$\gamma_{\Xi}(z)=\gamma_{\Xi}(z-b_2)$. Thus points of this type correspond
to points in $N_{0}(B,\Xi)$.

(iii) Consider a point $x\notin D$, i.e., in the $(u,z)$
coordinates one has $u_1\neq 0\neq u_2$. Then equations \eqref
{eq:vertrank2} mean that $z$ corresponds to a point in
$\Xi\cap\Xi_{b_1}\cap\Xi_{b_2}$ where $\Xi,\Xi_{b_1},\Xi_{b_2}$ are
tangentially degenerate. In other words points of this type correspond to
points in $N_{0,2}(B,\Xi )$.  

\begin{remark} \label {rem:coni3} A point $x$ in ${\rm Sing}^
v({\overline \Xi})$ determines a quadric $Q_x$ in $\PP^{g-1}$.  Remark \ref
{rem:coni2} is still valid here. \end{remark}

A variant of this second type of rank--$2$ degeneration is
obtained as follows. Given a $\GG_m^2$-extension $X=X_b$ of $B$
determined by a point $b=(b_1,b_2)\in \hat{B}^2$ we consider
two $\PP^2$-bundles $\PP$ and $\PP^{\prime}$ over $B$:
$$
\PP=\PP({\Ocal}_B\oplus L_{1}\oplus L_{2})
\qquad {\rm and}\qquad
\PP^{\prime}=\PP(L_{2}\oplus L_{1}\oplus (L_1\otimes L_2)),
$$
where we write as before $L_i$ for $L_{b_i}$.
We can glue these over the common $\PP^1$-subbundle
$\PP(L_{1}\oplus L_{2})$. Then we glue the $\PP^1$-subbundle
$\PP({\Ocal}_B\oplus L_{1})$ of $\PP$ with the $\PP^1$-subbundle
$\PP(L_{2}\oplus (L_1\otimes L_2))$ of $\PP^{\prime}$ via a shift over
$b_2$. Similarly, we glue the $\PP^1$-subbundle
$\PP({\Ocal}_B\oplus L_{2})$ of $\PP$ with the $\PP^1$-subbundle
$\PP(L_{1}\oplus (L_1\otimes L_2))$ of $\PP^{\prime}$ via a shift over
$b_1$. In this way we obtain a non-normal variety over $B$.
Both $\PP$ and $\PP^{\prime}$ come with a relatively ample bundle
${\Ocal}_{\PP}(1)$ and ${\Ocal}_{\PP^{\prime}}(1)$.
On $\PP$ we have the linear equivalences
$$
\Pi_1+\pi^*(\Xi_{b_1}) \equiv \Pi_2+\pi^*(\Xi_{b_2}) \equiv 
\Pi_3+\pi^*(\Xi_{b_1+b_2}).
$$
with $\Pi_i=\PP({\Ocal}_B\oplus L_{i})$ for $i=1,2$ and
$\Pi_3=\PP(L_{1}\oplus L_{2})$. We let $M$ be the line bundle
${\Ocal}(\Pi_3+\pi^*(\Xi_{b_1+b_2}))$ on $\PP$ and $M'$ the line bundle
${\Ocal}(\Pi_3^{\prime}+\pi^*(\Xi_{b_1+b_2}))$ on $\PP^{\prime}$,
where $\Pi_3^{\prime}$ is the bundle $\PP(L_{1}\oplus L_{2}$).
This descends to a line bundle $\overline{M}$
on $\overline{X}$. This line bundle has a $1$-dimensional space of sections
generated by 
$$
\theta(\tau,z)=\xi(\tau,z)+u_1\, \xi(\tau,z_1-\omega_1)+
u_2\, \xi(\tau,z-\omega_2)
$$
in suitable affine coordinates $(u_1,u_2)$ on $\PP^2$.
Again the vertical singular points of $\overline \Xi$ are
essentially of three types:

(i) A point $x\in D$ which is the image via $\varphi: X \to \overline X$ 
of a point in $\Pi_{1}\cap\Pi_{2}=\PP({\Ocal}_B) \cong B$ is a vertical
singularity if it corresponds to a singularity on $\Xi$, to a singularity
on $\Xi_{b_1}$ on $\Pi_1\cap \Pi_3$ and a singularity on $\Xi_{b_2}$
on $\Pi_2\cap \Pi_3$;

(ii) A point $x\in D$ which is the image via $\varphi$ of
a point on one $\Pi_{3}$ but not of a point in $\Pi_1$ or $\Pi_2$
is a vertical singular point of $\overline \Xi$ if and only if
$x \in \Xi_{b_1}\cap \Xi_{b_2}$ and $\gamma_{\Xi_{b_1}}(x)=
\gamma_{\Xi_{b_2}}(x)$. Thus points of this type correspond
to points in $N_{0}(B,\Xi)$. Something similar happens for the
points on exactly one of $\Pi_1$ or $\Pi_2$.

(iii) A point $x\notin D$ is a vertical singularity if 
$x\in \Xi\cap\Xi_{b_1}\cap\Xi_{b_2}$ and $\Xi,\Xi_{b_1},\Xi_{b_2}$ are
tangentially degenerate at $x$. 
In other words points of this type correspond to points in $N_{0,2}(B,\Xi)$.  

We see that the compactification depends on a choice, but in both cases
we can deal explicitly with the singularities of the theta divisors.

\begin{remark}\label{variant}
 The variant just given corresponds to a tesselation
of $\RR^2$ given by the lines $x=a$, $y=b$ and $x+y=c$ with $a,b,c \in \ZZ$.
To a triangle with integral vertices $(n,m), (n+1,m)$ and $(n,m+1)$ 
(resp.\ $(n-1,m)$, $(n,m-1)$ and $(n,m)$) we associate
the $\PP^2$-bundle $\PP(L_{nb_1+mb_2}\oplus L_{(n+1)b_1+mb_2}\oplus 
L_{nb_1+(m+1)b_2})$ (resp.\ $\PP(L_{(n-1)b_1+mb_2}\oplus L_{nb_1+(m-1)b_2}\oplus
L_{nb_1+mb_2}$) and we glue the bundles belonging to adjacent triangles
over the $\PP^1$-bundle defined by the common edge. Then the generator
$\tau_1$ (resp.\ $\tau_2$) of  $\ZZ^2$ acts
by glueing the $\PP^1$-bundle associated to an edge to the $\PP^1$-bundle
associated to the translate by $x \mapsto x+1$ (resp.\ $y \mapsto y+1$) 
of this edge using a translation over $b_2$ (resp.\ $b_1$).
The quotient under $\ZZ^2$ is the non-normal variety we just constructed.
Also the earlier compctifications thus correspond to tesselations (see \cite 
{Mum}, \S 7).
\end{remark}
\end{section}

\begin{section}{Compactification of $\Ag$}\label{sec:comp}
In order to study the Andreotti-Mayer loci we need to compactify $\Ag$. 
The moduli space $\Ag$ admits a `minimal' compactification,
the {\sl Satake compactification} constructed by Satake and Baily-Borel in
characteristic $0$ and by Faltings-Chai over the integers 
(see \cite {Sat}, \cite{FC}). This compactification
${\mathcal A}_g^*$ is an orbifold or stack which admits a stratification
$$
{\mathcal A}_g^*=\sqcup_{i=0}^g {\mathcal A}_{i}
$$
and the closure of ${\mathcal A}_{m}$ in ${\mathcal A}_g^*$ is 
${\mathcal A}_{m}^*=\sqcup_{i=0}^{m} {\mathcal A}_{i}$. This compactification
is highly singular for $g\geq 2$. Smooth compactifications can be constructed
by the theory developed by Mumford and his co-workers in characteristic $0$ 
and by Faltings-Chai in general. These compactifications depend on 
combinatorial data. We shall use the {\sl Voronoi compactification}
$\tilde{\mathcal A}_g= \tilde{\mathcal A}_g^{\rm Vor}$ 
as described by Namikawa, Nakamura and Alexeev (see \cite{Nam, Nakam, Br}). 
This compactification is a smooth orbifold with a 
natural map $q: \tilde{\mathcal A}_g \to {\mathcal A}_g^*$. 
It has the stratification induced by that of ${\mathcal A}_g^*$:
the stratum
$$
{\mathcal A}_g^{(r)}= q^{-1}({\mathcal A}_{g-r})
$$
is called the stratum of {\sl quasi-abelian varieties of torus rank} $r$.
The general fibre of 
$q: {\mathcal A}_g^{(r)} \to {\mathcal A}_{g-r}$ 
has a dense open subset consisting of standard compactifications 
as seen in \S \ref {sec:stand}. Hence the general fibre of 
$q: {\mathcal A}_g^{(r)} \to {\mathcal A}_{g-r}$
has dimension $gr-r(r+1)/2$, thus 
$\dim({\mathcal A}_g^{(r))})=g(g+1)/2-r$.

The closure of ${\mathcal A}_g^{(r)}$ is 
${\mathcal A}_g^{(\geq r)}=q^{-1}({\mathcal A}_{g-r}^*)$.   
We set $\dd\tilde \Ag:={\mathcal A}_g^{(\geq 1)}$, which is the
{\it boundary} of $\tilde \Ag$.  Moreover we set 
${\mathcal A}_g^{(\leq r)}:=\tilde{\mathcal A}_g-{\mathcal A}_g^{(\geq r+1)}$.

The fibres of the map 
$q: {\mathcal A}_g^{(\leq r)}\to {\mathcal A}_g^*$
are well behaved if $r\leq 4$. Indeed, the points of 
$\Ag':= {\mathcal A}_g^{(\leq 4)}$
correspond to so-called {\it stable} quasi--abelian varieties
which are compactifications of semi-abelian varieties, which can
be explicitely described (see \cite{Nam, Nakam, Br} and 
\S\S  \ref {sec:rank1}, \ref {sec:rank2}
for torus ranks 1 and 2). Thus one can define the vertical singular locus
and the Andreotti-Mayer loci on
the partial compactification $\Ag'$ (seeRemark \ref {rem:vsing} below).
For higher torus rank the situation is more complicated. For instance
the fibres of $q: {\mathcal A}_g^{(\geq r)} \to {\mathcal A}_{g-r}^*$
might be non-reduced. However we will not need $r>4$ here.

An alternative approach might be to use the idea of Alexeev and Nakamura 
(cf.\ \cite{AN}, \cite {Nakam}, \cite {Br}) who
have constructed canonical limits for one-dimensional families of abelian
varieties equipped with principal theta divisors.

The stable quasi-abelian varieties that occur in $\Ag'$ for torus rank $1$
and $2$ are exactly those described in Section \ref{sec:rank1}
and \ref{sec:rank2}. 
For torus rank $3$ these are described by the tesselations of $\RR^3$
occuring on p.\ 188 of \cite {Nakam}, cf.\ also Remark \ref{variant}. 
The open stratum of
${\Acal}_g^{(3)}$ over ${\Acal}_{g-3}$
corresponds to the standard compactifications 
(see  \S \ref {sec:stand}), obtained by
glueing six $\PP^3$-bundles
over a $(g-3)$--dimensional abelian variety $B$ generalizing the construction
for torus rank $2$ where two $\PP^2$--bundles were glued. These closed strata
correspond to degenerations of the matrix $T$ of the glueing data
on which the standard compactifications depend (see  \S \ref {sec:stand}).
For instance the codimension $3$
stratum corresponds to the fact that in $T$ two of the three elements above the
main diagonal tend to zero (or to $\infty$).

\begin{remark}\label{rem:vsing} 
As we remarked before, for stable quasi-abelian varieties corresponding
to points $(\overline X, \overline \Xi)$ of  ${\mathcal A}_g^{\prime}$
one can define the vertical singularities
${\rm Sing}_{\rm vert}(\overline \Xi)$ using $\Omega^1(\log D)$ as in the previous sections.
One checks that for these compactifications the analogue of
Lemma \ref{stand} still holds.
\end{remark}
We
will need the following result from \cite {vdG1}.

\begin{theorem}\label {thm:ampcond} Let $Z$ be an irreducible, closed
subvariety of $\tilde \Ag$. Then $Z\cap \dd\tilde \Ag$ is not empty as
soon as ${\rm codim}_{\tilde \Ag}(Z)< g$.
\end{theorem}

\end{section}
\begin{section}{The Andreotti--Mayer loci and the boundary}\label{sec:AMbound}

As we are working with a fixed compactification 
$\tilde{\mathcal A}_g={\mathcal A}_g^{\rm Vor}$ and, as indicated in \S 
\ref {sec:comp} above,  we
can define the Andreotti--Mayer loci over the part 
${\mathcal A}_g^{\prime}={\mathcal A}_g^{(\leq 4)}$ of 
$\tilde{\mathcal A}_g$.
We have $N_{g,k}$ as a subscheme of ${\mathcal A}_g$ and
we define $\tilde N_{g,k}$  as the schematic closure.
The support of $\tilde N_{g,k}$ contains the
set of points corresponding to pairs $(\overline X,\overline \Xi)$ such
that ${\rm Sing}_{\rm vert}(\overline X)$ has a component of dimension
at least $k$ (see \cite {Mumford}).

It is interesting to look at the case
$k=0$, which has been worked out by Mumford \cite{Mumford} and
Debarre \cite{D1}. In this case $\tilde N_{g,0}$ is a divisor and by
Theorem \ref  {thm:ampcond}, every irreducible component $N$ of this
divisor intersects $\dd\tilde \Ag$. Let $M$ be an irreducible component of
$N\cap \dd\tilde{\Ag}$, which has dimension $\binom{g+1}{2} -2$. 

First of all, notice that $M$ cannot be equal to  
${\mathcal A}_g^{(\geq 2)}$. This
follows by the results in \S \ref{sec:rank2}  and by Propositions 
\ref{prop:partial} and \ref{prop:singolare}. More generally, in the
same way, one proves that $M$ cannot contain $\Acal_{g}^{(r)}$ for any
$r=2,3$ and $4$.

Hence  $M$ intersects $\Acal_{g}^{(1)}$ in a non-empty 
open set of $M$, i.e., the intersection with the boundary has points
corresponding semi-abelian varieties of torus rank~$1$. 
If $M$ does not dominate $\Acal_{g-1}$ via the map 
$q$, then each fibre must have dimension $g-1$. 
By Proposition \ref{prop:partial} this
implies that $M$ dominates $N_{0,g-1}$.  

If $M$ dominates $\Acal_{g-1}$ via $q$, the fibre of
$q_{\vert M}$ over a general point $(B,\Xi)\in \Acal_{g-1}$ is
$N_0(B,\Xi)$. 

Recall now that Debarre proves in \cite{D1} that $N_{0,g}$ has two
irreducible components, one of which is the so--called {\it theta--null}
component $\theta_{0,g}$:  the general abelian variety $(X,\Theta_X)$ in 
$\theta_{0,g}$,
with $\Theta_X$ symmetric, is such that $\Theta_X$ has a unique double point
which is a $2$--torsion point of $X$ lying on $\Theta_X$.

Let $M_{0,g}$ be the other component.  The general abelian variety 
$(X,\Theta_X)$ in $M_{0,g}$, with $\Theta_X$ symmetric, is such that $\Theta_X$ 
has exactly two double points $x$ and $-x$.  It is useful to
remind that, by Corollary \ref {cor:sgsmooth}, at a general point of either
one of these component of $N_{g,0}$, 
the tangent cone to the theta
divisor at the singular points is a smooth quadric.

The component $\theta_{0,g}$ intersects the boundary in two components,
$\theta'_{0,g}, \theta''_{0,g}$, one dominating $\theta_{0,g-1}$, the other
dominating $\Acal_{g-1}$ with fibre over the general point $(B,\Xi)\in
\Acal_{g-1}$ given by the component $2\Xi$ of $N_0(B,\Xi)$ (see \S \ref
{sec:gauss}). Also $M_{0,g}$ intersects the boundary in two irreducible
divisors $M^{\prime}_{0,g}, M^{\prime \prime}_{0,g}$, 
the former is irreducible and dominates
$M_{0,g-1}$, the latter  dominates $\Acal_{g-1}$ with fibre over the
general point $(B,\Xi)\in \Acal_{g-1}$ given by the components of
$N_0(B,\Xi)$ different from $2\Xi$.

The main ingredient for Debarre's proof of the irreducibility of 
$M_{0,g}$ is a monodromy argument which implies that, if $(B,\Xi)$ is a
general principally polarized abelian variety of dimension $g$, then
$N_0(B,\Xi)$ consists of only two irreducible components.

\begin{remark}\label{rem:exttheta0}{\rm Let $(B,\Theta_B)$ be a general element
in $\theta_{0,g}$ and let $(\overline X,\overline \Xi)$ be a 
semi--abelian variety of torus rank one with abelian part $B$. Then there are
no points in ${\rm Sing}_{\rm vert}(\overline \Xi)$ with multiplicity larger than $2$.
This follows from the fact that 
$\Theta_B$ has a unique singular point 
and by Remark \ref {rem:triplepts}. }\end{remark}

We finish this section with the following
result which will be useful later on. It uses the notion of
asymptotic cone given in  \S \ref{higherorder}.

\begin{proposition}\label{prop:null} One has:

\begin{itemize}

\item [(i)] let $g\geq 3$, let $(B,\Xi)$ be
a general point of $\theta_{0,g}$ and let $x$ be the singular point of
$\Xi$. Then the asymptotic cone $TC^{(4)}_\xi$ is strictly contained in
the quadric tangent cone $Q_x$;

\item [(ii)] let $g\geq 4$, let $(B,\Xi)$ be
a general point of $M_{0,g}$ and let $x, -x$ be the singular points of
$\Xi$. Then the asymptotic cone $TC^{(3)}_\xi$ is strictly contained in
the quadric tangent cone $Q_x=Q_{-x}$.
\end{itemize}
\end{proposition}

\begin{proof} Degenerate to the jacobian locus and apply the results from
\cite{Kempf,KS}.
\end{proof}

\end{section}
\begin{section}{The Conjecture for $N_{1,g}$}\label{sec:mainth}

In this section we prove Conjecture \ref{conj:mainconj} for $k=1$.
We consider an irreducible component $N$ of $\tilde{N}_{g,1}$ which is
of codimension $3$. The first observation is that the assumption
about the codimension of $N$ implies that the generic principally
polarized abelian variety is simple since by Remark \ref{rem:fac} 
every irreducible component of ${\mathcal A}_g^{\rm NS}$
has codimension $\geq g-1>3$ if we assume $g\geq 5$.

\begin{proposition}\label{rank1exists}
Let $g \geq 6$ and let $N$ be an irreducible component of $\tilde{N}_{g,1}$
which is of codimension $3$ in $\tilde{\mathcal A}_g$. Then 
$N$ intersects the stratum ${\mathcal A}_g^{(1)}$.
\end{proposition}
\begin{proof}
We begin by remarking that $N$ cannot be complete in ${\Ag}$
in view of Theorem \ref{thm:ampcond}.
Therefore $N$ intersects $\partial \tilde{\Ag}$.
Since $\partial \tilde{\Ag}$ is a divisor in $\tilde{\Ag}$ it intersects
$N$ in codimension one. Let $M$ be an irreducible component
of $N \cap \partial \tilde{\Ag}$. It is our intention to prove that
$M$ has a non-empty intersection with ${\mathcal A}_g^{(1)}$.

Suppose that $M\subseteq {\mathcal A}_g^{(\geq 4)}$. 
For dimension reasons
we have $M={\mathcal A}_g^{(\geq 4)}$. 
Since we are using a compactification $\tilde{\mathcal A}_g$ 
such that the general point of ${\mathcal A}_g^{(4)}$
corresponds to a standard compactification $(\overline{X},
\overline{\Xi})$ of
torus rank $4$ with abelian part $(B,\Xi) \in {\Acal}_{g-4}$
we deduce from Lemma \ref{stand}
and Remark \ref {rem:vsing} that if
$\dim({\rm Sing}_{\rm vert}(\overline{\Xi}))\geq 1$ then 
$(B,\Xi) \in N_{g-4,0}$.
But for $g\geq 5$ the locus $N_{g-4,0}$ is a divisor in 
${\mathcal A}_{g-4}$ and we obtain the inequality
${\rm codim}_{\partial \tilde{\Ag}}(M)\geq 
{\rm codim}_{\partial \tilde{\Ag}}({\mathcal A}_g^{(4)})+1\geq 5$, 
a contradiction. 
Therefore we can assume that 
$M \cap {\mathcal A}_{g}^{(\leq 3)}\neq \emptyset$. 

Suppose that
$M\cap {\mathcal A}_{g}^{(3)}$ has codimension $1$ in 
${\mathcal A}_{g}^{(3)}$. Then either $M$ maps dominantly to
$\tilde{\mathcal A}_{g-3}$ via the map 
$q: {\Acal}_{g}^{(3)} \map \tilde{\Acal}_{g-3}$
and $M$ intersects the general fibre $F$ of 
$q$ in a divisor, 
or $M$ maps to a divisor in ${\mathcal A}_{g-3}$ under $q$
with full fibres $F$ contained in $M$. 

The former case is impossible by Proposition \ref {prop:singolare}.
In the latter case for a general $(B,\Xi)$ in $q(M)$
all the quasi--abelian varieties $(\overline{X},\overline \Xi)$ 
in the fibre $F$ over $(B,\Xi)$ must have a $1$-dimensional vertical 
singular locus of $\overline \Xi$.
Note that $(\overline{X},\overline \Xi)$  corresponds to a 
standard compactification as considered in \S \ref{sec:stand}. 
By the discussion given in  \S \ref {sec:stand} and
by Proposition \ref {prop:partial}, we see that $q(M)$ has 
to be contained in $N_{g-3,1}$, against Theorem \ref {th:amlower}.  
We thus conclude that
$M \cap {\mathcal A}_g^{(\leq 2)}\neq \emptyset$.

Suppose that $M \cap {\mathcal A}_g^{(\geq 2)} \neq \emptyset$
which has codimension $2$ in ${\mathcal A}_g^{(2)}$. 
As above we have that  $q(M)$ is contained in
$N_{g-2,0}$.

Suppose first $q(M)$ is dense in a component of $N_{g-2,0}$. 
If $(B,\Xi)$  is a general element of $q(M)$, then $M$ intersects 
the fibre of $q: {\mathcal A}_g^{(2)}\map \tilde{\mathcal A}_{g-2}$ over 
$(B,\Xi)$ in codimension one.  This gives a contradiction by the 
analysis \S \ref {sec:rank2} and Proposition \ref {prop:partial}.

Suppose that $q(M)$ is not dense  in a component of $N_{g-2,0}$. 
If $(B,\Xi)$  is a general element of $q(M)$, then $M$ contains 
the full fibre of $q: {\mathcal A}_g^{(2)}\to {\mathcal A}_{g-2}$ 
over $(B,\Xi)$.  By taking into account the analysis \S \ref {sec:rank2} and
Proposition \ref {prop:partial}, this implies $q(M)$ contained in 
$N_{g-2,1}$, giving again a contradiction.

This proves  that $M \cap {\mathcal A}_g^{(\leq 1)}\neq \emptyset$.
\end{proof}

Let $g\geq 4$ and let $N$ be an irreducible component of 
$N_{g,1}$ of codimension $3$ in $\Ag$. As in the
proof above we denote by $M$ an irreducible component of the intersection
of the closure of $N$ in $\tilde{\Ag}$ with the boundary
$\partial\tilde{\mathcal A}_g$. According to Lemma
\ref{rank1exists} the morphism $q: \tilde{\Ag} \to {\mathcal A}_g^*$ induces
a rational map $\alpha: M\map \tilde\Acal_{g-1}$, whose image is not
contained in $\partial \tilde {\mathcal A}_{g-1}$.

\begin{lemma}\label{cor:struc} In the above setting, the Zariski
closure of $q(M)$ in $\tilde \Acal_{g-1}$ is:
\begin{itemize}
\item [(i)] either an irreducible component $N_1$ of $\tilde{N}_{g-1,1}$
of codimension 3 in $\tilde \Acal_{g-1}$;
\item [(ii)] or an irreducible
component $N_0$ of $\tilde N_{g-1,0}$ and in this case:
\begin{itemize}
\item [(a)] if $\eta=(B,\Xi)\in N_0$ is a general point, then the
closure of $q^{-1}(\eta)$ in $B$ is an irreducible
component of $N_1(B,\Xi)$ of codimension~$2$ in $B$;
\item [(b)] if $\xi=(\overline X,\overline \Xi)\in M$ is a general
point, then ${\rm Sing}_{\rm vert}(\overline \Xi)$ meets the singular 
locus $D$ of $\overline X$ in one or two points, 
whose associated quadric has corank~$1$.
\end{itemize}
\end{itemize}
\end{lemma}
\begin{proof} 
If $q(M)\subseteq N_{g-1,1}$, then  Theorem \ref{th:amlower} implies
$
3 \leq {\rm codim}_{{\mathcal A}_{g-1}}(q(M)) \leq 
{\rm codim}_{\partial \tilde{\mathcal A}_{g}}(M) =3
$
and the closure of $q(M)$ must an irreducible component of 
$\tilde{N}_{g-1,1}$.
If $q(M) \not\subseteq N_{g-1,1}$ then by Proposition\ \ref{prop:singolare}
we have that $q(M)\subseteq N_{g-1,0}$ and the fibre $q^{-1}(B,\Xi)
\subseteq N_1(B,\Xi)$. By Proposition  \ref{prop:partial} we 
have ${\rm codim}_B(N_1(B,\Xi))
\geq 2$ and since $N_{g-1,0}$ is a divisor in ${\mathcal A}_{g-1}$ we see that
(iia) follows. The last statement (iib) follows from Remark \ref{rem:piu},
the analysis in Section \ref{sec:rank2}, the description of 
$N_{g,0}$ by Mumford and Debarre (see \cite{Mumford}, \cite{D1}, and 
\S \ref {sec:AMbound}) and  Corollary \ref {cor:sgsmooth}. \end{proof}

We are now ready for the proof of the conjecture for $N_{g,1}$.

\begin{theorem}\label{th:con1}
Let $g\geq 4$. Then the codimension of an irreducible component $N$ of
$N_{g,1}$ in ${\mathcal A}_g$ is at least $3$ with equality if and only if:
\begin{itemize}
\item [(i)]  $g=4$ and either $N={\Hcal}_4$ is the hyperelliptic locus or 
$N={\mathcal A}_{4,(1,3)}$;
\item [(ii)] $g=5$ and  $N={\Jcal}_5$ is the jacobian locus.
\end{itemize}
\end{theorem}
\begin{proof} By Theorem \ref {th:amlower}, the codimension of $N$ is at
least $3$. Suppose that $N$ has codimension $3$. It is well known
that the assertion holds true for $g=4$ and $5$ (see \cite {B2}, \cite{D0},
\cite{C-vdG}).  We may thus assume $g\geq 6$ and proceed by induction.

Let $\zeta=(X,\Theta_X)$ be a general point of $N$, so that $X$ is
simple  (see Remark \ref {rem:fac}).  
Let $S$ be a 1--dimensional component of ${\rm Sing}(\Theta_X)$. 
We can assume that the class of  $S$ in $X$ is a 
multiple $m\gamma_X$ of the {\it minimal class} $\gamma_X=
\Theta_X^{g-1}/(g-1)! \in H^2(X,\ZZ)$. If not so, then
${\rm End}(X)\neq \ZZ$ and this implies that ${\rm codim}_{\Ag}(N) \geq g-1$
(see Remark \ref{rem:fac}).

By Theorem \ref {th:amlower}, the general point in $S$ is a
double point for $\Theta_X$.  We let $R$ be the curve in $S_g$
whose general point is $\xi=(X,\Theta_X,x)$, with $x\in S$ a general point.
Note that $R$ is birationally equivalent to $S$.
Let $\Qcal$ be the linear system of all quadrics in $\PP(T_{X,0})$. One has
the map
$$
\phi: \xi\in R\map Q_\xi\in \Qcal.
$$
As in the proof of Theorem \ref{th:amlower}, the map $\phi$ is not
constant. Let $\Qcal_R$ be the span of the image of $\phi$.
As in the proof of Theorem \ref{th:amlower}, one has
$\dim(\Qcal_R)\geq 2$. By Proposition \ref{prop:tanam}, $\Qcal_R$
is contained in the linear system $\Ncal_{g,1}(\zeta)$ (see 
\S \ref{sec:amloci} for the definition), which has dimension at most $2$ 
since ${\rm codim}_{\Ag}(N)=3$. Thus $\Qcal_R= \Ncal_{g,1}(\zeta)$ has 
dimension~$2$.

By Lemma \ref{rank1exists}, the closure of $N$ in 
${\mathcal A}_g^{(\leq 1)}$ has non empty intersection with the boundary.
As in the proof of Lemma \ref{rank1exists}, we let $M$ be an irreducible
component of the intersection of the closure of $N$ in $\Acal_g^{(\leq 1)}$ 
with the boundary. Consider the rational map
$\alpha: M\map \Acal_{g-1}$ and the closure of the image $\alpha(M)$, for
which we have the possibilities described in Lemma \ref {cor:struc}.

\begin{claim}\label{cl:caso1}
Possibility (i) in Lemma \ref{cor:struc} does not occur.
\end{claim}

\begin{proof} [Proof of the claim]
By induction, one reduces to the case $g=6$ and $\alpha(M)={\Jcal}_5$. Let
$(\overline X,\overline \Xi)\in M$ be a general point. Then
$(\overline X,\overline \Xi)$ is a general rank one extension of the
jacobian $(J(C),\Theta_C)$ of a general curve $C$ of genus 5. Note that
if $x\in J(C)$ corresponds to the extension,
then $\Theta_C$ and $x+\Theta_C$ are not tangentially 
degenerate (see \cite{H}, Thm.\ 10.8, p.\ 273). Then the
analysis of \S \ref {sec:rank1} implies that the vertical singular locus
$S_0$ of $\overline \Xi$ sits on the singular locus $D\cong J(C)$ of
$\overline X$ and it is isomorphic to $S_C={\rm Sing}(\Theta_C)$
with cohomology class $\Theta^4/12$ (see \cite{ACGH}). 
Thus $\overline \Xi\cdot
S_0=\Theta_C\cdot S_C=10$. Hence, if $\zeta=(X,\Theta_X)$ is a general point
of $N$, then ${\rm Sing}(\Theta_X)$ is a curve $S$ such that
$\Theta_X\cdot S=10$. On the other hand $S$ is homologous to $m\gamma_X$
and one has $10=m\Theta_X\cdot \gamma_X=6m$, a contradiction. 
\end{proof}

Claim \ref{cl:caso1} shows that only possibility (ii) in Lemma \ref
{cor:struc} can occur. In particular, by (iib) of Lemma \ref
{cor:struc}, for $\xi=(X,\Theta_X,x)$ general in $R$, the quadric $Q_\xi$
has corank 1. Let $v_\xi\in \PP^ {g-1}$ be the vertex of $Q_\xi$.
Remember that $R$ is birational to $S$. Hence, by Proposition \ref {kerQ}, the map
$$
\gamma: \xi\in R\map v_\xi\in \PP^ {g-1}
$$
can be regarded as the Gauss map $\gamma_S$ of $S$.

\begin{claim}\label{cl:caso2} 
If the general quadric in the linear system $\Qcal_R$ is singular, 
then  for $\xi=(X,\Theta_X,x)$ general in $R$, the vertex 
$v_\xi$ of $Q_\xi$ is contained in the asymptotic cone $TC_\xi^{(4)}$. 
\end{claim}

\begin{proof} [Proof of the claim] Suppose the general quadric in
$\Qcal_R$ is singular. Then the general quadric in $\Qcal_R$ has corank 1
(see Lemma \ref {cor:struc}, (iib))
and, by Bertini's theorem, its vertex lies in the base locus of $\Qcal_R$.
In particular,  for $\xi=(X,\Theta_X,x)$ general in $R$, the vertex $v_\xi$
of the quadric $Q_\xi$ lies in all the quadrics of $\Qcal_R$.

Choose a local
parametrization $x=x(t)$ of $S$ around a general point of it, with $t$
varying in a disc $\Delta$. Then  $\xi(t)=(X,\Theta_X,x(t))\in R$ and we set
$Q_{\xi(t)}:=Q_t$, its equation being
$$
\sum_{ij}\dd_i\dd_j\theta(x(t))z_iz_j=0,
$$
where we set $\theta(z):=\theta(\t,z)$ for the theta function of $X$. The
main remark is that all the subsequent derivatives of $Q_t$ with respect to
$t$ lie in $\Qcal_R$ and actually $Q_t$ and its first two derivatives $Q'_t$
and $Q''_t$ span $\Qcal_R$, because $\dim(\Qcal_R)=2$. 
Hence Bertini's theorem implies that $x':=x'(s)$
sits on all these quadrics for $t$ and $s$ general in $\Delta$. The
equations of  $Q'_t$ and $Q''_t$ are respectively

$$\sum_{ijh}\dd_i\dd_j\dd_h\theta(x(t))x'_h(t)z_iz_j=0$$

$$\sum_{ijhk}\dd_i\dd_j\dd_h\dd_k \theta(x(t))x'_h(t)x'_k(t)z_iz_j+
\sum_{ijh}\dd_i\dd_j\dd_h\theta(x(t))x''_h(t)z_iz_j=0$$
Thus we have the relations
%%%
\begin{equation}\label
{eq:eqdue}
\begin{aligned}
\sum_{ij}\dd_i\dd_j\theta(x(t))x'_i(s)x'_j(t)=0, \qquad
\sum_{ijh}\dd_i\dd_j\dd_h\theta(x(t))x'_h(t)x'_i(s)x'_j(s)=0\cr
\sum_{ijhk}\dd_i\dd_j\dd_h\dd_k
\theta(x(t))x'_h(t)x'_k(t)x'_i(s)x'_j(s) %%+\cr 
+\sum_{ijh}\dd_i\dd_j\dd_h
\theta(x(t))x''_h(t)x'_i(s)x'_j(s)=0\cr \end{aligned}
\end{equation}
identically in $s,t\in \Delta$. 
The first of these relations says that the tangent hyperplane 
to $Q_t$ at $x'(t)$ contains the vertex $x'(s)$. From the second relation
we have
%%%%
\begin{equation}\label
{eq:equno}
\sum_{ijh}\dd_i\dd_j\dd_h\theta(x(t))x'_h(t)x'_i(t)x'_j(t)=0
\end{equation}
which shows that $v_\xi$ is contained in the asymptotic cone $TC_\xi^ {(3)}$.
By differentiating \eqref {eq:equno}, one finds
$$
\sum_{ijhk}\dd_i\dd_j\dd_h\dd_k
\theta(x(t))x'_h(t)x'_k(t)x'_i(t)x'_j(t)+ 3\sum_{ijh}\dd_i\dd_j\dd_h
\theta(x(t))x''_h(t)x'_i(t)x'_j(t)=0.
$$
By comparing with \eqref {eq:eqdue} for $s=t$, we deduce that
$$
\sum_{ijhk}\dd_i\dd_j\dd_h\dd_k
\theta(x(t))x'_h(t)x'_k(t)x'_i(t)x'_j(t)=0
$$
which proves that $v_\xi\in TC_\xi^ {(4)}$.
\end{proof}

The crucial step in our proof is the following claim.

\begin{claim}\label{cl:caso3} 
The general quadric in the linear system $\Qcal_R$ is non singular. 
\end{claim}
\begin{proof} [Proof of the  claim] 
Suppose this is not the case.
Again we consider  an irreducible component $M$ of
$\bar N\cap\dd\tilde{\Acal}_g'$ and let 
$(\overline X,\overline \Xi)$ be a general
point of $M$. By Claim  \ref {cl:caso1} and Lemma \ref {cor:struc},
$(\overline X,\overline \Xi)$ is a rank 1 extension of $(B,\Xi)$
corresponding to a general point in a component of $N_{g-1,0}$, with
extension datum $b$ varying in a codimension 2 component of $N_1(B,\Xi)$. 
We let $S_0$
be the vertical singular locus of
$\overline \Xi$. By the analysis of \S \ref {sec:rank1}, this corresponds
to a contact curve $C:=C_b$ of $\Xi$ with $\Xi_b$, which contains the
singular points of both $\Xi$ and $\Xi_b$ (see Remark \ref
{rem:piu}). This means that we have a point on $S_0$, corresponding
to a singular point $x$ of $\Xi$, where the tangent cone is the cone over
$Q_x$ with vertex the point of $\PP^ {g-1}$ corresponding to $b$ 
 (see \ S \ref {sec:rank1}).

Now we note that $C$ is smooth at $x$. Indeed locally around $x$, the
divisor $\Xi$ looks like a quadric cone of corank 1  in $\PP^ {g-2}$ and
$\Xi_b$ looks like a hyperplane, which touches it along a curve. This
implies that  $C$ locally at $x$ looks like a line along which a
hyperplane touches a quadric cone of corank 1 (see Remark \ref {rem:casi}).

Hence, the Gauss image $x_b:=\gamma_{C_b}(x)$ lies in $Q_x$ and actually,
by Claim \ref {cl:caso2}, the point $x_b$ lies in the asymptotic cone 
$TC_x^{(4)}$.

The Gauss map $\gamma_\Xi:\Xi\map \PP^ {g-2}$ of $\Xi$ has an indeterminacy
point at $x$ and $-x$, which can be resolved by blowing up $x$ and $-x$,
since we may assume $\Xi$ to be symmetric. Let
$p: \tilde \Xi\to \Xi$ be the blow--up and let $\tilde \gamma_\Xi$ be the
morphism which coincides with $\gamma_\Xi\circ p: \tilde \Xi\map \PP^
{g-2}$ on an open subset. The exceptional divisor at $x$ and $-x$ is isomorphic
to $Q_x$ and $\tilde \gamma_\Xi(x_b)$ is the tangent hyperplane to
$Q_x$ at $x_b$. The tangency property of $\Xi$ and $\Xi_b$ along $C$
implies that the tangent hyperplane to $Q_x$ at $x_b$
coincides with $\gamma_{\Xi_b}(x)$.

Now we let $b$ vary in a component $Z$ of $N_1(B,\Xi)$ of dimension
$g-3$, so that we have a rational map
$$
f: Z\map Q_x, \quad b\mapsto x_b.
$$
Note that  $Q_x$ also has dimension $g-3$ and we claim that $f$ has
finite fibres, hence it is dominant. If not, we would
have an irreducible curve $\Gamma$ in $Z$ such that for all $b\in \Gamma$,
$x_b$ stays fixed. But then for the general $b\in \Gamma$, the divisor $\Xi_b$
has a fixed tangent hyperplane at $x$, a contradiction by
Lemma \ref {lem:fin}.

On the other hand, by Claim \ref {cl:caso2} and by part (i) of Proposition \ref
{prop:null}, one has that $f$ cannot be dominant, a
contradiction.\end{proof}

By Lemma \ref  {cl:caso3}, the curve
$\phi(R):=\Sigma$ is an irreducible component of the discriminant
$\Delta\subset \Qcal_R$ of singular quadrics in $\Qcal_R\cong \PP^ 2$.
Note that $\Delta$ has degree $g$ in $\PP^2$.
The map $\phi$ has degree at least $2$ since we may assume $\Theta_X$ to
be symmetric, hence it factors through the multiplication by $-1_{X}$
on $X$.

\begin{claim}\label {cl:deg2}
The map $\phi: R\to \Sigma$ has degree 2.
\end{claim}

\begin {proof} Let $d\geq 2$ be the
degree of the map. Then for $\xi=(X,\Theta_X,x)$ general in $R$, we have
distinct points $\xi_i=(X,\Theta_X,x_i)\in R$, with $\xi=\xi_1$, such that
all the quadrics $Q_{\xi_i}$, $i=1,\ldots,d$, coincide with the quadric
$Q=Q_\xi$. If for $i=1,\ldots,d$ we let $\eta_i$ be
the tangent vector to $S$ at $x_i$, we have that 
$\eta_1=\cdots =\eta_d$ and $\dd_{\eta_i}Q_{\xi_i}$,
for $i=1,\ldots,d$, coincide with the quadric $Q_1=\dd_{\eta_1} Q$, 
which is linearly independent from $Q$. 
The analysis we made in \S \ref {sec: essgi}
shows that $S_g$ is smooth at each of the points $\xi_i$, $i=1,\ldots,d$, and
the tangent space there is determined by $Q$ and $Q_{\eta_1}$.  This implies
that a general deformation of $\xi=(X,\Theta_X,x)$ inside $S_g$ carries
with it a deformation of each of the points $\xi_i$ ($i=1,\ldots,d$) in $S_g$,
because the involved quadrics are the same at these points.
This yields that a general element of an irreducible component of $N_{g,0}$
containing $N$ has at least $d$ singular points. By Debarre's result in
\cite {D1} (see \S \ref {sec:AMbound}), one has $d=2$, proving our claim. 
\end{proof}

\begin{claim}\label{cl:morf}
The map $\phi$ is a morphism.
\end{claim}

\begin{proof} To prove the claim, it suffices to show that $N$ is not contained
in $N_{g,0,r}$, with $r\geq 3$. In order to prove this, one verifies that,
for a general point $(\overline X,\overline \Xi)$ in a component $M$ of the 
intersection of $N$ with the boundary, there are no points of 
multiplicity $r\geq 3$ in ${\rm Sing}_{\rm vert}(\overline \Xi)$. 
Recall that, by Proposition \ref{rank1exists},
we may assume that $(\overline X,\overline \Xi)$ is a semi--abelian variety
of torus rank $1$, with abelian part $\eta=(B,\Xi)$. Moreover, 
by Claim \ref{cl:caso1}, only case (ii) of Lemma \ref{cor:struc} can occur. 
Therefore we may assume that $\eta$ is either a general point of 
$\theta_{0,g-1}$ or a general point of $M_{0,g-1}$ and the
extension corresponds to a general point in an irreducible component of 
$N_1(B,\Xi)$ which has dimension $g-3>0$.  
By Remark \ref{rem:triplepts} we see that no triple points can occur on 
${\rm Sing}_{\rm vert}(\overline \Xi)$. 
\end{proof}

Note now that the morphism $\phi$ is defined on $R\cong S$ by sections of
$\Ocal_S(\Theta_X)$, since the points of $S$ verify the equations \eqref
{eq:singt} and, if $\xi=(X,\Theta_X,x)\in R$, the entries of the matrix of
$Q_\xi$ are $\dd_i\dd_j\theta(\tau,z)$, where $z$ corresponds to $x$.
We deduce from $\deg(\Delta)=g$ and from Claim \ref{cl:deg2}, that
\begin{equation}\label{eq:ge}
S \cdot \Theta_X\leq 2g.
\end{equation}
As we assumed at the beginning of the proof, the class of  $S$ in $X$ is a
multiple $m\gamma_X$ of the minimal class. In view of \eqref {eq:ge}, we
find $m \leq 2$. The Matsusaka-Ran  criterion \cite {Ran} and a result of
Welters \cite{We2} imply that $(X,\Theta_X)$ is either a Jacobian or a Prym
variety. Since $g\geq 6$ this is not possible in view of the dimensions. 
This ends the proof.
\end{proof}

\begin{remark}\label{rem:verra} 
A.\ Verra communicated to us an interesting
example of an irreducible component $M$ of codimension $6$ of $N_{6,1}$
contained in the Prym locus. We briefly sketch, without entering in any
detail, its construction and properties.
Let $C$ be the normalization of
a general curve of type $(4,4)$ on $\PP^1\times \PP^1$ 
with two nodes on a line of type $(0,1)$,
so that $C$ has genus $g=7$. Let $d,t$ be the linear series formed
by pull-back divisors  on $C$
of the rulings of type $(1,0), (0,1)$ respectively. Consider a non trivial,
unramified double cover $f:\tilde C\to C$ of $C$ and let $(P,\Xi)$ be the
corresponding Prym variety. Then $\Xi$ has a $1$-dimensional unstable
singular locus $R$ (see \cite {LB}), 
homologically equivalent to twice the minimal
class, described by all classes in ${\rm Pic}^{12}(\tilde C)$ of
divisors of the form $f^*(d)+M$, with $f_*(M)\in t$. One proves that the
map $\phi$ described in the proof of Theorem \ref{th:con1} sends $R$ to a
plane sextic of genus $7$ which is tetragonally associated to $C$
(see \cite {LB} for the tetragonal construction). 
The divisor
$\Xi$ has $24$ further isolated singular points, which are pairwise exchanged
by the multiplication by $-1$.  One shows that the
corresponding $12$ tangent cones span the same linear system $\Qcal$ of
dimension $2$ spanned by $\phi(R)$.  The linear system $\Qcal$ is the
tangent space to the Prym locus $\Pcal_6$ at $(P,\Xi)$, which is therefore
a smooth point for $\Pcal_6$. By contrast, $M$ is a non-reduced component
of $N_{6,1}$ of codimension $6$ such that the projectivized normal space
$\Qcal$ at a general point has dimension $2$ rather than $5$. This shows
that the hypotheses of Theorem \ref{th:con1}  cannot be relaxed by
assuming only that the projectivized normal space to $M$ at a general point
has dimension $2$. 
\end{remark}

\end{section}

\begin{section}{Appendix: A Result on Pencils of Quadrics}\label{sec:segre}

One of the ingredients of the proof Theorem \ref {th:amlower} is a
classical result of Corrado Segre from \cite{Seg} on pencils of
quadrics. 

First we recall the following:

\begin{proposition}\label{prop:pencil}  Let $\Lcal$ be a 
pencil of quadrics in $\PP^n$ with $n\geq 1$
whose general member $Q$ is smooth. Then:

\begin{itemize}
\item[(i)] the number of singular quadrics $Q\in \Lcal$ is
$n+1$, where each such quadric $Q$ has to be counted with a
suitable multiplicity $\mu(Q)\geq n+1-{\rm rk}(Q)$; 
\item [(ii)] for a singular quadric $Q\in \Lcal$ one has $\mu(Q)\geq 2$ if
and only if  either ${\rm rk}(Q)<n$ or
the singular point of $Q$ is also a base point of
$\Lcal$;
\item [(iii)] for a singular quadric $Q\in \Lcal$ with rank $n$ 
one has $\mu(Q)= 2$ if
and only if any other quadric $Q'\in \Lcal$ is smooth at $p$ and the tangent
hyperplane to $Q'$ at $p$ is not tangent to $Q$ along a line.
\end{itemize}
\end{proposition}

\begin{proof} Consider the linear system $\Qcal_n$ of dimension $n(n+3)/2$
all quadrics in $\PP^ n$. Inside $\Qcal_n$ we have the discriminant locus
$\Delta_n$ of singular quadrics, which is a hypersurface of degree $n+1$,
defined by setting the determinant of a general quadric equal to zero.
The differentiation rule for determinants implies that the locus
$\Delta_{n,r}$ of quadrics of rank $r<n+1$ has multiplicity $n+1-r$ for
$\Delta_n$. By intersecting $\Delta_n$ with the line corresponding to
$\Lcal$ we have (i). 

As for assertion (ii), we may assume ${\rm rk}(Q)=n$, so that $Q$ has a unique
double point $p$, which we may suppose to be the point $(1,0,\dots,0)$. Thus
the matrix of $Q$ is of the form
$$
\left( \begin{matrix}0 &&& 0_n& \cr
0_n^ t &&& A& \cr\end{matrix}
\right)
$$
where $0_n\in \CC^ n$ is the zero vector and $A$ is a symmetric matrix of
order $n$ and maximal rank.  Let $Q'$ be another quadric in $\Lcal$, with matrix 
$$
\left( \begin{matrix}\beta &&& b& \cr
b^ t &&& B& \cr\end{matrix}
\right)
$$
with $\beta\in \CC$, $b=(b_1,\ldots,b_n)\in \CC^n$ and $B$ is a 
symmetrix matrix of order $n$. By intersecting $\Lcal$ with $\Delta_n$, 
we find the equation

\begin{equation}\label {eq:penceq}
\det\left( \begin{matrix}t\beta  &&& t b& \cr
tb^ t &&& A+tB& \cr\end{matrix}
\right)=0.
\end{equation}
The constant term in the left hand side is $0$. The coefficient of
the linear term is 
$$
\det\left( \begin{matrix}\beta  &&& b& \cr
0_n^ t &&& A& \cr\end{matrix}
\right)=\beta\det(A)
$$
which proves (ii). 

Let us prove (iii). Suppose ${\rm rk}(Q)=n$,
so that $Q$ has a unique double point $p$, which is a base point
of $\Lcal$. Again we may
suppose $p$ is the point $(1:0:\ldots:0)$ and we
can keep the above notation and continue the above analysis.
The left hand side in \eqref{eq:penceq} is
$$
t^ 2 \det\left(\begin{matrix}0&&&b&\cr
b^t&&&A+tB&\cr\end{matrix}\right)=0
$$
hence the coefficient of the third order term is
\begin{equation}\label{eq:term}
\det\left(\begin{matrix}0&&&b&\cr
b^ t&&&A&\cr\end{matrix}\right).
\end{equation}
One has$\mu(Q)=2$ if and only if this determinant is not zero, 
hence $b\not= 0$, which is equivalent
to saying that all quadrics in the pencil different from $Q$ are smooth at $p$.
Note that there is a vector $c=(c_1,\ldots,c_n)\in \CC^n$ such that
$b=c\cdot A$. Now the determinant in \eqref{eq:term} vanishes if and
only if $c\cdot b^t=0$, i.e.\ $c\cdot A\cdot c^ t=0$.
This means that the line $L$ joining $p$ with
$(0:c_1:\ldots:c_n)$ sits on $Q$ and that the tangent hyperplane to $Q'$
at $p$, which has equation $b_1x_1+\cdot+b_nx_n=0$,
is tangent to $Q$ along $r$.
\end{proof}

Next we prove Segre's theorem.

\begin{theorem}\label{thm:segre} Let $\Lcal$ be a linear
pencil of singular quadrics in $\PP^n$ with $n\geq 2$
whose general member $Q$ has rank $n+1-r$, i.e.\  ${\rm Vert}(Q)\cong
\PP^{r-1}$. We assume that ${\rm Vert}(Q)$ is not constant when $Q$ varies
in $\Lcal$ with rank $n+1-r$. Then:
\begin{itemize}
\item [(i)] the Zariski closure
$$
V_\Lcal = \overline {\big( \bigcup_{Q\in \Lcal, \, {\rm rk}(Q)= n+1-r} {\rm
Vert}(Q)\big)} 
$$
is a variety of dimension $r$ spanning a linear subspace $\Pi$ of dimension
$m$ in $\PP^n$ with $r\leq m \leq (n+r-1)/2$;
\item [(ii)] $V_\Lcal$ is a variety of minimal degree $m-r+1$ in $\Pi$;
\item [(iii)] if 
$$
\dim \bigl( { \bigcap_{Q\in \Lcal, \, {\rm rk}(Q)= n+1-r} 
{\rm Vert}(Q)} \bigr)=s
$$
then $r\leq (n+2s+3)/3$; 
\item[(iv)] the number of quadrics $Q\in \Lcal$ of rank ${\rm rk}(Q)<n+1-r$ is
$n+r-2m-1\leq n-r-1$, where each such quadric $Q$ has to be counted with a
suitable multiplicity $\nu(Q)\geq n+1-r-{\rm rk}(Q)$. 
\end{itemize}
\end{theorem}

\begin{proof} We start with the proof of part (i). Notice that, by
iteratedly restricting to a general hyperplane, we can reduce to the case $r=1$. 

In this case $V_\Lcal$ is a rational curve which, by Bertini's theorem is
contained in the base locus $B$ of $\Lcal$. Let $p,q$ be general points on it
and let $L$ be the line joining them.  There is a quadric $Q_p\in \Lcal$
with vertex at $p$. Hence $Q_p$ contains $L$. Similarly there is a
different quadric $Q_q$ with vertex at $q$, and it also contains $L$. Since
$Q_p$ and $Q_q$ span $\Lcal$, we see that $L$ is contained in $B$, i.e., the
secant variety to $V_\Lcal$ is contained in $B$. Take now three general
points $p,q,r$ on $V_\Lcal$. Since the lines $pq$, $pr$, $qr$ are contained
in $B$ also the plane spanned by $p,q,r$ is contained in $B$.
Continuing this way, we see that $\Pi=\langle V_\Lcal \rangle$ is
contained in $B$. Since the general quadric in $\Lcal$ has rank $n$
(recall we are assuming $r=1$ now), the
maximal dimension of subspaces on it is $n/2$. Thus $\dim(\Pi)\leq n/2$
which proves part (i).

Also for (ii) we can reduce ourselves to the case $r=1$, in
which we have to prove that  $V_\Lcal$ is a rational normal curve in 
$\Pi=\langle V_\Lcal \rangle$. Set $\dim(\Pi)=m$.

Let $p\in V_\Lcal$ be a general point. The polar hyperplane
$\pi_p$ of $p$ with respect to $Q\in \Lcal$ does not depend on $Q$, since
there is a quadric in $\Lcal$ which is singular at $p$ (see the proof of 
Proposition \ref {prop:pencil}). Note that $\pi_p$ has to contain all vertices
of the quadrics in $\Lcal$, hence it contains $\Pi$. By the linearity of
polarity, we have that polarity with respect to all quadrics in $\Lcal$ is
constant along $\Pi$ and for a general point $x\in \Pi$, the polar hyperplane
$\pi_x$ with respect to all quadrics in $\Lcal$ contains $\Pi$. Furthermore
the linear system of hyperplanes $\Pcal=\{\pi_x\}_{x\in \Pi}$ has dimension
$m-1$. 

Now, let $p\in V_\Lcal$ be a general point and let $Q_p$ be the unique quadric
in $\Lcal$ with a double point at $p$. We denote by ${\rm Star}(p)$ the $\PP^
{m-1}$ of all lines in $\Pi$ containing $p$. Let $\pi\in \Pcal$ be a general
hyperplane, which is tangent to $Q_p$ along a line $L$ containing $p$. Moreover
$L$ sits in $\Pi$, because this is the case if $\pi=\pi_q$ with $q$ another
general point on $V_\Lcal$, in which case $L$ is the line $\langle p, q
\rangle$. Thus we have a
linear map $\phi_p: \Pcal\to {\rm Star}(p)$, which is clearly injective and
therefore an isomorphism. 

Fix now another general point $q$ on $V_\Lcal$. The two maps $\phi_p$ and
$\phi_q$ determine a linear isomorphism $\phi: {\rm Star}(p)\to {\rm
Star}(q)$. Note that $L$ meets $\phi(L)$ if and only if $L\cap\phi(L)$ is a
point of $V_\Lcal$. This implies that $V_\Lcal$ is a rational normal curve in
$\Pi$, proving part (ii). 

Let us prove part (iii). It suffices to prove the assertion if
$s=-1$. The variety $V_\Lcal$ is swept out by a
$1$--dimensional family of projective spaces of dimension $r-1$, i.e.,
the vertices of the quadrics in $\Lcal$. Under the assumption $s=-1$ no
two of these vertices can intersect. Thus we must have $2(r-1)<m$. Using
part (i), the assertion follows. 

Finally we come to part (iv). Let us restrict $\Lcal$ to a general subspace
$\Lambda$ of dimension $n-r$. We get a pencil $\bar\Lcal$ of quadrics in
$\Lambda$ whose general member is smooth. 

We get a singular quadric in $\bar\Lcal$ when we intersect
$\Lambda$ with a quadric in $\Lcal$ whose vertex intersects $\Lambda$. We
claim that this is the only possibility for getting a singular quadric in
$\bar\Lcal$. Indeed, let $Q\in \Lcal$ and suppose that its
intersection $\bar Q\in \bar\Lcal$ with $\Lambda$ is singular at $p\in
\Lambda$, but $Q$ is not singular at $p$. Then $\Lambda$ is tangent to $Q$ at
$p$ and therefore also intersects the vertex of $Q$. 

In conclusion we have only two possibilities for getting singular quadrics
in $\bar \Lcal$:
\begin{itemize}
\item [(a)]  there is quadric of rank $n+1-r$ in $\Lcal$ whose vertex
intersects $\Lambda$;
\item [(b)]  there is quadric of rank $n+1-h<n+1-r$ in $\Lcal$ giving rise
to a quadric of rank $n-h$ in $\bar\Lcal$.
\end{itemize}
Case (a) occurs as many times as the degree of $V_\Lcal$, that is,
$m-r+1$ times. According to part (ii) of Proposition 
\ref{prop:pencil}, each quadric $Q$ in case (a) 
contributes with multiplicity at
least $2$ in the counting of singular quadrics in $\bar\Lcal$. We
claim that, because of the generality of $\Lambda$, this
multiplicity is exactly $2$. To prove this, by part (iii) of Proposition
\ref{prop:pencil}, we will prove that for each quadric $Q$ in case (a),
with vertex $p\in \Lambda$ and for any other quadric $Q'\in \Lcal$, the
tangent hyperplane $\pi$ to $Q'$ at $p$ is not tangent to $Q$ along a line
contained in $\Lambda$. To see this we can, by first cutting with a general
subspace of dimension $n-r+1$ through $\Lambda$, 
reduce ourselves to the case $r=1$, in which
$V_\Lcal$ is a rational normal curve.  Choose then a general point $q\in
V_\Lcal$ and let $Q'=Q_q$ be the unique quadric in $\Lcal$ with a double point at
$q$. The hyperplane $\pi$ is tangent to $Q_q$ along the line 
$\langle p,q \rangle$. This
implies that $\pi$ is tangent to $Q$ only along the tangent line $L_p$ to
$V_\Lcal$ at $p$ (see the proof of part ii)). 
By the generality assumption, $\Lambda$ is not
tangent to $V_\Lcal$ at $p$. Thus the assertion follows.   

As for quadrics in case (b), again by part (i) of
Proposition \ref {prop:pencil}, each such quadric contributes
to the same count with multiplicity $h-r$. Since, by part (i) of Proposition \ref
{prop:pencil},  the number of singular quadrics in $\bar\Lcal$,
counted with appropriate multiplicity, is $n-r+1$, the assertion follows.
\end{proof}

One has the following consequence:

\begin{corollary}\label{cor:segre} Let $\Lcal$ be a linear
pencil of quadrics in $\PP^n$ with $n\geq 2$. Then the 
general member $Q\in \Lcal$ has rank $n+1-r$ (i.e., ${\rm Vert}(Q)\cong
\PP^{r-1}$) if and only if  the base locus of $\Lcal$ contains a linear
subspace $\Pi$ of dimension $m$ with $r\leq m \leq (n+r-1)/2$, along which
all the quadrics in $\Lcal$ have a common tangent subspace of dimension
$n+r-m-1$. In this case $\Pi$ is the span of the variety $V_\Lcal$.
\end{corollary}

\begin{proof} 
As usual it suffices to prove the assertion for $r=1$. 
If the general quadric $Q\in \Lcal$ has rank ${\rm rk}(Q)=n$, 
the assertion follws from the proof of Theorem \ref{thm:segre}. 
The converse is trivial, since a smooth quadric in $\PP^n$ has a 
tangent subspace of dimension $n-m-1$, and
not larger, along a subspace of dimension $m$. 
\end{proof}

These results imply the existence of {\it canonical forms} for pencils of
singular quadrics, originally due to Weierstrass \cite{Wei} and Kronecker 
\cite{K}. 
This is explained in some detail in \cite{Seg}, \S\S 20-25, and we will not
dwell on this here.

It would be desirable to
have an extension of the results in this Appendix  
to higher dimensional linear families of quadrics.
\end{section}


\begin{thebibliography}{99}

\bibitem{AN} V.\ Alexeev, I.\ Nakamura: On Mumford's construction of
degenerating abelian varieties. 
{\sl Tohoku Math.\ J.\ \bf 51}  (1999),  399--420.

\bibitem {AM} A.\ Andreotti, A.L.\ Mayer: 
On period relations for abelian integrals on algebraic curves. 
{\sl Ann.\ Sc.\ Norm.\ Pisa \bf 3} (1967), p.\ 189--238.

\bibitem {ACGH}  A.\ Arbarello, M.\ Cornalba, Ph.\ Griffiths, J.\ Harris:
 Geometry of Algebraic Curves, vol. I, 
{\sl Grundlehren der math. Wissenschaft, Springer Verlag, \bf 267}, 1985.

\bibitem{B2} A.\ Beauville:
Prym varieties and the Schottky problem.
{\sl Invent.\ Math.\  \bf 41} (1977), p.\ 149--196.

\bibitem{Br} M.\ Brion: 
Compactification de l'espace des modules des vari\'et\'es 
ab\'eliennes principalement polaris\'ees [d'apr\`es V. Alexeev] 
S\'eminaire Bourbaki 2005.

\bibitem{Ca-M} S.\ Casalaina-Martin: 
Prym varieties and the Schottky problem for cubic threefolds.
ArXiv: {\tt math.AG/0605666}.

\bibitem{cvdg} C.\ Ciliberto, G.\ van der Geer: Subvarieties of the moduli
space of curves parametrizing jacobians with non trivial endomorphisms.
{\sl Amer.\ J.\ Math., \bf 114} (1991), 551--570.

\bibitem{C-vdG} C.\ Ciliberto, G.\ van der Geer:
The moduli space of Abelian varieties and the singularities of 
the theta divisor. In: Yau, S.T. (ed.), 
Surveys in differential geometry. 
Papers dedicated to Atiyah, Bott, Hirzebruch and Singer. 
Somerville, MA: International Press. Surv.\ Differ.\ Geom., 
Suppl.\ J.\ Differ.\ Geom.\ 7, 61-81 (2000).

\bibitem{De1} O.\ Debarre: 
$\Scal_g$ est Cohen--Macaulay.
Private communication, unpublished.

\bibitem{D0} O.\ Debarre: Sur les vari\'et\'es ab\'eliennes dont
le diviseur th\^eta est singulier en codimension $3$.
{\sl Duke Math.\ J.,\bf 25} (1988), 221--273.

\bibitem{D1} O.\ Debarre: 
Le lieu des vari\'et\'es ab\'eliennes dont le diviseur th\^eta est 
singulier a deux composantes.
{\sl Ann.\ Sci.\ \'Ecole Norm.\ Sup., \bf 25} (1992), p.\  687--707.

\bibitem{EL} L.\ Ein, R.\ Lazarsfeld:
Singularities of theta divisors, and the birational geometry of irregular
varieties. 
{\sl J.\ Am.\ Math.\ Soc.\ \bf 10}, (1997), p.\ 243-258. 

\bibitem{FC} G.\ Faltings, C. L.\ Chai:
Degeneration of abelian varieties. 
{\sl Ergebnisse der Math., 3 Folge}
{\bf 22}, Berlin, Springer--Verlag, 1990. 

\bibitem{vdG1} G.\ van der Geer: Cycles on the moduli space of abelian
varieties. In: Moduli of Curves and Abelian Varieties. C.\ Faber, E.\ Looijenga, Ed., Aspects of Math.\ 33, (1999) Vieweg.

\bibitem{GH} Ph.\ Griffiths, J.\ Harris: 
Principles of algebraic geometry,  A.\ Wiley, 1978.

\bibitem{GS} S.\ Grushevski, R.\ Salvati Manni: 
Jacobians with a vanishing theta--null in genus $4$. ArXiv:
{\tt arXiv:math. AG/0505160}.

\bibitem{H} R. \ Hartshorne: Algebraic Geometry. 
{\sl Graduate texts in Math., \bf 52}, 1977.

\bibitem{Kempf} G.\ Kempf:
A Torelli theorem for osculating cones to the theta divisor.
{\sl International Centre for Theoretical Physics, Trieste;
Course on Riemann Surfaces}, 1987.

\bibitem{KS} G.\ Kempf, F.\ O.\ Schreier:
A Torelli theorem for osculating cones to the theta divisor. 
{\sl Compositio Math., \bf 67} (1988), 343--353.

\bibitem{Ko} J.\ Koll\'ar:
Shafarevich maps and automorphic forms. 
{\sl Princeton Univ. Press}, (1995).

\bibitem{K} L.\ Kronecker:
Ueber Schaaren von quadratischer Formen. {\sl Monatsber. Akad.
Wiss. Berlin}, (1868), 339--346. 

\bibitem{Ku} E.\ Kunz: 
Holomorphe Differentialformen auf algebraischen Variet\"aten 
mit Singularit\"aten. I.
{\sl Manuscr.\ Math.\ \bf 15} (1975), 91--108.

\bibitem{LB} H.\ Lange, Ch. Birkenhake:
Complex Abelian Varieties. {\sl Grundlehren der math.\ Wiss.}, 302,
Springer--Verlag, 1992.

\bibitem{Li} J.\ Lipman:
Dualizing sheaves, differentials and residues on algebraic varieties. 
{\sl Ast\'erisque} 117, (1984).

\bibitem{MumfordAV} D.\ Mumford: Abelian Varieties.
Studies in Mathematics, 5. Tata Institute of Fundamental Research, Bombay. 
Oxford etc.: Oxford University Press. XII, 279 p. 

\bibitem{Mum} D.\ Mumford: An analytic construction of degenerating
abelian varieties. {\sl Comp. Math.}, {\bf 24} (1972), 239--272.

\bibitem{Mumford} D.\ Mumford: On the Kodaira dimension of the Siegel
modular variety. In: Algebraic geometry--open problems (Ravello, 1982),  
348--375, Lecture Notes in Math., 997, Springer, Berlin, 1983. 

\bibitem{Nakam} I.\ Nakamura: On moduli of stable quasi-abelian varieties.
{\sl Nagoya Math.\ J.\ \bf 58} (1975), 149--214. 

\bibitem{Nam} Y.\ Namikawa: A new compactification of the Siegel space
and degeneration of abelian varieties. {\sl Math. Ann.}, I, II, {\bf 221}
(1976), 97--141, 201--241.

\bibitem{PP} G. Pareschi, M. Popa:  Generic vanishing and minimal 
cohomology classes on abelian varieties. ArXiv: {\tt math. AG/0610166}.

\bibitem{Ran} Z.\ Ran: On subvarieties of abelian varieties.
{\sl Inventiones Math., {\bf 62}}, (1981), 459--479.

\bibitem{Sat} I.\ Satake: On the compactification of the Siegel space. {\sl
J. Indian math. Soc.}, {\bf 20} (1956), 259--281.

\bibitem{Seg} C.\ Segre: Ricerche sui fasci di coni quadrici in uno
spazio lineare qualunque. {\sl Atti della R.\ Accademia delle
Scienze di Torino  \bf XIX}, (1883/4), p.\ 692-710.

\bibitem{serre} J.\ P.\ Serre: Algebraic Groups and Class Fields. {\sl
Graduate texts in Math  \bf 117}, Springer--Verlag, 1988. 



\bibitem{SV1} R.\ Smith, R.\ Varley:
Components of the locus of singular theta divisor of genus $5$. {\sl
Springer  Lecture Notes in Math., Proceedings of the Conference in Algebraic
Geometry, Sitges, 1983}, {\bf 1124}, (1985), 338--416.

\bibitem{SV} R.\ Smith, R.\ Varley:
Multiplicity $g$ points on theta divisors. 
{\sl  Duke Math.\ J.\  \bf 82}  (1996), p.\ 319--326.

\bibitem{Wei} K.\ Weierstrass:
Zur Theorie der bilineaeren un quadratischen Formen. {\sl Monatsber. Akad.
Wiss. Berlin}, (1868), 310--338.

\bibitem{We1} G.\ Welters:
Polarized abelian varieties and the Heat Equation. {\sl Compositio Math.\ \bf
49}, (1983), no. 2, 173--194.

\bibitem {We2} G.\ Welters: Curves of twice the minimal class on principally
polarized abelian varieties. {\sl Indag. Math., {\bf 94}}, (1987),
87--109.

\bibitem {zak} F.\ Zak: Determinants of projective varieties., in Algebraic
transformation groups and algebraic varieties, V.\ L.\ Popov (ed.),
Encyclopedia of Math.\  Invariant theory and transformation groups III, 
Springer Verlag, 207--238 (2004).


\end{thebibliography}
 \end{document}